\documentclass[12pt]{amsart}
\usepackage{amssymb}
\usepackage{latexsym}
\usepackage{graphicx}
\usepackage{subfigure}
\usepackage{tikz}
\usepackage[linktocpage=true]{hyperref}
\usepackage{ulem}
\usepackage[left=3cm, right=3cm, top=3cm, bottom=3cm]{geometry}
\usepackage{color}

\newcommand{\re}[1]{(\ref{#1})}

\newcommand{\rl}[1]{Lemma~\ref{#1}}

\newcommand{\rp}[1]{Proposition~\ref{#1}}
\newcommand{\rt}[1]{Theorem~\ref{#1}}
\newcommand{\rd}[1]{Definition~\ref{#1}}
\newcommand{\rrem}[1]{Remark~\ref{#1}}

\hypersetup{ colorlinks   = true, 
urlcolor  = blue, 
linkcolor    = blue, 
citecolor   = blue 
}

\def\I{{\mathrm I}}

\def\SS{{\mathbb S}}
\def\om{{\omega}}

\def\lam{{\lambda}}
\def\del{{\delta}}

\def\al{{\alpha}}
\def\Ga{{\Gamma}}
\def\De{{\Delta}}
\def\ti{\tilde}
\def\ka{{\kappa}}

\def\beq{\begin{equation}}
\def\eeq{\end{equation}}

\def\Exp{\mathrm{Exp}}
\def\IM{\mathrm{Im}}

\def\Ker{\mathrm{Ker}}

\def\det{\mathrm{det}\ }

\def\id{\mathrm{Id}}
\def\bm{\begin{matrix}}
\def\em{\end{matrix}}

\newcommand{\Z}{{\mathbb Z}}
\newcommand{\R}{{\mathbb R}}

\newcommand{\C}{{\mathbb C}}
\newcommand{\D}{{\mathbb D}}
\newcommand{\T}{{\mathbb T}}

\newcommand{\N}{{\mathbb N}}
\newcommand{\PP}{{\mathbb P}}

\newcommand{\CD}{{\mathcal D}}
\newcommand{\CE}{{\mathcal E}}

\newcommand{\CG}{{\mathcal G}}

\newcommand{\CH}{{\mathcal H}}
\newcommand{\CI}{{\mathcal I}}

\newcommand{\CW}{{\mathcal W}}
\newcommand{\CP}{{\mathcal P}}
\newcommand{\CS}{{\mathcal S}}

\newcommand{\CZ}{{\mathcal Z}}
\newcommand{\CU}{{\mathcal U}}

\newcommand{\CR}{{\mathcal R}}

\newtheorem{thm}{Theorem}[section]
\newtheorem{cor}[thm]{Corollary}
\newtheorem{prop}[thm]{Proposition}
\newtheorem{lemma}[thm]{Lemma}
\newtheorem{remark}[thm]{Remark}

\newtheorem{defi}[thm]{Definition}
\newtheorem{hypo}[thm]{H}

\newcommand{\la}{\langle}
\newcommand{\ra}{\rangle}

\definecolor{gainsboro}{rgb}{0.86, 0.86, 0.86}
 \definecolor{gray-light}{rgb}{0.75, 0.75, 0.75}
  	\definecolor{gray}{rgb}{0.5, 0.5, 0.5}

\newcommand{\RE}{\operatorname{Re}}

\begin{document}
\setcounter{thm}{0}\setcounter{equation}{0}

\title[]{Local rigidity of group actions of isometries on compact Riemannian manifolds}

\author{Laurent Stolovitch}
\address{Universit\'e C\^ote d'Azur, CNRS, Laboratoire J. A. Dieudonn\'{e}, 06108 Nice, France}
\email{Laurent.stolovitch@univ-cotedazur.fr}
\thanks{}

\author{Zhiyan Zhao}
\address{Universit\'e C\^ote d'Azur, CNRS, Laboratoire J. A. Dieudonn\'{e}, 06108 Nice, France}
\email{zhiyan.zhao@univ-cotedazur.fr}
\thanks{}

\begin{abstract}
In this article, we consider perturbations of isometries on a compact Riemannian manifold $M$.
We investigate the smooth (resp. analytic) rigidity phenomenon of groups of these isometries. As a particular case, we prove that if a finite family of smooth (resp. analytic) small enough perturbations is simultaneously conjugate to the family of isometries via a finitely smooth diffeomorphism, then it is simultaneously smoothly (resp. analytically) conjugate to it whenever the family of isometries satisfies a Diophantine condition.
Our results generalize the rigidity theorems of Arnold, Herman, Yoccoz, Moser, etc. about circle diffeomorphisms which are small perturbations of rotations as well as Fisher-Margulis's theorem on group actions satisfying Kazhdan's property (T).
\end{abstract}

\maketitle

\tableofcontents

\section{Introduction and main results}

The aim of this article is to study the smooth\footnote{Through this paper, ``smooth" is interpreted as in the $C^\infty$ category unless otherwise specified.} and analytic rigidity of a group action by isometries on a compact Riemannian manifold $M$ (supposed to be connected and without boundary). Both the manifold $M$ and its Riemannian metric $g$ are assumed to be smooth or analytic. Let $G$ be a finitely presented group acting on $M$ via smooth or analytic isometries $\pi$. We also consider a group action $\pi_0$ of $G$ on $M$ by diffeomorphisms, which represents a small smooth or analytic perturbation of $\pi$. Our goal is to provide conditions under which $\pi_0$ is smoothly or analytically conjugate to $\pi$.

This problem originates from the seminal works of Arnold \cite{arnold-cercle}, Herman \cite{Herman-ihes} and Yoccoz \cite{yoccoz-cercle1998} on circle diffeomorphisms. It was demonstrated that if a diffeomorphism $F$ is a smooth or analytic small enough perturbation of a rotation $R_\alpha$ by a {\it Diophantine} angle $\alpha$, and if the {\it rotation number} of $F$ is $\alpha$, then $F$ is smoothly or analytically conjugate to $R_\alpha$. A similar statement was established in the smooth category by Moser \cite{Moser-cercle} for abelian groups of smooth circle diffeomorphisms. We refer to these results as ``local" rigidity theorems, meaning that we only consider small perturbations of rotations.
Global results (i.e., without the smallness assumption on the perturbation) were achieved by Herman \cite{Herman-ihes}, Yoccoz \cite{yoccoz-cercle}, and Fayad-Khanin \cite{FK-cercle} for single circle diffeomorphisms and abelian groups of such diffeomorphisms, respectively. These rigidity problems have a long history. For recent developments in slightly different contexts, see \cite{deW23, DK-duke, FK-ICM18, Khanin-ICM18}.
This problem, although not directly related, also echoes questions arising in Zimmer's program (e.g., \cite{fisher-zimmer}), as well as those in hyperbolic dynamics (e.g. \cite{DWX-duke, katok-hyperbolic}).

In the present work, we consider a general compact manifold $M$ of arbitrary dimension, which plays the role of the circle, on which a group of isometries acts. This group is defined by a finite number of generators and relations, with the isometries playing the role of rotations. We introduce an appropriate notion of ``Diophantiness" for the isometries, which heavily depends on the geometry and metric of the manifold $M$, particularly involving the spectrum of the Laplace-Beltrami operator.
 The condition analogous to requiring that ``the rotation number of the perturbation equals that of the rotation being perturbed" can be expressed as ``the perturbation can be conjugated as closely as desired to the unperturbed isometries".
In other contexts, this might be phrased as ``the perturbation is {\it almost conjugate} to the unperturbed one" or `` the unperturbed isometries are {\it almost rigid} ".
The point is then to prove that one can effectively achieve a genuine smooth or analytic conjugacy between the unperturbed and the sufficiently small perturbed actions.

One corollary of our main results can be stated in a non-technical way as follows.

\begin{thm}\label{thmimpressive}
Let $M$ be a smooth (resp. analytic) compact Riemannian manifold of dimension $n$ (connected and without boundary). Let finitely many smooth (resp. analytic) isometries $\pi$ on $M$ satisfy a simultaneous Diophantine condition. Then there exists $R>0$ such that any sufficiently small smooth (resp. analytic) perturbations of these isometries into diffeomorphisms of $M$, which are simultaneously conjugate to the original isometries through a $C^0$ near-identity $C^R$ transformation, are smoothly (resp. analytically) simultaneously conjugate to the original isometries.
\end{thm}

\noindent In other words, finitely smooth local rigidity implies smooth or analytic local rigidity:
	\begin{center}
Diophantine + $C^R{\rm -rigid} \ \Longrightarrow \ C^{\infty}{\rm -rigid} \ / \ C^{\omega}{\rm -rigid}.$
	\end{center}

\subsection{Diophantine properties of group action by isometries}\label{sec_Dio-intro}

In what follows, $M$ denotes a smooth (resp. analytic) compact Riemannian manifold of dimension $n$ (connected and without boundary).

Given a finitely presented group $G$, with the set of generators $\CS=\{\gamma_1,\cdots,\gamma_k\}$ satisfying a finite collection $\CR$ of relations, let $\pi:G\to {\rm Isom}^\infty(M)$ be a $G-$action by smooth isometries of the smooth compact Riemannian manifold $M$. 
Let $L^2(M,TM)$ be the $L^2$ vector field on the tangent bundle $TM$. Define the linear operator $d_0: L^2(M,TM)\to L^2(M,TM)^k$ as
$$ d_0 v =(v-\pi(\gamma_l)_*v)_{1\leq l\leq k},\qquad v\in L^2(M,TM). $$
where $\pi(\gamma_l)_*$ is the push-forward introduced by the isometry $\pi(\gamma_l)$.
The adjoint of $d_0$ through the $L^2-$scalar product on $M$, $d_0^*: L^2(M,TM)^k\to L^2(M,TM)$, is defined as
$$ d_0^*V=\sum_{l=1}^k(v_l-\pi(\gamma_l^{-1})_*v_l),\qquad V=(v_1,\cdots, v_k)\in L^2(M,TM)^k.$$
Hence we obtain the self-adjoin non-negative operator $d_0\circ d_0^*$ on $L^2(M,TM)^k$:
$$(d_0\circ d_0^*)V=\left(\sum_{l=1}^k(v_l-\pi(\gamma_l^{-1})_*v_l-\pi(\gamma_j)_*v_l+\pi(\gamma_j\gamma_l^{-1})_*v_l)\right)_{1\leq j\leq k}.$$

For the Laplace-Beltrami operator on the tangent bundle $\Delta_{TM}$, let $\{\lambda_j\}_{j\in\N}$ be the eigenvalues of $|\Delta_{TM}|^{\frac12}$, strictly increasing w.r.t. $j$ and tending to $\infty$ as $j\to\infty$. Then we have the orthogonal decompositions
\begin{equation}\label{decompL2-intro}
L^2(M,TM)=\bigoplus_{j\geq 0} E_{\lam_j},\qquad L^2(M,TM)^k=\bigoplus_{j\geq 0} E_{\lam_j}^k,  \
\end{equation}
where $E_{\lam_j}$ is the eigenspace associated to the eigenvalue $\lambda_j$.
It is known that each $E_{\lam_j}$ is finite-dimensional (see Section \ref{secVB}) and any isometry leaves invariant  $E_{\lambda_j}$\footnote{These decompositions are deduced from Peter-Weyl theorem, with details stated in Section \ref{IsoTMsmooth}.}.

By relating $d_0$ to the spectral properties of the Laplace-Beltrami operator on the tangent bundle $\Delta_{TM}$, we define the Diophantine condition for the generating smooth isometries as follows.

\begin{defi}\label{defi-dio-gnrt}
The $G-$action $\pi$ by isometries on $M$ is said to be $\bf d_0-${\bf Diophantine} if there exist $\sigma>0$ and $\tau\geq 0$ such that
\begin{equation}\label{diophantine-gnrt-intro}
\|(d_0\circ d_0^*)u\|_{L^2}\geq \frac{\sigma}{(1+\lambda_j)^\tau}\|u\|_{L^2},\qquad u\in \IM d_0\cap E_{\lam_j}^k.
\end{equation}
The $G-$action $\pi$ satisfying (\ref{diophantine-gnrt-intro}) is also called ${\bf (\sigma,\tau)-}{\bf d_0-}${\bf Diophantine}.
\end{defi}

\begin{remark}\label{rmk-d0Dio}
	The $(\sigma,\tau)-d_0-$Diophantine condition is weaker than the simultaneous Diophantine condition for the generating isometries $(\pi(\gamma))_{\gamma\in \CS}$, corresponding to Dolgopyat's definition (see Definition \ref{dioph-dolg} and Lemma \ref{lemma_d_0-dio} in Section \ref{sec_gener_iso}). Indeed, in the later, it is implicitly assumed that $\Ker d_0\subset E_{0}$ (in particular $\dim \Ker d_0<+\infty$) and inequality akin to \re{diophantine-gnrt-intro} are supposed to hold on the full space, while the lower bound (\ref{diophantine-gnrt-intro}) is required to hold only on $\IM d_0$ and the vectors in $\Ker d_0^*=(\IM d_0)^\perp$ are not involved.
	
	The Diophantine property in Dolgopyat sense (Definition \ref{dioph-dolg}) coincides with that for rotations on the circle ((1.3) in \cite{Moser-cercle}), and for translations on the torus  ((2.3) in \cite{petko-tore}).
	Moreover, this property is widely used for group actions by isometries when the group is a discrete group with Kazhdan's property (T) or an irreducible lattice in a semi-simple Lie group with rank at least $2$
	(see Section \ref{sec_groups}).
\end{remark}

Following the concepts outlined in Section \ref{sec_GroupAction}, we observe that an action of a finitely presented group $G$ by smooth (or analytic) isometries on a smooth (or an analytic) manifold $M$ induces an action on the tangent bundle  $TM$. This action gives rise to the Hochschild complex of cochains of $L^2$ vector fields:
\begin{equation}\label{complex}
L^2(M,TM)\stackrel{d_0}{\longrightarrow}C^1(G, L^2(M,TM))\stackrel{d_1}{\longrightarrow}C^2(G, L^2(M,TM))\longrightarrow \cdots .
\end{equation}
The spaces $C^1(G, L^2(M,TM))$ and $C^2(G, L^2(M,TM))$ in the above complex can be identified with $L^2(M,TM)^k$ and $L^2(M,TM)^p$ respectively, where $k$ denotes the cardinal of $\CS$ (the set of generators) and $p$ that of $\CR$ (the set of relations; see Section \ref{sec_GroupAction} for more details).
The operators $d_0$ and $d_1$ have explicit expressions in terms of the actions of the generators and their relations.
We then introduce the self-adjoint Box operator
$$\square=d_0\circ d_0^*+ d_1^*\circ d_1: L^2(M,TM)^k\to L^2(M,TM)^k,$$
which is fundamental for our purpose,
with the adjoint being defined upon the $L^2-$scalar product on $M$.
Recalling the orthogonal decompositions in (\ref{decompL2-intro}), it will be shown in Section \ref{sec_GroupAction} that $\square$ is invariant on every $E_{\lam_j}^k$.

Relating the spectral properties of $\square$ to that of the Laplace-Beltrami operator on the tangent bundle $\Delta_{TM}$, we define the Diophantine condition for the action by smooth isometries as follows.
\begin{defi}\label{defi-dio-intro}
The $G-$action $\pi$ by isometries on $M$ is said to be {$\bf\square-$}{\bf Diophantine} if there exist $\sigma>0$ and $\tau\geq 0$ such that
\begin{equation}\label{dio-box}
\|\square V\|_{L^2}\geq \frac{\sigma}{(1+\lambda_j)^\tau}\|V\|_{L^2} ,\quad \forall \ V\in \IM \square\cap E_{\lam_j}^k.
\end{equation}
The $G-$action $\pi$ satisfying (\ref{dio-box}) is also called $\bf (\sigma,\tau)-$$\bf\square-${\bf Diophantine}.
\end{defi}

\begin{remark}
We emphasize again that, as in Definition \ref{defi-dio-gnrt}, the lower bound (\ref{dio-box}) holds only for $u \in \operatorname{Im} \square$, and the vectors in $\operatorname{Ker} \square$ are not involved.
The $\square$-Diophantine condition does not imply that the subset of generators of $G$ is a Diophantine subset (in the sense of Dolgopyat, see Definition \ref{dioph-dolg}) for $\pi$. Indeed, Proposition \ref{prop_box-not-dolg} provides an example of $G$-action satisfying $\square$-Diophantine condition while its set of generators is not a Diophantine set.
\end{remark}

\begin{remark}\label{rmk_DioConds}
For the complex (\ref{complex}), since $d_1\circ d_0=0$, we have
$\IM \square=\IM d_0\bigoplus\IM d^*_1$.
The Diophantine condition in Definition \ref{defi-dio-gnrt} provides an asymptotic lower bound, polynomially decaying w.r.t. $\lambda_j$, for the non-vanishing eigenvalues of $d_0\circ d^*_0$ (as shown in Lemma \ref{lemma_d_0-dio}), whereas the condition in Definition \ref{defi-dio-intro} ensures such a lower bound for the non-vanishing eigenvalues of both $d_0\circ d^*_0$ and $d^*_1\circ d_1$. Hence, the $\square-$Diophantine condition of $\pi$ means that it has Diophantine generating isometries and Diophantine relations among them.
\end{remark}

\begin{defi}\label{defi-dio-relation}
The $G-$action $\pi$ is called to have ${\bf (\sigma,\tau)-}${\bf Diophantine relations} if there exist $\sigma>0$ and $\tau\geq 0$ such that
\begin{equation}\label{dio-box}
\|(d_1^*\circ d_1) V\|_{L^2}\geq \frac{\sigma}{(1+\lambda_j)^\tau}\|V\|_{L^2} ,\quad \forall \ V\in \IM d_1^*\cap E_{\lam_j}^k.
\end{equation}
\end{defi}

\begin{remark}\label{rmk_Diobox-equi}
As an equivalent definition, $\pi$ is $(\sigma,\tau)-\square-$Diophantine if and only if it is $(\sigma,\tau)-d_0-$Diophantine and has $(\sigma,\tau)-$Diophantine relations.
\end{remark}

\subsection{Almost conjugacy and local rigidity}

Let us denote by $\Ga^{k}=\Ga^{k}(M,TM)$, $k\in\N\cup\{\infty\}$, (resp. $\Ga^{\omega}=\Ga^{\omega}(M,TM)$) the space of $C^k$ (resp. analytic) sections of $TM$ (resp. provided that $M$ is real analytic).
Let $\Exp$ be the exponential map defined upon the Riemannian connection (see Section \ref{secExp} for more details). Given $P\in\Ga^{k}$, which is $C^0$ sufficiently small, we denote by $\pi_P$ the $G-$action by $C^k$ smooth diffeomorphisms on $M$, satisfying
\begin{equation}\label{expPpi}
\pi_P(\gamma)=\Exp\{P(\gamma)\}\circ \pi(\gamma),\quad \gamma\in \CS,\quad for \ some \ P:\CS\to \Gamma^k(M,TM).
\end{equation}
The almost conjugacy of the $G-$action by diffeomorphisms is defined as follows.

\begin{defi}\label{defi_almost_conj-intro}
Given $0<\zeta<1$, $R\in\N^*$, $\pi_P$ is said to be {$\bf \zeta -$}$\bf C^R$ {\bf almost conjugate} to $\pi$ on $M$, if, for any $\varepsilon>0$, there exists $y^{\varepsilon} \in \Ga^R(M,TM)$ with
\beq\label{esti-y} \|y^{\varepsilon}\|_{C^0}<\zeta, \qquad  \|y^{\varepsilon}\|_{C^R}<\zeta^{-1},\eeq
and $z^{\varepsilon}:\CS\to \Ga^{0}(M,TM)$ with
\beq\label{esti-z}\|z^{\varepsilon}\|_{\CS,C^0}:=\left(\sum_{1\leq l\leq k}\|z_{m}^{\varepsilon}(\gamma_l)\|^2_{ C^0}\right)^{\frac12}<\varepsilon,\eeq
such that, for $\gamma\in \CS$,
\begin{equation}\label{eq_form_conj}
\Exp\{y^\varepsilon\}^{-1} \circ \Exp\{P(\gamma)\}\circ \pi(\gamma) \circ  \Exp\{y^\varepsilon\}=\Exp\{z^\varepsilon(\gamma)\} \circ \pi(\gamma).
\end{equation}
\end{defi}

 \begin{remark}\label{formalC0} If the action $\pi_P$ is $\zeta-C^R$ almost conjugate to $\pi$ on $M$ with $R\geq 3$, then it is topologically conjugate to $\pi$, since, through interpolation (Lemma \ref{interpol-cs}), (\ref{esti-y}) implies that $ \|y^{\varepsilon}\|_{C^1}<\zeta^\frac13$.
 Indeed, given a sequence $\varepsilon_m$ tending to $0$, we have a sequence $\{y_m\}_m$ bounded by $\zeta^\frac13$ in the $C^{1}$-topology and a sequence $\left\{\{z_m(\gamma)\}_{\gamma \in \CS}\right\}_m$, bounded in the $C^0$-topology by $\varepsilon_m$, (hence tending to 0), such that
$$\Exp\{y_m\}^{-1} \circ \pi_P(\gamma)\circ  \Exp\{y_m\} =\Exp\{z_m(\gamma)\} \circ \pi(\gamma).$$ By the Arzel\`a-Ascoli Theorem, we can extract a subsequence $\{y_{m_k}\}_{k}$ converging in the $C^0$-topology to a Lipschitz continuous vector field $y$. Hence, as in \cite{Mos69}, this defines a homeomorphism $\Exp\{y\}$ of $M$ such that
$\Exp\{y\}^{-1} \circ \pi_P(\gamma)\circ  \Exp\{y\}=\pi(\gamma)$.

In the specific case where $M=\T^1$ (the circle) and $\pi$ is a rotation $R_\alpha$,
if $\pi_P$ is almost conjugate to $R_\alpha$ in the above sense, then it must have the same rotation number $\alpha$ as $R_\alpha$, since the rotation number is an invariant under topological conjugacy for circle diffeomorphisms.
\end{remark}

We can state our first main result~:
\begin{thm}\label{thmmain}
Let $\pi$ be a $(\sigma,\tau)-d_0-$Diophantine $G-$action by smooth (resp. analytic) isometries on $M$.
There exists $\widehat R >0$ (depending on $n$ and $\tau$) such that any $G-$action $\pi_{P_0}$ by sufficiently small smooth (resp. analytic) perturbations, which is $\|P_0\|_{\CS,C^0}^\frac34 - C^{\widehat R}$ almost conjugate to $\pi$, is smoothly (resp. analytically) conjugate to $\pi$, i.e., there exists $W\in \Ga^{\infty}$ (resp. $\Ga^{\omega}$) such that, for every
$\gamma\in G$,
$\Exp\{W\}^{-1} \circ \pi_{P_0}(\gamma) \circ \Exp\{W\} = \pi(\gamma) $.
\end{thm}

In the special case of groups satisfying Kazhdan's property (T), Theorem \ref{thmmain} covers Fisher-Margulis's theorem \cite{fisher-margulis-invent} as well as its (new) analytic counterpart 
(see Corollary \ref{corKazhdan}).

\begin{remark} The conclusion of Theorem \ref{thmmain} naturally also holds for one single free diffeomorphism on $M$.
\end{remark}


Recalling Remark \ref{rmk-d0Dio}, and since $C^{\widehat R}$ conjugacy implies $C^{\widehat R}$ almost conjugacy, we obtain Theorem \ref{thmimpressive} as corollary of Theorem \ref{thmmain}.


\subsection{First cohomology and local rigidity}

Recalling the complex \re{complex}, its first cohomology group is defined as
$$H^{1}(G,L^2(M,TM)):=\Ker\, d_1/\IM\, d_0.$$
Our second main statement concerns the case where the perturbation is a priori not assumed to be almost conjugate to $\pi$~:
\begin{thm}\label{thm-geom0}
Let $\pi$ be a $\square-$Diophantine $G-$action by smooth (resp. analytic) isometries on $M$. 
Assume that $H^{1}(G,L^2(M,TM))=0$.
Then any $G-$action by smooth (resp. analytic) diffeomorphisms on $M$ which is sufficiently close to $\pi$ is smoothly (resp. analytically) conjugate to $\pi$.
\end{thm}

The smooth version of the above theorem can be seen as Fisher's local rigidity result \cite{Fisher}[Theorem 1.1] under the extra {\it Diophantine} assumption on the action by isometries.  Lemma \ref{lemma-cyclic} provides a simple example of vanishing cohomology.

\subsection{Idea of proof}

The main purpose is to establish a smooth or analytic conjugacy between the unperturbed group action by isometries $\pi$ and its perturbation $\pi_{P_0}$.

\subsubsection{KAM scheme for smooth rigidity}
In the smooth context, as a natural strategy, we proceed by setting up an iterative Kolmogorov-Arnold-Moser scheme (see Section \ref{secKAMsmooth}).
Regarding rigidity, there are two main challenges:
\begin{itemize}
  \item Solving and estimating a solution to the {\it cohomological equation} in order to build up a ``controlled" conjugacy that transforms the perturbed action into a perturbation much closer to the unperturbed action,
  \item Proving that {\it obstructions} - terms that cannot be eliminated by solving the cohomological equations - are in fact much smaller than expected.
\end{itemize}

To estimate solutions of the cohomological equation, we need the Diophantine conditions in Definition \ref{defi-dio-gnrt} (for Theorem \ref{thmmain}) or in Definition \ref{defi-dio-intro} (for Theorem \ref{thm-geom0}), and make use of the spectral properties of the Laplace-Beltrami operator on the tangent bundle $\Delta_{TM}$.
Both Diophantine conditions ensure that the non-vanishing eigenvalues of the operator $d_0\circ d_0^*$, associated with the group action $\pi$ by isometries, do not accumulate too rapidly zero, providing us with polynomial lower bounds with respect to the eigenvalues of $|\Delta_{TM}|^\frac12$.
At the $m-$step of the iteration, we have a perturbation $\pi_{P_m}$ of $\pi$ (which is conjugate to the original perturbation $\pi_{P_0}$). We then proceed a {\it KAM step}, meaning that we solve a linearized conjugacy equation - the {\it cohomological equation} - in order to build a conjugacy of $\pi_{P_m}$ into some $\pi_{P_{m+1}}$ which would be expected to be much closer to $\pi$ than $\pi_{P_m}$.
The Diophantine conditions guarantee that the solutions to these equations exist and are well-behaved. This is formulated in Section \ref{sec_cohomoeq-smooth}.

However, in general, there are {\it obstructions}, that cannot be handled by solving a cohomological equation as they do not lie in the range of the linearized operator. Our other assumptions will allow to obtain a much closer $\pi_{P_{m+1}}$.
For instance, in the case of close-to-rotation circle diffeomorphisms, the {\it obstruction} is merely a constant Fourier mode (which lies in a finite-dimensional space). Assuming that the rotation number of the circle diffeomorphism is the same as the rotation it is a perturbation of, ensures that this obstruction is in fact much smaller and even more, it is of the size of the expected $P_{m+1}$. Hence, it does not require to be taken care of at the $m-$KAM step.
As for a group action $\pi$ by isometries on a general compact Riemannian manifold, the obstructions are included in $\Ker(d_0\circ d_0^*)$, which is generally of infinite dimension. If the group action $\pi_{P_m}$ is {\it almost conjugate to $\pi$}, then these obstructions are also much smaller than expected (in case of Theorem \ref{thmmain}). In order to iterate this scheme, we need to prove that the new $\pi_{P_{m+1}}$ obtained through $m$-KAM-steps is still {\it almost conjugate to $\pi$}. This is presented in detail in Section \ref{sec_refine-smooth} and \ref{sec_verify_alconj}.

If the first cohomology group $H^1(G,L^2(M,TM))$ vanishes, then we obtain the same conclusion as in the previous case as soon as group action $\pi$ by isometries has Diophantine relations (as in Theorem \ref{thm-geom0}).

The smooth rigidity is then established by showing the composition of the $m$ first conjugacies obtained at first $m$ KAM steps, converges in $C^{\infty}$ topology.

\subsubsection{Hardy interpolation and analytic rigidity}

In the analytic context, we consider the complex neighborhoods $M_{r}$, with $r>0$ sufficiently small, of the analytic compact Riemannian manifold $M$, known as {\it Grauert tubes} \cite{Gra58} (see Section \ref{sec_Grauert}). The complex structure on the tangent bundle of Grauert tubes was studied in a series of works by Sz\"oke \cite{Szo91, Szo95}, Lempert-Sz\"oke \cite{LeSz91}, Guillemin-Stenzel \cite{GS91, GS92}.
Additionally, Boutet de Monvel \cite{BdM78} (see also \cite{Leb18, Zel12}) defined them as domains for the holomorphic extension of the eigenvectors of the Laplace-Beltrami operator on $M$.

For the group action $\pi_{P_0}$ by analytic diffeomorphisms on $M$ as in (\ref{expPpi}), the vector field $P_0(\gamma)$ extends to a holomorphic vector field on some Grauert tube $M_{r_0}$. Through the work of Boutet de Monvel \cite{BdM78}, we consider vector fields belonging to the {\it Hardy space of $M_{r_0}$} (see Definition \ref{def-Hardy}).
These are holomorphic vector fields over $M_{r_0}$ with $L^2$ boundary values, characterized by the decay of the coefficients of its ``Fourier-like decomposition" and depending on the ``radius" $r_0$ of the Grauert tube.

Instead of presenting a parallel KAM scheme\footnote{In previous work \cite{SZ2023}, a KAM scheme on analytic perturbation of $\pi$ is given under the supplementary hypothesis $\dim\Ker\square<\infty$ and $\pi$ is $\square-$Diophantine. Through this KAM scheme,
the analytic rigidity is shown for any analytic perturbation which is $C^1$ almost conjugate to $\pi$.} as in the smooth context, we take advantage of the fact that the sequence of group actions from the smooth KAM scheme remains $\{\pi_{P_m}\}$ analytic.
Through Hardy interpolation (see Lemma \ref{interpol-hardy}), it is shown in Section \ref{sec_seq-anal} that there exists a subsequence $\{\pi_{\widehat P_m}\}\subset \{\pi_{P_m}\}$ with the analytic $\widehat P_m$ extending to the Grauert tube $M_{r_m}$ and converging to zero on $M_{\frac{r_0}2}$, with the sequence of radii $\{r_m\}$ satisfying $\frac{r_0}2<r_{m}< r_0$.
Then the analytic rigidity is established by verifying the convergence of the sequence of approximate conjugacies on $M_{\frac{r_0}2}$.

\subsection{Description of the remaining of paper} The remaining of paper is organized as follows.
Section \ref{section-geom} is devoted to the smooth geometric setting, including definitions concerning isometries and properties of smooth norms and Sobolev norms. In Section \ref{sec_GroupAction} we introduce the group action $\pi$ by isometries, as well the associated operators $d_0$, $d_1$ and $\square$.
We also discuss the group action by diffeomorphisms that is a perturbation of $\pi$. Several examples of smooth and analytic rigidity, as applications of the main theorems, are provided in Section \ref{sec_exemple}.
In Section \ref{secKAMsmooth}, we present the KAM scheme in the smooth context, which proves the smooth rigidity. Section \ref{sec_Grauert_Hardy} is devoted to the analytic geometric setting, including definitions of Grauert tubes and Hardy spaces, as well as the equivalence between different norms. The proof of analytic rigidity is then provided in Section \ref{secKAManal}.

In Appendices \ref{app_proof-s1smooth} and \ref{app_proof}, we prove Proposition \ref{prop_s1-smooth} and Proposition \ref{propMoser-anal}, respectively. In Appendix \ref{app_interpolation}, we recall various interpolation inequalities that are used in the proofs.

\medskip

\noindent{\it Acknowledgment.} This work was stimulated by a work of David Fisher \cite{Fisher} with whom the first author had discussions around 2007 about it. Although not published, this article contains lot of interesting examples, in the smooth category. The first author thanks Charlie Epstein for having pointed out Boutel de Monvel's theory of holomorphic extension of eigenvectors, L\'{a}zl\'{o} Lempert, Matthew Stenzel and Robert Sz\"{o}ke for exchanges about Grauert tubes and also David Fisher, Bassam Fayad and Jonhattan DeWitt for exchanges on group actions.

\section{Vector field on a smooth compact Riemannian manifold}\label{section-geom}

As a preliminary step, we introduce in this section some basic notions and crucial results about vector fields on a smooth compact Riemannian manifold.

\smallskip

Let $M$ be a compact $C^\infty$ Riemannian manifold of dimension $n\geq 1$. The Riemannian metric, which is assumed to be smooth, is defined by a scalar product $\langle\cdot,\cdot\rangle_m$ on the tangent space $T_mM$ for every $m\in M$. In local coordinates $(x_1,\cdots,x_n)$ over which the tangent bundle is trivialized, we write
$$\langle v,w\rangle_m=\sum_{1\leq i,j\leq n}g_{i,j}(x(m))v^i w^j,\qquad {\rm for} \  \  v=\sum_{i=1}^nv^i\frac{\partial}{\partial x_{i}}, \quad  w=\sum_{i=1}^nw^i\frac{\partial}{\partial x_{i}},$$
with the matrix $g(m)=(g_{i,j}(x(m)))_{1\leq i,j\leq n}$ positive definite at every $m\in M$, if $(m,v)$ and $(m,w)$ belong to $T_mM$.
Define the isomorphism $\al_m : T_mM\rightarrow T_m^*M$ by
$$\al_m(v)w:=\langle v,w\rangle_m,\qquad {\rm for} \  \  v,w\in T_mM.$$
This induces a scalar product on the cotangent bundle $T^*_mM$~:
$$\langle v^*,w^*\rangle_m:=\langle v,w\rangle_m, \qquad {\rm for} \  \ v^*=\al_m(v),\ w^*=\al_m(w).$$
The above is extended to an isomorphism
$\al_m: \wedge^pT_mM\rightarrow \wedge^pT^*_mM$ by
$$\alpha_m(v_1\wedge \cdots\wedge v_p):=\alpha_m(v_1)\wedge \cdots\wedge \alpha_m(v_p). $$
The scalar product induced on $\wedge^pT_mM$ is defined to be
$$\langle v_1\wedge \cdots\wedge v_p,w_1\wedge \cdots\wedge w_p\rangle_m:=\det (\langle v_i,w_j\rangle_m)_{1\leq i,j\leq p}.$$

Let us denote by $\Ga^{k}=\Ga^{k}(M,TM)$, $k\in\N\cup\{\infty\}$, (resp. $\Ga^{\omega}=\Ga^{\omega}(M,TM)$) the space of $C^k-$smooth (resp. analytic, provided that $M$ is analytic) sections of $TM$.
Let $\CU=\{(U_i, h_i,\phi_i)\}_i$ be an atlas of $M$ over which $TM$ is trivialized: $\{U_i\}_i$ is a covering of $M$ by finite open sets, $h_i:U_i\rightarrow \R^n$ is a $C^\infty$ diffeomorphism and $\{\phi_i\}_i$ is a partition of unity subordinated to the covering. For $v\in \Gamma^\infty(M,TM)$, its $C^R-$norm, $R\in\N$, is defined by
 \begin{equation}\label{normsCk}
\|v\|_{C^R}:=\max_{0\leq l\leq R} |v|_l \quad {\rm with} \quad |u|_l:=\sum_{i}|(h_i)_*(\phi_i v)|_l,
\end{equation}
where the latter is the sup-norm over $\R^n$ of the $l^{\rm th}-$order derivative of $\R^n$-valued functions. Moreover, for $V=(v_1,\cdots,v_\nu)\in\Ga^{R}(M,TM)^\nu$, we define
 \begin{equation}\label{normsCk-vector}
 \|V\|_{C^R}:=\left(\sum_{j=1}^\nu \|v_j\|^2_{C^R}\right)^\frac12. \end{equation}

\subsection{Exponential map of smooth vector field}\label{secExp}

For a compact  $C^\infty$ Riemannian manifold  $M$, let us recall the {\it exponential map} $\Exp\{\cdot\}$ defined upon the Riemannian connection (see \cite{Hel01}[Chap. I, Section 6], \cite{Mos69}[Section 2-a)]).

Consider a fixed point  $m \in M$  and a tangent vector  $\xi \in T_m M$ at  $m$.
The exponential map $\Exp_m$ is defined such that $\Exp_m\{\xi\}$ is the point on the geodesic starting at  $m$  in the direction $\xi$, at a distance  $|\xi|_{g(m)}$, the length of  $\xi$  in the given metric.
It is well known that, for sufficiently small $r > 0$, the exponential map $\Exp_m$ transforms the set $\{ \xi \in T_m M : |\xi|_{g(m)} < r \}$  diffeomorphically onto a neighborhood of  $m$  in  $M$.
For a sufficiently small continuous vector field  $v$ on $M$, the map $\psi: m \mapsto \Exp\{v(m)\}$ represents a continuous mapping of  $M$  into itself. Conversely, every continuous mapping  $\psi(\cdot)$ sufficiently close to the identity can be represented in the form $\Exp\{v(\cdot)\}$ for a vector field $v$, and $v$ has the same smoothness property as $\psi$.

For every $m\in M$, let $B_m(0,r)\subset T_mM$ be the ball in the tangent space $T_mM$ centered at $0$ with radius $r$. There exists $r(m)>0$ such that the mapping
$$\Exp_m : B_m(0,r(m))\subset T_mM \rightarrow M$$
is a $C^\infty$ diffeomorphism onto its image and it is analytic if $M$ is analytic. Moreover, the mapping $m\mapsto r(m)$ can be chosen lower semi-continuous.

About the composition of exponential map, we have the following lemma.
\begin{lemma}\label{lem_moser-s1}(Moser \cite{Mos69}[Lemma 1])
Let $v\in \Gamma^0$, $w\in\Gamma^1$ and $\|v\|_{C^0}$, $\|w\|_{C^1}$ sufficiently small. Then there exists $s_1(w,v)\in \Gamma^0$ such that,
\begin{equation}\label{eq_lemMoser1-smooth}
\Exp \{w\}\circ \Exp \{v\} = \Exp \{w+v+s_1(w,v)\},
\end{equation}
with $s_1(w,0)=s_1(0,v)=0$ and $\|s_1(w,v)-s_1(w,v')\|_{C^0}\leq c_0\|w\|_{C^1}\|v-v'\|_{C^0}$,
where $c_0$ is a constant which depends only on the manifold and the metric.
\end{lemma}

We have the generalized $C^R$ estimate for $s_1(\cdot,\cdot)$:

\begin{prop}\label{prop_s1-smooth} For $w_1,w_2\in \Gamma^1$ and $v\in \Gamma^0$ with $\|w_1\|_{C^1}$, $\|w_2\|_{C^0}$ and $\|v\|_{C^0}$ sufficiently small, we have $s_1 (w_1+w_2,v)\in\Gamma^0$ with
\begin{equation}\label{s1-C0-decomp}
\|s_1 (w_1+w_2,v)\|_{C^0}\lesssim \|w_1\|_{C^1}\|v\|_{C^0}+\|w_2\|_{C^0}.
\end{equation}
For $w, v\in \Gamma^R $ with $\|w\|_{C^1}$ and $\|v\|_{C^1}$ sufficiently small,
we have $s_1 (w,v)\in\Gamma^R$ with
\begin{equation}\label{s1-Cr}
 \|s_1(w,v)\|_{C^R}\lesssim_R \|w\|_{C^R}+\|w\|_{C^1}\|v\|_{C^R},
\end{equation}
and, if $w\in\Gamma^{R+1} $, we have
\begin{equation}\label{s1-Cr+}
    \|s_1(w,v)\|_{C^R} \lesssim_R \|w\|_{C^2}\|v\|_{C^R} + \|w\|_{C^{R+1}}  \|v\|_{C^0}.
\end{equation}
\end{prop}

\smallskip

In the inequalities (\ref{s1-C0-decomp}) -- (\ref{s1-Cr+}) and afterward in the smooth setting, ``$\lesssim$" means boundedness from above by a positive constant depending on the manifold $(M,g)$,
and ``$\lesssim_R$" means that the implicit constant depends also on the differential order $R$.

\begin{remark}
The a priori non symmetric estimate \re{s1-C0-decomp} is due to non symmetric assumptions. The particular situation with $w_2=0$ is stated in Lemma \ref{lem_moser-s1}.
\end{remark}

\begin{remark}\label{rmk_vect_s1}
The vector field $s_1(w,v)$ in (\ref{eq_lemMoser1-smooth}) can be presented in a more general form.
For $w\in \Gamma^1$, $v\in \Gamma^0$ and $W=(w_j)_{1\leq j\leq \nu}\in (\Gamma^1)^\nu$, $V=(v_j)_{1\leq j\leq \nu}\in (\Gamma^0)^\nu$, let us define
$$ s_1 (W,v):=(s_1 (w_j,v))_{j},\quad
 s_1 (w,V):=(s_1 (w,v_j))_{j},\quad
 s_1 (W,V):=(s_1 (w_j, v_j))_{j},$$
all of which belong to $(\Gamma^0)^\nu$. With $\|W\|_{C^R}$ and $\|V\|_{C^R}$ defined in (\ref{normsCk-vector}), the estimates (\ref{s1-C0-decomp}) -- (\ref{s1-Cr+}) hold for $s_1 (W,v)$, $s_1 (w,V)$ and $s_1 (W,V)$ if the hypotheses in Proposition \ref{prop_s1-smooth} are satisfied.
\end{remark}

The proof of Proposition \ref{prop_s1-smooth} is postponed in Appendix \ref{app_proof-s1smooth}.


\begin{lemma}\label{lemma_exp-inverse} 
Let  $R\in\N^*$. For $w\in\Ga^R$ with a sufficiently small $\|w\|_{C^1}$, there exists $\tilde w\in \Gamma^R$ such that $\Exp\{w\}^{-1}=\Exp\{\tilde w\}$ with $\|\tilde w\|_{C^R}\lesssim_R \|w\|_{C^R}$.
\end{lemma}
\proof Let us find $\tilde w\in \Gamma^R$ solving the equation
$$\Exp\{w\}\circ\Exp\{\tilde w\}={\rm Id.}=\Exp\{0\}.$$
Through Lemma \ref{lem_moser-s1}, it is sufficient to find a small enough $C^0$-norm $\tilde w$ such that $w+\tilde w+s_1(w,\tilde w)=0$. According to (\ref{s1-Cr}), we have 
	$\|s_1(w,\tilde w)\|_{C^R}\leq g_R( \|w\|_{C^R}+ \|w\|_{C^1}  \|\tilde w\|_{C^R})$. If $\|w\|_{C^1}g_R\ll1$, then the Implicit function theorem yields the existence of $\tilde w\in \Gamma^R$, such that
	$\|\tilde w\|_{C^R}\leq (1-g_R \|w\|_{C^1} )^{-1}  g_R \|w\|_{C^R}$.\qed

\subsection{Spaces of sections of vector bundles}\label{secVB}
Let $\rm E$ be a smooth vector bundle over the smooth compact Riemannian manifold $M$.
We denote by $\Ga^{\infty}(M,{\rm E})$ be the space of $C^\infty-$smooth sections of ${\rm E}$.
If ${\rm E}$ admits a smooth scalar product $\langle\cdot,\cdot\rangle_{{\rm E}}$, then we define the scalar product on the space of section to be
$$
\langle v,w\rangle:=\int_M\langle v(x),w(x)\rangle_{{\rm E},x} \,  d{\rm vol}(x),\quad  v,w\in \Ga^{\infty}(M,{\rm E}),
$$
where $d{\rm vol}$ is a volume element which can be expressed in local coordinates~:
$$
d{\rm vol}(x)= \sqrt{\det g_{i,j}(x)} \,  dx_1\cdots dx_n.
$$
Let $L^2(M,{\rm E})$ denote the completion of $\Ga^{\infty}(M,{\rm E})$ with respect to this scalar product. It is the Hilbert space of $L^2$ sections of ${\rm E}$.

From the de Rham complex, we construct a complex on the space of smooth sections of multi-vector fields as follows:
$$
\begin{array}{cccccc}
	\Ga^{\infty}(M,\Bbb R) & \stackrel{d_0}{\longrightarrow}& \Ga^{\infty}(M,T^*M) & \stackrel{d_1}{\longrightarrow}& \Ga^{\infty}(M,\wedge^2 T^*M) & \stackrel{d_2}{\longrightarrow}\cdots\\
	& & & & & \\
	\parallel & &\downarrow\alpha^{-1} & & \downarrow\alpha^{-1}& \\
		& & & & & \\
	\Ga^{\infty}(M,\Bbb R) & \stackrel{\ti d_0}{\longrightarrow}& \Ga^{\infty}(M,TM) & \stackrel{\ti d_1}{\longrightarrow}& \Ga^{\infty}(M,\wedge^2 TM) & \stackrel{\ti d_2}{\longrightarrow}\cdots\\
\end{array}  \  \  \  \  .
$$
The first differentials are defined to be $\ti d_0 := \al^{-1}\circ d_0$ and $\ti d_1 := \al^{-1}\circ d_1\circ\alpha $. We shall call the $L^2$ extension of this complex the
``tangential complex" of $M$.
Since the de Rham complex is elliptic (the complex of the associated symbols is exact), the same holds true for the tangential complex. We then define the Laplacian on the tangent bundle to be the self-adjoint operator
$
\De_{TM} := \ti d_0\circ \ti d_0^*+ \ti d_1^*\circ \ti d_1$.
According to the classical generalized Hodge theory for elliptic complex \cite{Wel80}[Chapter IV, Theorem 5.2], \cite{Dem96}[Theor\`eme 3.10, Corrolaire 3.16] or \cite{Ros97}, there exists an orthonormal basis $({\bf e}_i)_{i\geq 0}$ of eigenvectors of $-\De_{TM}$ in $L^2(M, TM)$, with the associated eigenvalues $(\tilde\lam_i^2)_{i\geq 0}$ satisfying\begin{equation}\label{tildelambda}
\tilde\lambda_0=0\leq\tilde\lam_1\leq\tilde\lam_2 \leq \cdots, \quad \lim_{i\to\infty}\tilde\lam_i=+\infty.\end{equation}
 Each eigenvalue is of finite multiplicity and $+\infty$ is the only accumulation point. Moreover, according to Weyl asymptotic estimate, we have
\begin{equation}\label{asymp-lamdak}
\#\left\{i\in\N:\tilde\lambda^2_i \leq \lambda\right\}\sim a_0 \lambda^{\frac{n}{2}},\quad \tilde\lambda^2_i\sim b_0 \cdot i^{\frac2n}  \quad {\rm as}  \   \ i\to \infty
\end{equation}
for some constants $a_0, b_0>0$ depending only on the Riemannian manifold $(M,g)$.
(see e.g., \cite{B86}[Page 70]).
The eigenvectors ${\bf e}_i$ are smooth on $M$
and real analytic if $M$ is real analytic \cite{treves-book-analytic}[Theorem 4.1.2].

\subsection{Action of isometry on $TM$}\label{IsoTMsmooth}

By a smooth {\it isometry} of $M$, we mean a $C^\infty$ diffeomorphism of $M$ which preserves the distance induced by the Riemannian metric $g$.
Recall that the group of smooth isometries $\text{Isom}^\infty(M)\subset {\rm Diff}^\infty(M)$ of $M$ is compact \cite{Kob72}[Chap. 2, Theorem 1.2].
As diffeomorphisms, elements of $\text{Isom}^\infty(M)$ act on the space of sections $L^2(M,TM)$ by push-forward. For a.e. $x\in M$, the action is given by
$$f_*v(x)=Df(f^{-1}(x))v(f^{-1}(x)),\quad  f\in \text{Isom}^\infty(M),\quad v\in L^2(M,TM).$$
This action is {\it unitary}~:
\beq\label{scal-tg}
\langle f_*v,f_*w\rangle=\int_M \langle f_*v,f_*w\rangle_{g,x}d{\rm vol}(x)=\int_M \langle v,w\rangle_{f^{-1}_*g,f^{-1}(x)} d{\rm vol}(f^{-1}(x))=\langle v,w\rangle,
\eeq
where $\langle \cdot, \cdot\rangle_{g,x}$ denotes the inner product at $x$ induced by the metric $g$, and $f^{-1}_*g$ is the pullback metric under $f^{-1}$.
This action of the group of isometries commutes with the Laplacian on $TM$, that is
\beq \label{iso-laplacian}
\Delta_{TM}(f_*v)=f_*\De_{TM} v.
\eeq
This follows from the general fact that $f_*\De_{TM,g}(v)=\De_{TM,f_*g}(f_*v)$.

According to the Peter–Weyl theorem \cite{Zim90}[p.61], the Hilbert space $L^2(M,TM)$ can be decomposed into an orthogonal sum of finite-dimensional subspaces $V_i$ which are irreductible with respect to the action of $\text{Isom}^\infty(M)$~:
\beq\label{peter-weyl}
L^2(M,TM)=\bigoplus_{i\geq 0} V_i.
\eeq
In particular, each $V_i$ is invariant under the action, that is, $f_*V_i\subset V_i$, for all $f\in \text{Isom}^\infty(M)$ and all indices $i$. Moreover, each $V_i$ is contained in an eigenspace $E_{\lam_j}$ associated to the eigenvalue $\lam_j$ of $|\De_{TM}|^\frac12$, which is also the eigenspace associated to the eigenvalue $-\lam_j^2$ of $\De_{TM}$. Hence, we have
\beq\label{E_lamj}
E_{\lam_j}=\bigoplus_{i\in  J_j} V_i \quad {\rm with} \quad J_j:=\left\{i\in \N: V_{i}\subset E_{\lambda_j} \right\},\eeq
which gives rise to the first decomposition in (\ref{decompL2-intro}). Recalling the orthonormal basis $({\bf e}_i)_{i\geq 0}$ of eigenvectors of $|\De_{TM}|^\frac12$ (as well as $\De_{TM}$) in $L^2(M, TM)$, we have the generalized ``Fourier expansion" on the smooth Riemannian manifold $M$.

\begin{thm}\label{theo-boutet-L2}\cite{BdM78}
Given a global $L^2$ section $u\in L^2(M,TM)$, we have
\beq\label{generalFourier}
u=\sum_{i\geq 0} \hat u_i{\bf e}_i=\sum_{j\in\N}\sum_{i\in I_j} \hat u_{i} {\bf e}_i \quad with \quad  \hat u_i\in \R,\eeq
 where, for $j\in\N$, $I_j:=\{i\in \N:\tilde \lambda_i=\lambda_j\}$.
\end{thm}
\noindent Based on the decomposition (\ref{generalFourier}), let $\PP_j$ be the projection onto $E_{\lam_j}$ or $E_{\lam_j}^\nu$, $\nu\in\N^*$, (depending on the context), i.e.,
\begin{equation}\label{Proj_i}
\PP_j u=\left\{\begin{array}{cl}
  \displaystyle \sum_{i\in I_j} \hat u_i {\bf e}_i, &   \displaystyle  u= \sum_{i\in \N} \hat u_i {\bf e}_i\in L^2(M,TM)\\[2mm]
  \displaystyle \left(\sum_{i\in I_j} \hat u_{l,i} {\bf e}_i\right)_{1\leq l\leq \nu}, &  \displaystyle  u= \left(\sum_{i\in \N}  \hat u_{l,i} {\bf e}_i\right)_{1\leq l\leq \nu}\in L^2(M,TM)^\nu
\end{array}\right.  .\end{equation}

\smallskip

Isometries act on the exponential map of smooth vector field as follow (see e.g \cite{oneill-book}[p.91])
\begin{lemma}\label{lem_exp-isometry-smooth}
For $\pi\in{\rm Isom}^\infty(M)$ and $w\in \Ga^0$ with $\|w\|_{C^0}$ sufficiently small, we have
$$
\pi\circ \Exp\{w\} \circ \pi^{-1}=\Exp \left\{(D\pi\cdot w)\circ\pi^{-1}\right\}= \Exp \left\{\pi_* w\right\}.
$$
\end{lemma}

\subsection{Sobolev norm}\label{sec-Hr_space}
For $R\geq 0$, based on the decomposition (\ref{generalFourier}), we define the {\bf $\bf \CH^R-$space} or {\bf Sobolev space} as
$$\CH^R :=\left\{u=\sum_{j\in\N}\sum_{i\in I_j} \hat u_{i} {\bf e}_i\in L^2(M,TM): \sum_{j\in\N}(1+\lambda_j)^{2R}\sum_{i\in I_j} |\hat u_{i}|^2 <\infty\right\}.$$
According to \cite{gilkey}, we have $\Ga^\infty\subset \CH^R $ for any $R\geq0$.
 We define the $\bf \CH^R-${\bf norm} or {\bf Sobolev norm} as
  $$\|u\|_{\CH^R}:= \left(\sum_{i\in\N}(1+\tilde \lambda_i)^{2R} |\hat u_{i}|^2 \right)^{\frac12}=\left(\sum_{j\in\N}(1+\lambda_j)^{2R}\sum_{i\in I_j} |\hat u_{i}|^2 \right)^{\frac12}.$$
More generally, for $u=(u_l)_{1\leq l\leq \nu}\in (\CH^R)^\nu$, we define
$\displaystyle \|u\|_{\CH^R}:=\left(\sum_{1\leq l\leq \nu}\|u_{l}\|_{\CH^R}^2 \right)^{\frac12}$.

\begin{prop}\label{propNorm-CrHr} For $R\in\N$, we have $\CH^{R+\frac{3}{2}n+1} \subset \Ga^R\subset \CH^R$ and
\beq\label{propnormr}
\|u\|_{\CH^R} \lesssim_R \|u\|_{C^R} \lesssim_R \|u\|_{\CH^{R+\frac{3}{2}n+1}},\quad \forall \ u\in \CH^{R+\frac{3}{2}n+1} .\eeq
\end{prop}
\proof Since $|\De_{TM}|^\frac12$ is a self-adjoint elliptic partial differential operator of order $1$
on the smooth compact Riemannian manifold $M$ without boundary of dimension $n$, according to \cite{gilkey}[Lemma 1.6.3 (b)], for $R\in\N$, there exists $l_R>0$ such that
$$\|{\bf e}_i\|_{C^R}\lesssim_R (1+\tilde\lambda_i)^{l_R}.$$
For $u\in \CH^{l_R+n+\frac12}$, we have
$$|\hat u_{i}|\leq \|u\|_{\CH^{l_R+n+\frac12}}(1+\tilde\lambda_i)^{-(l_R+n+\frac12)},$$
hence, recalling that $\tilde\lambda_i\sim i^{\frac{1}{n}}$,
$$\|u\|_{C^R} \lesssim_R \sum_{i\in\N}(1+\tilde\lambda_i)^{l_R} |\hat u_{i}|\lesssim_R \|u\|_{\CH^{l_R+n+\frac12}}  \sum_{i\in\N}(1+\tilde\lambda_i)^{-(n+\frac12)}\lesssim_R  \|u\|_{\CH^{l_R+n+\frac12}}.$$
More precisely, according to \cite{gilkey}[Proof of Lemma 1.6.3 (b)], we can choose $l_R=R+\frac{n}{2}+\frac12$.
On the other hand,
$$\|u\|_{\CH^R}= \left(\sum_{j\in\N}(1+\lambda_j)^{2R}\sum_{i\in I_j} |\hat u_{i}|^2 \right)^{\frac12} \lesssim_R \||\Delta|^{\frac{R}2} u\|_{\CH^0}\lesssim_R \|u\|_{C^R}.\qed $$

Recall the projection defined in (\ref{Proj_i}), we have
\begin{cor}\label{cor-CR-truncation}
Given any subspace $E$ of $L^2(M,TM)^\nu$, $\nu\in \N^*$, let ${\bf P}_E$ be the orthogonal projection from $L^2(T,TM)^\nu$ onto $E$.
If ${\bf P}_E\circ\PP_j=\PP_j\circ{\bf P}_E$ for every $j\in\N$, then, for any $N\in\N$, for any $u\in (\Ga^R)^\nu$,
$$\left\|\sum_{\lam_j\leq N}( {\bf P}_E\circ\PP_j) u \right\|_{C^R}, \ \left\|\sum_{\lam_j> N}( {\bf P}_E\circ\PP_j) u \right\|_{C^R}\lesssim_R N^{\frac32n+1}\| u \|_{C^R}.$$
\end{cor}
\proof By computations with (\ref{propnormr}), we have, for $u\in (\Ga^R)^\nu$,
\begin{eqnarray*}
\left\|\sum_{\lam_j\leq N}( {\bf P}_E\circ\PP_j) u \right\|_{C^R}&\lesssim_R& \left\|\sum_{\lam_j\leq N}( {\bf P}_E\circ\PP_j) u \right\|_{\CH^{R+\frac32n+1}} 
\lesssim_R N^{\frac32n+1} \left\|\sum_{\lam_j\leq N}(\PP_j\circ{\bf P}_E) u \right\|_{\CH^{R}}\\
&\lesssim_R& N^{\frac32n+1} \left\|u\right\|_{\CH^{R}} \ \lesssim_R \  N^{\frac32n+1} \left\|u\right\|_{C^{R}}.
\end{eqnarray*}
On the other hand,
\begin{eqnarray*}
\left\|\sum_{\lam_j> N}( {\bf P}_E\circ\PP_j) u \right\|_{C^R}\leq \|u\|_{C^R} + \left\|\sum_{\lam_j\leq N}( {\bf P}_E\circ\PP_j) u \right\|_{C^R}\lesssim_R  N^{\frac32n+1} \left\|u\right\|_{C^{R}}.\qed
\end{eqnarray*}

\section{Group action by isometries on $M$}\label{sec_GroupAction}

For a {\it finitely presented group} $G$, let us fix its presentation, i.e.,
\begin{itemize}
\item a finite collection $\CS=\{\gamma_1, \cdots, \gamma_k\}$ of generators,
\item a finite collection of relations $\CR=\{\CW_1, \cdots, \CW_p\}$, where each $\CW_i$ is a finite word of generators in $\CS$ and their inverses.
\end{itemize}
In other words, with $\gamma_{l+k}:=\gamma_{l}^{-1}$, $1\leq l\leq k$, we can view each $\CW_j$ as a word over the alphabet of the $2k$ letters $\{\gamma_l\}_{1\leq l\leq 2k}$:
$$\CW_j =\gamma_{l^{(j)}_1}\cdots \gamma_{l^{(j)}_{m_j}},\quad 1\leq l^{(j)}_1, \cdots , l^{(j)}_{m_j}\leq k, \quad j=1,\cdots,p.$$
Furthermore, we impose the relations $\CW_j=e$ in $G$, where $e$ is the identity element, for each $1\leq j\leq p$.
See Section \ref{sec_exemple} for concrete examples of finitely presented group.

From now on, in any inequalities involving the $G$-actions denoted by ``$\lesssim$" (as in inequality (\ref{s1-C0-decomp})), it is understood that the implicit constants also depend on the number of generators $k$, the number of relations $p$, and the lengths $\{ m_j \}$ of the relations.

\medskip

For the $C^\infty-$smooth compact Riemannian manifold $(M,g)$, let
$\pi:G\to {\rm Isom}^\infty(M)$ be a morphism group, which defines a {\it $G-$action by $C^\infty-$smooth isometries}.
The $G-$action $\pi$ induces the representation on $\Gamma^\infty=\Gamma^\infty(M,TM)$:
\begin{equation}\label{repre_pi}
\pi(\gamma)_*v :=
\left(D\pi(\gamma)\cdot v\right)\circ \pi\left(\gamma^{-1}\right) ,\quad v\in\Gamma^\infty,
\end{equation}
which means the differential of $\pi(\gamma)$ evaluated at $\pi\left(\gamma^{-1}\right)$ and applied on $v\left(\pi\left(\gamma^{-1}\right)\right)$.
In this section lest us introduce some elementary properties related to the $G-$action by smooth isometries $\pi$
and the orthogonal decomposition induced by on the $L^2(M,TM)^k$.

\subsection{Generating isometries and $d_0-$Diophantineness}\label{sec_gener_iso}

With $\CS=\{\gamma_l\}_{1\leq l \leq  k}$ the set of generators of $G$,
let us define $d_0:L^2(M,TM) \to L^2(M,TM)^k $ by
\begin{equation}\label{defi_d0}
d_0v=\left(v-\pi(\gamma_l)_*v\right)_{1\leq l \leq  k},\quad v\in L^2(M,TM).
\end{equation}
We have hence its adjoint $d_0^*:L^2(M,TM)^k\to L^2(M,TM)$, with respect to the scalar product from $L^2(M,TM)$ induced on $L^2(M,TM)^k$,
\beq\label{d0*}d_0^*V =\sum_{i=1}^k(v_i-\pi(\gamma_i^{-1})_*v_i),\quad V=(v_i)_{1\leq i\leq k}\in L^2(M,TM)^k. \eeq
Both of $d_0$ and $d_0^*$ are the same with that introduced in Section \ref{sec_Dio-intro}.
According to (\ref{iso-laplacian}) -- (\ref{E_lamj}), we see that
\beq\label{invar-d0}d_0E_{\lam_j}\subset E_{\lam_j}^k,\qquad d_0^*E_{\lam_j}^k\subset E_{\lam_j}.
\eeq
Then the operator $(d_0\circ d_0^*)_j:=d_0\circ d_0^*\circ \PP_j$, with the projection $\PP_j$ defined in (\ref{Proj_i}), is self-adjoint operators on the finite-dimensional vector space $E_{\lam_j}^k$.

Related to the $d_0-$Diophantine condition in Definition \ref{defi-dio-gnrt}, let us introduce the definition of Diophantine property of group actions given by Dolgopyat \cite{Dolgo2002}, who defined a small-divisor condition for the subset of the group $G$.
\begin{defi}(Dolgopyat \cite{Dolgo2002}[Appendix A])\label{dioph-dolg}
Given a finitely presented group $G$ acting transitively on a smooth compact manifold $M$ by isometries $\pi$, the subset $S\subset G$ is said to be a {\bf Diophantine subset for $\pi$} if there are constants $\sigma>0$, $\tau\geq 0$ such that,  for all $j\in \N$ with $\lambda_j\neq 0$ and all $u\in E_{\lambda_j}$, there exists $\gamma\in S$ (depending on $u$) such that\begin{equation}\label{dioph-cond}
\|u -\pi(\gamma )_*u\|_{L^2} \geq\frac{\sigma}{\lambda_j^{\tau}}\|u\|_{L^2}.
\end{equation}
\end{defi}

If the set of generators $\CS$ of $G$ is a Diophantine subset for $\pi$ in the sense of Definition \ref{dioph-dolg}, then  the $G-$action $\pi$ is $d_0-$Diophantine in the sense of Definition \ref{defi-dio-gnrt}. Indeed, we have the following lemma.

\begin{lemma}\label{lemma_d_0-dio}
If there exist $\sigma>0$, $\tau\geq 0$ such that, for $u\in \IM(d_0^*\circ \PP_j)\subset E_{\lam_j}$, we have
\beq\label{gener-d0}
\|u-\pi(\gamma_{l(u)})_*u\|_{L^2}\geq \frac{\sigma}{(1+\lam_j)^\tau}\|u\|_{L^2},\eeq
for some $1\leq l(u)\leq k$, then the $G-$action $\pi$ is $(\sigma^2,2\tau)-d_0-$Diophantine.
\end{lemma}
\proof For the $G-$action by isometries $\pi$ with (\ref{gener-d0}) satisfied, we have
\begin{equation}\label{d0-Dio-condition}
\|d_0 u\|_{L^2} =\left(\sum_{l=1}^k \|u-\pi(\gamma_l)_* u\|^2_{L^2}\right)^\frac12\geq \frac{\sigma}{(1+\lam_j)^\tau}\|u\|_{L^2}, \qquad \forall \ u\in \IM(d_0^*\circ \PP_j).
\end{equation}

Since $(d_0\circ d_0^*)_j$ is self-adjoint, it is diagonalizable, and all its non-zero eigenvalues (if any exist), denoted by
$\mu_{j,1},\cdots,\mu_{j,K_j}$, are all positive.
Otherwise, since $\Ker d_0$ and $\IM d_0^*$ are orthogonal to each other, we have $E_{\lam_j}^k=\Ker(d_0\circ d_0^*)_j=\Ker(d_0^*\circ \PP_j)$.
Let $\{{\CE}_{j,1}, \cdots, {\CE}_{j,K_j}\}\subset E_{\lambda_j}^k$ be an orthonormal basis of the eigenvectors of $(d_0\circ d_0^*)_j$ associated with the non-vanishing eigenvalues $\{\mu_{j,1}, \cdots, \mu_{j,K_j}\}$.
Then every $V\in\IM(d_0 \circ\PP_j)$ can be decomposed along these eigenvectors, i.e.,
\begin{equation}\label{u_EquiForms}
V=\sum_{1\leq i\leq K_j}V_{j,i} \CE_{j,i}.
\end{equation}
According to (\ref{d0-Dio-condition}), we have, for $1\leq i \leq K_j$,
$$\mu_{j,i} =\|(d_0\circ d_0^*) \CE_{j,i} \|_{L^2}   \geq \frac{\sigma}{(1+\lam_j)^\tau}\|d_0^* \CE_{j,i} \|_{L^2} . $$
On the other hand, we have
$$\|d_0^*  \CE_{j,i} \|_{L^2}=\la  (d_0\circ d_0^*)  \CE_{j,i},  \CE_{j,i} \ra^\frac12=\mu_{j,i}^\frac12.$$
Hence we obtain that $\mu_{j,i}\geq  \sigma (1+\lam_j)^{-\tau} \mu_{j,i}^\frac12 $, which implies that
$$\mu_{j,i}\geq  \frac{\sigma^2}{(1+\lam_j)^{2\tau}},\qquad 1\leq i \leq K_j . $$
Therefore, for $V=\sum_{1\leq i\leq K_j}V_{j,i} \CE_{j,i}\in\IM(d_0 \circ \PP_j)$, we have
\beq\label{d0d0*DIO-proof} \|(d_0\circ d_0^*) V\|_{L^2} \ = \ \left(\sum_{1\leq i\leq K_j}\mu^2_{j,i} |V_{j,i}|^2\right)^\frac12
\geq  \frac{\sigma^2}{(1+\lam_j)^{2\tau}} \|V\|_{L^2}. \eeq
Hence $\pi$ is $(\sigma^2,2\tau)-d_0-$Diophantine.\qed

\begin{remark}\label{rmk_DioDol}
 If the set of generators $\CS$ of $G$ is a Diophantine subset for $\pi$ in the sense of Definition \ref{dioph-dolg}, it implicitly implies that $\Ker d_0 \subset E_0$. Hence the inequalities (\ref{gener-d0}) and (\ref{d0-Dio-condition}) hold for any $u\in E_{\lam_j}$, $\lambda_j\neq 0$. With the same proof, we obtain (\ref{d0d0*DIO-proof}) for all $V\in E_{\lam_j}^k$. Hence $\pi$ is $d_0-$Diophantine in the sense of Definition \ref{defi-dio-gnrt}. However, the converse of the above argument is not true. Examples are provided in Section \ref{sec_cyclic} and \ref{sec_sphere}.
\end{remark}

\subsection{Self-adjoint Box operator}\label{sec_op-box}

Let $C^i(G,L^2(M,TM))$ be the $i-$cochain on $G$ with values in $L^2(M,TM)$.
One can identify $0-$cochains $C^0(G,L^2(M,TM))$ with $L^2(M,TM)$, $1-$cochains $C^1(G,L^2(M,TM))$ with maps from $\CS$ to $L^2(M,TM)$, or equivalently $L^2(M,TM)^k$, and $2-$cochains $C^2(G,L^2(M,TM))$ with maps from $\CR$ to $L^2(M,TM)$, or equivalently $L^2(M,TM)^p$.
To the representation (\ref{repre_pi}), we can associate the Hoshchild complex (as introduced in (\ref{complex})):
$$L^2(M,TM) \stackrel{d_0}{\longrightarrow} C^1(G,L^2(M,TM)) \stackrel{d_1}{\longrightarrow} C^2(G,L^2(M,TM))  \stackrel{d_2}{\longrightarrow}  \cdots ,$$
or equivalently,
\begin{equation}\label{HoshchildComplex}
L^2(M,TM) \stackrel{d_0}{\longrightarrow} L^2(M,TM)^k  \stackrel{d_1}{\longrightarrow} L^2(M,TM)^p  \stackrel{d_2}{\longrightarrow}  \cdots ,
\end{equation}
where the differentials $d_0$ and $d_1$ can be written explicitly, i.e., $d_0$ is defined as in (\ref{defi_d0}) and,
for the finite words $\CW_j=\gamma_{l^{(j)}_1}\cdots \gamma_{l^{(j)}_{m_j}}$, $1\leq j\leq p$, in $R$,
\begin{equation}\label{defi_d1}
d_1 V=\left(\sum_{1\leq z \leq m_j}\pi\left(\prod_{i=1}^{z-1}\gamma_{l^{(j)}_i}\right)_*v_{l^{(j)}_z}\right)_{1\leq j\leq p},\quad V=(v_j)_{1\leq j\leq k}\in L^2(M,TM)^k
\end{equation}
where, for $1\leq l\leq k$, $v_{l+k}:=-\pi(\gamma_l^{-1})_*v_{l}$.

Let $d_1^*$ denote the adjoint of $d_1$, with respect to the scalar product from $L^2(M,TM)$ induced on $L^2(M,TM)^p$.
By direct computations, we have that $d_1^*:L^2(M,TM)^p\to  L^2(M,TM)^k$ is defined as
\beq\label{d1*}
d_1^*W= \left( \sum_{j=1}^p \sum_{w\in S_{j,l}}\pi(w^{-1})_*W_j \right)_{1\leq l \leq k},\quad W=(W_i)_{1\leq i\leq p}\in L^2(M,TM)^p,
\eeq
where, for a fixed relation $\CW_j=\gamma_{l^{(j)}_1}\cdots \gamma_{l^{(j)}_{m_j}}\in \CR$, we define
$$S_j:=\left\{\gamma_{l^{(j)}_1}\cdots \gamma_{l^{(j)}_q}: 1\leq q\leq m_j\right\}$$ as the set of subwords of $\CW_j$ that begin with $\gamma_{l^{(j)}_1}$, and, for $1\leq l\leq k$, we define $$S_{j,l}:=\left\{w\in S_j\cup\{e\}: w\gamma_{l}\in S_j\right\},$$
which is the set of words $w$ in $S_j$ (including the identity element $e$) such that
$w\gamma_{l}$ is also a word in $S_j$.
In view of (\ref{iso-laplacian}) -- (\ref{E_lamj}), we have that
\beq\label{invar-d1} d_1E_{\lam_j}^k\subset E_{\lam_j}^p,\qquad d_1^*E_{\lam_j}^p\subset E_{\lam_j}^k.\eeq
Let $(d_1^*\circ d_1)_j:= d_1^*\circ d_1\circ \PP_j$, with the projection $\PP_j$ defined in (\ref{Proj_i}), be the self-adjoint operator on the finite-dimensional vector space $E_{\lam_j}^k$.
The connection between the independent generating isometries is expressed through the group relations in $\CR$.
More properties of $d_1$ will be introduced in Section \ref{sec_GAdiffeo}.

Following the idea of Lombardi-Stolovitch \cite{LS10} developed in the context of germs of holomorphic vector fields at a singular point, we introduce the self-adjoint Box operator
\beq\label{BoxOp}
\square=d_0\circ d_0^*+ d_1^*\circ d_1: L^2(M,TM)^k\to L^2(M,TM)^k.\eeq
Combining (\ref{invar-d0}) and (\ref{invar-d1}), we see that $\square$ is invariant on $E_{\lambda_j}^k$ for every $j\in\N$.
Hence $\square_j:=\square\circ \PP_j$ is a self-adjoint operator on the finite-dimensional vector space $E_{\lam_j}^k$.
Hence it is diagonalizable and its non-zero eigenvalues $\tilde\mu_{j,1},\cdots,\tilde\mu_{j,\tilde K_j}$ (if existing) are all positive.
Through elementary properties of Hilbert space, we have
\begin{equation}\label{decomVi_ImKer}
L^2(M,TM)^k= {\rm Ker} \square \bigoplus {\rm Im} \square,\qquad E_{\lam_j}^k  = {\rm Im}\square_j \bigoplus {\rm Ker}\square_j.
\end{equation}
Let $\CH$ be the projection operator from $L^2(M,TM)^k$  onto ${\rm Ker} \square$.


\begin{lemma}\label{lemmaH_d0} We have ${\CH}\circ d_0=0$ and ${\CH}\circ d_1^*=0$ on $L^2(M,TM)^k$.
\end{lemma}
\proof
In view of the Peter-Weyl decomposition (\ref{peter-weyl}) and (\ref{E_lamj}), it is sufficient to show that, for $\square_j=\square\circ\PP_j$, $j\in \N$, $d_0v ,d^*_1W \in{\rm Im}\square_j$ for any $v \in E_{\lam_j}$ and $W\in E_{\lam_j}^p$.

Since  $E_{\lam_j}^k={\rm Ker}\square_j\bigoplus{\rm Im}\square_j$,
we have the unique decomposition of $d_0v$ for $v\in E_{\lam_j}$:
$$d_0v=w^{\rm Ker}+w^{\rm Im},\qquad w^{\rm Ker} \in {\rm Ker} \square_j,\quad w^{\rm Im} \in {\rm Im} \square_j,$$
which implies that
$$
\square(d_0v-w^{\rm Im})=\left(\left(d_0\circ d_0^*+d_1^*\circ d_1\right)\circ d_0\right)v-\square w^{\rm Im}
=\left(d_0\circ d_0^*\circ d_0\right)v-\square w^{\rm Im}=0.
$$
Recalling that $v$ is arbitrarily chosen in $E_{\lam_j}$, we have that ${\rm Im} (d_0\circ d_0^*\circ d_0\circ \PP_j)\subset{\rm Im} \square_j.$

Now let us show that $d_0 v\in {\rm Im} (d_0\circ d_0^*\circ d_0\circ \PP_j)$ for $v\in E_{\lam_j}$, which will complete the proof of ${\CH}\circ d_0=0$.
Since $E_{\lam_j}={\rm Ker}(d_0\circ \PP_j)\bigoplus{\rm Im}(d_0^*\circ \PP_j)$, we have the unique decomposition for $v$:
$$v=v^{\rm Ker}+v^{\rm Im},\qquad v^{\rm Ker}\in {\rm Ker}(d_0\circ \PP_j), \quad v^{\rm Im}\in {\rm Im}(d_0^*\circ \PP_j).$$
Then, there exists some $ u\in E_{\lam_j}^k$ satisfying $d_0^*u=v^{\rm Im}$ such that
$$d_0v=d_0v^{\rm Im}=(d_0\circ d_0^*)u.$$
In view of the fact that $E_{\lam_j}^k={\rm Ker}(d_0^*\circ \PP_j)\bigoplus{\rm Im}(d_0\circ \PP_j)$, we have the unique decomposition for $u$:
$$u=u^{\rm Ker}+u^{\rm Im},\qquad u^{\rm Ker}\in {\rm Ker}(d_0^*\circ \PP_j), \quad u^{\rm Im}\in {\rm Im}(d_0\circ \PP_j).$$
Hence, there exists some $ z\in E_{\lam_j}$ satisfying $d_0 z=u^{\rm Im}$ such that
$$ d_0v=(d_0\circ d_0^*)u=(d_0\circ d_0^*)u^{\rm Im}=(d_0\circ d_0^*\circ d_0)z.$$
This implies that $d_0v\in {\rm Im}(d_0\circ d_0^*\circ d_0\circ \PP_j)\subset{\rm Im} \square_j$.

In a similar way, we have $d_1^*W\in \IM\square$, which implies $\CH\circ d_1^*=0$. \qed

\smallskip

Since $d_1\circ d_0=0$, according to Lemma \ref{lemmaH_d0}, we have $\IM\square=\IM d_0 \bigoplus  \IM d_1^*$, which gives the decomposition of $L^2(M,TM)^k$ w.r.t. the $G-$action $\pi$:
\begin{equation}\label{L2k-decomp}
L^2(M,TM)^k=\Ker \square \bigoplus \IM d_0 \bigoplus  \IM d_1^*.
\end{equation}
Let us define $\D_0$ and $\D_1$ as the projection form $L^2(M,TM)^k$ onto $\IM d_0$ and $\IM d_1^*$ respectively.
%

\smallskip

\begin{lemma}\label{lem-H1}
Given a $G-$action $\pi$ by isometries on a smooth manifold $M$, we have
$$H^1(G, E_{\lambda_j}):=\Ker \,(d_{1}\circ\PP_j )/ \IM \, (d_{0}\circ\PP_j )\approxeq {\rm Ker}\square_j.$$
\end{lemma}
\proof Given $f\in E_{\lambda_j}^k$, we have that $\la \square_j f, f\ra = \|d_0^*f\|_{L^2}^2+\|d_1f\|_{L^2}^2$.
	Hence, $\square_j f=0$ if and only if $d_0^*f=0$ and $d_1f=0$. The latter means that $f$ is a 1-cocycle, i.e., $f\in Z^1(G,E_{\lambda_j})$. The former means that $f$ is orthogonal to $\text{Im}\,d_0$. Hence, $f$ belongs to a space isomorphic to $Z^1(G,E_{\lambda_j})/\text{Im}(d_{0}\circ\PP_j )=H^1(G, E_{\lambda_j})$.
\qed

\subsection{Group action by diffeomorphisms on $M$}\label{sec_GAdiffeo}

Based on the properties of the $G-$action by smooth isometries $\pi$ on $M$ introduced in Section \ref{sec_gener_iso} and \ref{sec_op-box}, let us focus on some elementary properties of the $G-$action $\pi^u$ by smooth diffeomorphisms on $M$, considered as the perturbation of $\pi$, of the form
$$\pi^u(\gamma)=\Exp\{u(\gamma)\}\circ\pi(\gamma),\quad \gamma\in \CS,$$
where $\Exp$ is the exponential map introduced in Section \ref{secExp} and $u:\CS\to \Ga^\infty(M,TM)$ is $C^0$ sufficiently small.

Let us define the subset, in the neighborhood of origin in $L^2(M,TM)^k$,
$${\bf G}_\pi:=\left\{\CU=(u_l)_{1\leq l\leq k}\in \Ga^\infty(M,TM)^k: \begin{array}{l}
\pi^u  \ {\rm is \ a  } \ G- {\rm action \ by \ smooth  }\\[1mm]
{\rm  \ diffeomorphisms \  if } \ u(\gamma_l)=u_l
\end{array}\right\}.$$
It is obvious that ${\bf G}_\pi$ is non-empty since $0\in{\bf G}_\pi$.

\begin{lemma}\label{lem-d1-quadra}
For $\CU\in {\bf G}_\pi$ with $\|\CU\|_{C^1}$ sufficiently small, we have
$$\|(d_1^*\circ d_1)\CU\|_{L^2}\lesssim \|\CU\|_{C^{1}}\|\CU\|_{C^{0}}. $$
\end{lemma}
\proof For the $G-$action by smooth diffeomorphisms $\pi^u$ with $(u(\gamma_l))_{1\leq l\leq k}=\CU\in {\bf G}_\pi$,
we have, for the relation word
$\CW_j=\gamma_{l_1} \cdots \gamma_{l_{m_j}}=\gamma_{l^{(j)}_1} \cdots \gamma_{l^{(j)}_{m_j}}\in \CR$,
\begin{equation}\label{Gpi-Id}
\pi^u(\gamma_{l_1} \cdots \gamma_{l_{m_j}})= \pi(\gamma_{l_1} \cdots \gamma_{l_{m_j}})={\rm Id},
\end{equation}
which implies that $u(\CW_j)=u(\gamma_{l_1} \cdots \gamma_{l_{m_j}})=0$.

For $1\leq q\leq m_j$, let $w_q=w^{(j)}_q$ be the sub-word of $\CW_j$ satisfying $w_q:=\gamma_{l_{q}}\cdots\gamma_{l_{m_j}}$.
 For $2\leq q\leq m_j$, we have $w_{q-1}=\gamma_{l_{q-1}} w_q$, and
 \begin{eqnarray*}
 \pi^u(w_{q-1})
&=&\Exp\{u(\gamma_{l_{q-1}})\}\circ\pi(\gamma_{l_{q-1}})\circ\Exp\{u(w_q)\}\circ\pi(w_q)\\
&=&\Exp\{u(\gamma_{l_{q-1}})\}\circ\Exp\{\pi(\gamma_{l_{q-1}})_*u(w_q)\}\circ \pi(\gamma_{l_{q -1}} w_q)\\
&=&\Exp\{u(\gamma_{l_{q-1}})+\pi(\gamma_{l_{q-1}})_*u(w_q)+s_1(u(\gamma_{l_{q-1}}),\pi(\gamma_{l_{q-1}})_*u(w_q))\} \circ \pi(w_{q-1}).
\end{eqnarray*}
 Hence, we obtain
\begin{equation}\label{rec-uw}
u(w_{q-1})=u(\gamma_{l_{q-1}})+\pi(\gamma_{l_{q-1}})_*u(w_q)+s_1(u(\gamma_{l_{q-1}}),\pi(\gamma_{l_{q-1}})_*u(w_q)),
\end{equation}
where, by (\ref{s1-C0-decomp}) in Proposition \ref{prop_s1-smooth} (with $w_1=u(\gamma_{l_{q-1}})$ and $w_2=0$),
$$
\|s_1(u(\gamma_{l_{q-1}}),\pi(\gamma_{l_{q-1}})_*u(w_q))\|_{C^0}\lesssim  \|u(\gamma_{l_{q-1}})\|_{C^1}\|u(w_q)\|_{C^0}.
$$
Therefore, we have
 $$\|u(w_{q-1})\|_{C^0}\lesssim \|u(\gamma_{l_{q-1}})\|_{C^0}+(1+\|u(\gamma_{l_{q-1}})\|_{C^1})\|u(w_q)\|_{C^0}\lesssim \|u(\gamma_{l_{q-1}})\|_{C^0}+ \|u(w_q)\|_{C^0}. $$
Hence, by induction, for $q=2,\cdots, m_j+1$,
\begin{equation}\label{estim-wq-Ck}
\|u(w_{q-1})\|_{C^0}\lesssim \sum_{r=q-1}^{m_j}\|u(\gamma_{l_{r}})\|_{C^0}\lesssim \|\CU\|_{C^0}.
\end{equation}
On the other hand, the induction with (\ref{rec-uw}) implies that
 \begin{eqnarray}
 u(w_{1})&=&u(\gamma_{l_{1}})+\pi(\gamma_{l_{1}})_*u(w_2)+s_1(u(\gamma_{l_{1}}),\pi(\gamma_{l_{1}})_*u(w_2))\nonumber\\
 &=& u(\gamma_{l_{1}})+\pi(\gamma_{l_{1}})_*u(\gamma_{l_{2}})+\pi(\gamma_{l_{1}}\gamma_{l_{2}})_*u(w_3)\nonumber\\
& & + \, s_1(u(\gamma_{l_{1}}),\pi(\gamma_{l_{1}})_*u(w_2))+ \pi(\gamma_{l_{1}})_*s_1(u(\gamma_{l_{2}}),
\pi(\gamma_{l_{2}})_*u(w_3))\nonumber\\
 &\vdots &  \nonumber \\
  &=&u(\gamma_{l_{1}})+\pi(\gamma_{l_{1}})_*u(\gamma_{l_{2}})+\cdots +\pi(\gamma_{l_{1}}\cdots\gamma_{l_{m_j-1}})_*u(w_{m_j})\label{d1Uj}\\
& &+ \,  s_1(u(\gamma_{l_{1}}),\pi(\gamma_{l_{1}})_*u(w_2))
+   \sum_{q=2}^{m_j-1} \pi\left(\prod_{\ell=1}^{q-1}\gamma_{l_{\ell}}\right)_*s_1(u(\gamma_{l_q}), \pi(\gamma_{l_q})_*u(w_{q+1})).\label{Fj-1}
\end{eqnarray}
Recalling (\ref{defi_d1}), we see that the sum of terms in (\ref{d1Uj}) is $(d_1\CU)_j$. Denote the sum of the terms in (\ref{Fj-1}) by $\tilde\CU_j$. By (\ref{estim-wq-Ck}), as well as (\ref{s1-C0-decomp}) in Proposition \ref{prop_s1-smooth}, we obtain
$$   \|\tilde\CU_j\|_{C^0}\lesssim\|\CU\|_{C^1}\|\CU\|_{C^0}.$$

Now, for $1\leq j\leq p$, we have
$u(\CW_j)=u(w_{1})= (d_1\CU)_j+ \tilde\CU_j$.
Since $u(\CW_j)=0$, we have $(d_1\CU)_j=- \tilde\CU_j$, and hence, for $R\in\N$,
$$\|d_1\CU\|_{C^0}=\left(\sum_{j=1}^p\|(d_1\CU)_j\|^2_{C^0}\right)^\frac12= \left(\sum_{j=1}^p\|\tilde\CU_j\|_{C^0}\right)^\frac12\lesssim  \|\CU\|_{C^{1}}\|\CU\|_{C^0}. $$
 Then, in view of (\ref{d1*}) and through Proposition \ref{propNorm-CrHr}, we have
 $$\|(d_1^*\circ d_1)\CU\|_{L^2}\lesssim \| d_1\CU\|_{L^2}\lesssim  \| d_1\CU\|_{C^0} \lesssim   \|\CU\|_{C^{1}}\|\CU\|_{C^0}. \qed$$

\section{Examples of local rigidity -- Applications of Theorems}\label{sec_exemple}

Let us provide some examples of local rigidity in concrete situations of finitely presented group $G$ and of compact Riemannian manifold $M$.

\subsection{Abelian group action by isometries}

Let $G$ be an abelian group, which means, for the generators $\{\gamma_l\}_{1\leq l \leq k}$, the relations are presented by the $\frac{k(k-1)}{2}$ words
\beq\label{abelian-relation}
\CW_{i, l}:= \gamma_i\gamma_l\gamma_{i}^{-1}\gamma_{l}^{-1},\quad 1\leq i < l \leq k.
\eeq
Consider the $G-$action $\pi$ by smooth isometries on $M$.
\begin{prop}\label{propDioAbelian} For the abelian group action $\pi$, we have
\beq\label{box-abelian}
\square u=((d_0^*\circ d_0 ) u_i)_{1\leq i\leq k},\quad u=(u_i)_{1\leq i\leq k}\in (\Ga^\infty)^k.\eeq
Moreover, if $\pi$ is $d_0-$Diophantine, then $\pi$ is $\square-$Diophantine.
\end{prop}

\proof Let us define $L_i:L^2(M,TM)\to L^2(M,TM)$, $1\leq i\leq k$, by
$$L_iv:=v-\pi(\gamma_i)_*v,\qquad v\in L^2(M,TM). $$
It is easy to verify that its adjoint is given by $L_i^*v=v-\pi(\gamma_i^{-1})_*v$, since $\pi$ is an action by isometries.
The group relation (\ref{abelian-relation}) guarantees that $L_iL_l^*=L_l^*L_i$ for $1\leq i,l\leq k$.

For $u=(u_l)_{1\leq l\leq k}\in L^2(M, TM)^{k}$, we have, by computations with (\ref{defi_d0}) and (\ref{d0*}),
\beq\label{d0d0*}
(d_0\circ d_0^*)u= \left(\sum_{l=1}^kL_iL_l^*u_l\right)_{1\leq i\leq k}.
\eeq
For one relation $\CW_{i, l}$ given in (\ref{abelian-relation}), for $u=(u_l)_{1\leq l\leq k}\in L^2(M, TM)^{k}$, by the definition (\ref{defi_d1}), we have, for any $1\leq i<l\leq k$,
\begin{eqnarray*}
(d_1u)_{i, l}&=&u_i+\pi(\gamma_i)_*u_l-\pi(\gamma_i\gamma_l\gamma_i^{-1})_*u_{i}-\pi(\gamma_i\gamma_l\gamma_{i}^{-1}\gamma_{l}^{-1})_*u_{l}\\
&=& (u_i-\pi(\gamma_l)_*u_i)-(u_l-\pi(\gamma_i)_*u_l)
 \ = \ L_l u_i -L_i u_l .\end{eqnarray*}
For $W=(W_{i, l})_{1\leq i< l \leq k}\in L^2(M,TM)^p$ with $p=\frac{k(k-1)}{2}$,
\begin{eqnarray*}
\la d_1u, W \ra&=& \sum_{1\leq i<l\leq k}\la L_l u_i -L_i u_l ,W_{i, l}\ra\\
&=&\sum_{1\leq i<l\leq k}\left(\la u_i, L_l^* W_{i, l} \ra - \la u_l, L_i^* W_{i, l}\ra \right)\\
&=&\sum_{i=1}^k\left\la u_i, \sum_{1\leq l\leq k\atop{l>i}} L_l^* W_{i, l} -\sum_{1\leq l\leq k\atop{l<i}} L_l^* W_{l, i}\right\ra,
\end{eqnarray*}
which implies that
$$d_1^* W =\left(\sum_{1\leq l\leq k\atop{l>i}} L_l^* W_{i, l}-\sum_{1\leq l\leq k\atop{l<i}} L_l^* W_{l, i}\right)_{1\leq i \leq k}.$$
Then, by direct computations, we obtain
\begin{eqnarray}
(d_1^*\circ d_1)u&=&\left(\sum_{1\leq l\leq k\atop{l>i}}\left( L_l^*L_l u_i -L_l^*L_i u_l \right)
 -  \sum_{1\leq l\leq k\atop{l< i}}   \left(L_l^*L_i u_l -L_l^*L_l u_i\right)\right)_{1\leq i\leq k}\nonumber\\
&=&\left(\sum_{1\leq l\leq k\atop{l\neq i}}L_l^*L_l u_i-\sum_{1\leq l\leq k\atop{l\neq i}} L_l^*L_i u_l\right)_{1\leq i\leq k}.\label{d1*d1-abelian}
\end{eqnarray}
Combining (\ref{d0d0*}) and (\ref{d1*d1-abelian}), and recalling that $L_l^*L_i=L_iL_l^*$, we obtain (\ref{box-abelian}) by
$$
\square u=(d_0\circ d_0^*+ d_1^*\circ d_1)u=\left(\sum_{1\leq l\leq k}L_l^*L_l u_i\right)_{1\leq i\leq k}=\left((d_0^*\circ d_0 ) u_i\right)_{1\leq i\leq k}.
$$

If $\pi$ is $(\sigma,\tau)-d_0-$Diophantine, then, for every $j\in\N$, all the non-zero eigenvalues of $(d_0\circ d_0^*)_j$ are greater than $\sigma(1+\lambda_j)^{-\tau}$, which implies that $\pi$ is $(\sigma,\tau)-\square-$Diophantine, since $d_0^*\circ d_0\circ \PP_j$ has the same non-zero eigenvalues with $(d_0\circ d_0^*)_j$.\qed

\medskip

According to Theorem \ref{thm-geom0} and Proposition \ref{propDioAbelian}, for the $d_0-$Diophantine abelian action $\pi$ by smooth (or analytic) isometries on $M$, if its first cohomology $H^1(G,L^2(M,TM))$ vanishes, then it is smoothly (or analytically) rigid.

As a concrete corollary of Theorem \ref{thmmain}, we obtain an analytic version of results by Moser \cite{Moser-cercle}, 
and Petkovic \cite{petko-tore} relative to simultaneous conjugacy of a commutative family of perturbations of rotations on the torus to rotations.

\begin{thm}\label{thm-tore0} Let $\CG=\{e_1,\cdots, e_{k}\}$ be the canonical basis of $\Z^{k}$. Let $\pi$ be a $\Z^{k}-$action by translations on the torus $\T^d$ defined by
\beq\label{GA_rot_tori}
\pi(e_i):x\mapsto x+\al_i,\qquad i=1,\cdots, k, \eeq
 with the translation vectors $\al_i\in \R^d$ satisfying the simultaneous Diophantine condition: there exist $c$, $\tau>0$, such that for all $({\bf k},l)\in \Z^d\setminus \{0\}\times \Z$,
\beq\label{simuDio_tori}
\max_{1\leq i\leq k}|\la {\bf k},\al_i\ra- l\pi|\geq \frac{c}{|{\bf k}|^{\tau}}.
\eeq
Then any $\Z^{k}-$action $\pi_P$ by analytic diffeomorphisms on $\T^d$, which is sufficiently small perturbation of $\pi$ and isotopic to the identity, is analytically conjugate to $\pi$, if, for each $i$, the rotation vector $\al_i$ belongs to the convex hull of rotation set of $\pi_P(e_i)$.
\end{thm}

\proof
For $u\in L^2(\T^d, T \T^d) $, it can be presented as
$$u(x)=\sum_{m=1}^d u^m(x) \partial_{x_m},\quad x=(x_m)_{1\leq m\leq d},$$
where $u^m$, $1\leq m\leq d$, are scalar functions on $\T^d$, with the Fourier expansion
$$u^m(\bullet)=\sum_{{\bf k}\in\Z^d}\left(\hat u^m_{\bf k} e^{{\rm i}\la{\bf k},\bullet\ra}+\bar{\hat u}^m_{\bf k} e^{-{\rm i}\la{\bf k},\bullet\ra}\right),\quad u^m\in L^2(\T^d).$$
For $|\Delta_{T \T^d}|^\frac12$ on $L^2(\T^d, T \T^d)$, the eigenvalues are $\lambda_j=j$, $j\in\N$, associated with the eigenspaces
\beq\label{eigenspace-tore}
E_{j}={\rm Vect}\left(\cos\la{\bf k},\bullet\ra\partial_{x_m} , \sin\la{\bf k},\bullet\ra\partial_{x_m}\right)_{{\bf k}\in\Z^d, |{\bf k}|=j \atop{1\leq m\leq d}}.\eeq
Then, for $i=1,\cdots, k$, we have
\begin{eqnarray*}
\left(u -\pi(e_i )_*u\right)(\bullet)&=&u(\bullet)-u(\bullet+\alpha_i)\\
&=&\sum_{m=1}^d\sum_{{\bf k}\in\Z^d \atop{|{\bf k}|=j}}\left(\left(1-e^{{\rm i}\la{\bf k},\alpha_i\ra}\right)\hat u^m_{\bf k} e^{{\rm i}\la{\bf k},\bullet\ra}+\left(1-e^{-{\rm i}\la{\bf k},\alpha_i\ra}\right)\bar{\hat u}^m_{\bf k} e^{-{\rm i}\la{\bf k},\bullet\ra}\right)\partial_{x_m},
\end{eqnarray*}
which implies that $\left\|u -\pi(e_i )_*u\right\|_{L^2}^2=
4\sum_{m=1}^d\sum_{{\bf k}\in\Z^d \atop{|{\bf k}|=j}}\sin^2\frac{\la{\bf k},\alpha_i\ra}{2}|\hat u^m_{\bf k}|^2.$
Under the simultaneous Diophantine condition (\ref{simuDio_tori}), the set $\CG=\{e_1,\cdots, e_{k}\}$, generators of $\Z^{k}$, is a Diophantine subset of $\Z^{k}$ in the sense of Definition \ref{dioph-dolg}.


For the $\Z^{k}-$action $\pi_P$ by analytic diffeomorphisms on $\T^d$, $\{\pi_P(e_i)\}$ is commuting.
Then, according to \cite{petko-tore}[Theorem 5], the family of analytic diffeomorphisms $\{\pi_P(e_i)\}$ is simultaneously smoothly conjugate to $\{\pi(e_i)\}$ through a near-identity transformation in $C^\infty$ topology. Hence, for any $R\in\N^*$, $\pi_P$ is $C^R$ conjugate to $\pi$.
According to \rt{thmimpressive}, $\{\pi_P(e_i)\}$ is simultaneously analytically conjugate to $\{\pi(e_i)\}$.\qed

\begin{remark} The commutativity of $\{\pi_P(e_i)\}$, as well as the properties of the Box operator $\square$, are not used in the above proof. Indeed, it is only used in \cite{petko-tore} to show smooth rigidity, which ensures the hypothesis of \rt{thmimpressive}.
\end{remark}

\subsection{Cyclic group action by isometries}\label{sec_cyclic}

Let $G$ be a cyclic group of order $n$, $n\geq 2$, which means, for the generator $\gamma$, the only relation is presented by the word $\gamma^n$.
In this situation we have $k=p=1$.

Let $\pi$ be a $G-$action by smooth isometries on $M$. We have the specific representation of $L^2(M,TM)$ regarding the decomposition (\ref{L2k-decomp}).

\begin{lemma}\label{lemma-cyclic}
For such a $G-$action $\pi$, we have $\Ker \square=0$, and, for $u\in L^2(M,TM)$, we have the unique decomposition
\begin{equation}\label{decomp-cyclic}
u=\frac1{n^2}\sum_{1\leq l\leq \lfloor\frac{n}2\rfloor} y_l   (d_0\circ d_0^*)^{l}u + \frac1{n^2} (d_1^*\circ d_1) u\in \IM(d_0\circ d_0^*)\bigoplus\IM(d_1^*\circ d_1),
\end{equation}
with the coefficients $\{y_l\}_{1\leq l\leq \lfloor\frac{n}2\rfloor}\subset \Z$ depending only on $n$.\end{lemma}
\begin{remark} In view of Lemma \ref{lem-H1}, for the $G-$action $\pi$,
the first cohomology $H^{1}(G,L^2(M,TM))=0$.
\end{remark}
\proof According to (\ref{defi_d0}), (\ref{d0*}) and (\ref{defi_d1}), (\ref{d1*}), we have, for $u\in L^2(M,TM)$,
$$(d_0\circ d_0^*)u=2u-\pi(\gamma)_*u-\pi(\gamma^{n-1})_*u,\quad
(d_1^*\circ d_1)u=n(u+\pi(\gamma)_*u+\cdots+\pi(\gamma^{n-1})_*u). $$
By direct computations for $n= 2$ and $3$, we have the decomposition (\ref{decomp-cyclic}):
$$u=\left\{\begin{array}{cl}
\displaystyle \frac14 (d_0\circ d_0^*)u+\frac14 (d_1^*\circ d_1)u , & n=2\\[4mm]
\displaystyle \frac13 (d_0\circ d_0^*)u+\frac19 (d_1^*\circ d_1)u , & n=3
\end{array}
 \right. . $$
 For general $n\geq 3$, for $j\leq \lfloor\frac{n}2 \rfloor-1$, assume that
 \begin{eqnarray*}
  (d_0\circ d_0^*)^j u&=&2c_0^ju+c_1^j\left(\pi(\gamma)_*u+\pi(\gamma^{n-1})_*u\right)+\cdots +c_j^j\left(\pi(\gamma^j)_*u+\pi(\gamma^{n-j})_*u\right),
 \end{eqnarray*}
 for some suitable coefficients $\{c_l^j\}_{0\leq l \leq j}\subset \Z$.
Then we have
 \begin{eqnarray}
  (d_0\circ d_0^*)^{j+1}u&=&2(2c_0^j- c_1^j)u+(2c_1^j-2c_0^j-c_2^j)\left(\pi(\gamma)_*u+\pi(\gamma^{n-1})_*u\right)\nonumber\\
  & & + \, (2c_2^j-c_1^j-c^j_3)\left(\pi(\gamma^2)_*u+\pi(\gamma^{n-2})_*u\right)\nonumber\\
 & & + \, \cdots + (2c_j^j-c_{j-1}^j)\left(\pi(\gamma^j)_*u+\pi(\gamma^{n-j})_*u\right)\nonumber\\
 & & - \, c_j^j\left(\pi(\gamma^{j+1})_*u+\pi(\gamma^{n-j-1})_*u\right),\label{termcjj}
 \end{eqnarray}
where the term in (\ref{termcjj}) can be $-2c_j^j\pi(\gamma^{j+1})_*u$ if $n$ is even and $j+1=\frac{n}2$.
For the coefficients $\{c_l^j\}_{0\leq l \leq j}$, $1\leq j\leq \lfloor\frac{n}2 \rfloor$,
they obey the recurrence rules
$$\begin{array}{lllll}
c_0^{j+1}= 2c_0^j- c_1^j,&c_1^{j+1}= 2c_1^j-2c_0^j-c_2^j,& & &  \\[2mm]
c_l^{j+1} = 2c_l^j-c_{l-1}^j-c_{l+1}^j,& 2\leq l\leq j-1,& c_j^{j+1} = 2c_j^j-c_{j-1}^j,& &c_{j+1}^{j+1} =- c_{j}^{j}.
\end{array}
$$
It is easy to verify that, for any $1\leq j\leq \lfloor\frac{n}2 \rfloor$,
\begin{equation}\label{recurrence_cjl}
c_{j}^{j}=(-1)^j,\qquad \sum_{l=0}^j  c_l^{j} =0.
\end{equation}

Now for $J=\lfloor\frac{n}2 \rfloor$, let us find $\{y_{l}\}_{1\leq l\leq J}$ such that
\begin{equation}\label{decomp_u-n}
(d_1^*\circ d_1) u+\sum_{1\leq l\leq J}   y_l   (d_0\circ d_0^*)^{l}u = n^2 u.\end{equation}
 Then Eq. (\ref{decomp-cyclic}) is of that form.
 Recalling $(d_1^*\circ d_1)u=n(u+\pi(\gamma)_*u+\cdots+\pi(\gamma^{n-1})_*u)$, 
 we see that \re{decomp_u-n} reads
 $$ \alpha_0 u+\sum_{j=1}^J\alpha_j(\pi(\gamma^j)_*u+\pi(\gamma^{n-j})_*u)=n^2u,$$
 with the coefficients $\{\alpha_j\}_{0\leq j\leq J}$ defined as
 $$\alpha_0:=n+2\sum_{l=0}^Jc_0^l y_l,\qquad  \alpha_j:=n+\sum_{l=j}^Jc_j^ly_l,\quad j\geq 1.$$
It is sufficient to solve the upper-triangular linear system:
\begin{itemize}
\item if $n$ is odd,
\begin{align*}
 c_1^1 y_1 +  c_1^2 y_2+\cdots +c_{1}^{J-1} y_{J-1}+c_1^J y_J=& -n ,\\
  c_2^2 y_2+\cdots +c_{2}^{J-1} y_{J-1}+c_2^J y_J=& -n, \\
  \cdots \cdots \cdots \cdots & \\
   c_{J-1}^{J-1} y_{J-1}+c_{J-1}^J y_J=& -n, \\
   c_{J}^J y_J=& -n,
  \end{align*}
 \item if $n$ is even, the last equation is replaced by $c_{J}^J y_J= -\frac{n}2$ and the others are the same with the above system,
\end{itemize}
that is $\alpha_1=\cdots=\alpha_J=0$ together with $\alpha_0=n^2$.
Since the upper-triangular system has $c_{l}^l=(-1)^l$ as the diagonal coefficients, the solution exists in $\Z^J$ and is unique.
Moreover, by summing all equations in the system, combining with (\ref{recurrence_cjl}) we obtain
$$c_0^1 y_1 +  c_0^2 y_2+\cdots +c_{0}^{J-1} y_{J-1}+c_0^J y_J=\left\{
\begin{array}{cl}
Jn ,& n \ {\rm is \ odd}\\
Jn-\frac{n}2,& n \ {\rm is \ even}
\end{array}
\right. = \frac{(n-1)n}2.$$
Hence we have 
$$\alpha_0=n+2\left(c_0^1 y_1 +  c_0^2 y_2+\cdots +c_{0}^{J-1} y_{J-1}+c_0^J y_J\right)=n^2.$$
The decomposition (\ref{decomp_u-n}) is shown for $u$, which implies that
 $$u=\frac1{n^2}\sum_{1\leq l\leq J} y_l   (d_0\circ d_0^*)^{l}u + \frac1{n^2}(d_1^*\circ d_1) u \in \IM(d_0\circ d_0^*)\bigoplus\IM(d_1^*\circ d_1) .$$
This decomposition shows that $u\in \IM\square$ for any $u\in L^2(M,TM)$, hence $\Ker\square=0$.\qed

\smallskip

According to Theorem \ref{thm-geom0}, any $\square-$Diophantine cyclic group action by smooth (or analytic) isometries on the smooth (or analytic) compact Riemannian manifold $M$ is smoothly (or analytically) rigid. Such isometry is usually characterized by the periodic feature.
Let us give a concrete example for $M=\T^d$, $d\geq 1$, where the generator is not Diophantine in the sense of Definition \ref{dioph-dolg} but the action is $\square-$Diophantine.
\begin{prop}\label{prop_box-not-dolg}
There exists an action by isometry on $\T^d$ the generator of which does not define a Diophantine set (in the sens of Dolgopyat \rd{dioph-dolg}) whereas the action is $\square-$Diophantine.
\end{prop}

\proof
Let $\alpha=\left(\alpha_1,\cdots, \alpha_d\right)=2\pi\left(n_1^{-1},\cdots, n_d^{-1}\right)\in\T^d$ with the integers $n_1,\cdots, n_d\geq 2$ and pairwise coprime. Then $(\prod_{1\leq l\leq d}n_l )\alpha\in 2\pi \Z^d$.

For $n:=\prod_{1\leq l\leq d}n_l $, let us define the $n-$periodic translation $\pi:\T^d\circlearrowleft$ as $\pi:x\mapsto x+\alpha$.
Consider the $n-$periodic $F\in {\rm Diff}^\infty(\T^d)$ (resp. ${\rm Diff}^\omega(\T^d)$) with
\begin{equation}\label{F_defi}
F(x)=x+ \alpha + f(x),\quad F^{\circ n}(x) =x,\qquad x\in \T^d,
\end{equation}
where $f\in {\rm Diff}^\infty(\T^d)$ (resp. ${\rm Diff}^\omega(\T^d)$) is sufficiently small. Both $\pi$ and $F$ are $G-$actions by diffeomorphims on $\T^d$,
where $G$ is the cyclic group of order $n$ as introduced in the beginning of subsection.

For $u\in L^2(\T^d, T\T^d)$,
according to (\ref{defi_d0}), (\ref{d0*}) and (\ref{defi_d1}), (\ref{d1*}), we have
$$((d_0\circ d_0^*)u)(x)=2u(x)-u(x+\alpha)-u(x-\alpha),\quad
((d_1^*\circ d_1)u)(x)=n\sum_{l=0}^{n-1}u(x+l \alpha). $$
With the decomposition (recalling (\ref{eigenspace-tore})),
$$u(x)=\sum_{j\in\N}\sum_{{\bf k}\in\Z^d\atop{|{\bf k}|=j}}\sum_{m=1}^d\left(a_{\bf k} \cos\la {\bf k},x \ra+b_{\bf k} \sin\la {\bf k},x \ra\right)\partial_{x_m}\in \bigoplus_{j\in\N} E_{j}, $$
let $u_j=\PP_j u$. On $E_{j}={\rm Vect}\left(\cos\la{\bf k},\bullet\ra\partial_{x_m} , \sin\la{\bf k},\bullet\ra\partial_{x_m}\right)_{{\bf k}\in\Z^d, |{\bf k}|=j \atop{1\leq m\leq d}}$, we have
\begin{eqnarray*}
((d_0\circ d_0^*)u_{j})(x)&=& \sum_{{\bf k}\in\Z^d\atop{|{\bf k}|=j}}\sum_{m=1}^d a_{\bf k} \left(2\cos\la {\bf k},x \ra-\cos\la {\bf k},x+\alpha \ra-\cos\la {\bf k},x-\alpha \ra\right) \partial_{x_m} \\
& & + \, \sum_{{\bf k}\in\Z^d\atop{|{\bf k}|=j}}\sum_{m=1}^d  b_{\bf k} \left(2\sin\la {\bf k},x \ra-\sin\la {\bf k},x+\alpha \ra-\sin\la {\bf k},x-\alpha \ra\right)\partial_{x_m},\\
((d_1^*\circ d_1)u_{j})(x)&=& n \sum_{{\bf k}\in\Z^d\atop{|{\bf k}|=j}}\sum_{m=1}^d  a_{\bf k} \sum_{l=0}^{n-1}\cos\la {\bf k},x+l\alpha \ra \partial_{x_m}
 +    n\sum_{{\bf k}\in\Z^d\atop{|{\bf k}|=j}} \sum_{m=1}^d b_{\bf k}\sum_{l=0}^{n-1} \sin\la {\bf k},x+l\alpha \ra \partial_{x_m}.
\end{eqnarray*}
For ${\bf k}=({\bf k}_1,\cdots, {\bf k}_d)$,
if $\la {\bf k},\alpha\ra={\bf k}_1\alpha_1+\cdots +{\bf k}_d\alpha_d\in 2\pi\Z $, which means that ${\bf k}_1\in n_1 \Z$, $\cdots$, ${\bf k}_d\in n_d \Z$, we have
\begin{eqnarray}
2\cos\la {\bf k},x \ra-\cos\la {\bf k},x+\alpha \ra-\cos\la {\bf k},x-\alpha \ra&=&0,\label{d0res_cos}\\
2\sin\la {\bf k},x \ra-\sin\la {\bf k},x+\alpha \ra-\sin\la {\bf k},x-\alpha \ra&=&0,\label{d0res_sin}
\end{eqnarray}
\begin{equation}
\sum_{l=0}^{n-1}\cos\la {\bf k},x+l\alpha \ra =n\cos\la {\bf k},x\ra,\quad \sum_{l=0}^{n-1}\sin\la {\bf k},x+l\alpha \ra =n\sin\la {\bf k},x\ra.\label{d1nonres_cossin}
\end{equation}
If $\la {\bf k},\alpha\ra={\bf k}_1\alpha_1+\cdots +{\bf k}_d\alpha_d\not\in 2\pi\Z $, which means that at least one ${\bf k}_l\not\in n_l \Z$, then
\begin{eqnarray}
2\cos\la {\bf k},x \ra-\cos\la {\bf k},x+\alpha \ra-\cos\la {\bf k},x-\alpha \ra
&=&4\sin^2\frac{\la {\bf k},\alpha \ra}{2}\cos\la {\bf k},x \ra,\label{d0nonres_cos}\\
2\sin\la {\bf k},x \ra-\sin\la {\bf k},x+\alpha \ra-\sin\la {\bf k},x-\alpha \ra
&=&4\sin^2\frac{\la {\bf k},\alpha \ra}{2}\sin\la {\bf k},x \ra,\label{d0nonres_sin}
\end{eqnarray}
\begin{eqnarray}
\sum_{l=0}^{n-1}\cos\la {\bf k},x+l\alpha \ra \ = \ \frac{\sin(\frac{n\la {\bf k},\alpha\ra}{2})}{\sin(\frac{\la {\bf k},\alpha\ra}{2})}\cos\left\la {\bf k},x +\frac{n-1}{2} \alpha\right\ra &=& 0,\\
\sum_{l=0}^{n-1}\sin\la {\bf k},x+l\alpha \ra  \ = \  \frac{\sin(\frac{n\la {\bf k},\alpha\ra}{2})}{\sin(\frac{\la {\bf k},\alpha\ra}{2})}\sin\left\la {\bf k},x +\frac{n-1}{2} \alpha\right\ra &=& 0.
\end{eqnarray}
In view of (\ref{d0res_cos}) and (\ref{d0res_sin}), the generator of $G$ is not Diophantine for the action $\pi$ on $\T^d$ in the sense of Definition \ref{dioph-dolg}.

Let us show that $\pi$ is $\square-$Diophantine as in Definition \ref{defi-dio-intro}.
Applying Lemma \ref{lemma-cyclic}, we have $\Ker \square=0$, then the subspace $\IM(d_0\circ d_0^*)_j=\Ker(d_1^*\circ d_1)_j$ is
$${\rm Vect}\left\{\cos\la {\bf k},x \ra\partial_{x_m}, \sin\la {\bf k},x \ra \partial_{x_m}:\begin{array}{l}
     {\bf k}\in\Z^d, \ |{\bf k}|=j, \ \la {\bf k},\alpha\ra\not\in 2\pi\Z\\[1mm]
      1\leq m\leq d
    \end{array}
\right\}.$$
According to (\ref{d0nonres_cos}) and (\ref{d0nonres_sin}), the eigenvalues of $(d_0\circ d_0^*)_j$ are
$$4\sin^2\frac{\la {\bf k},\alpha \ra}{2},\quad |{\bf k}|=j, \quad \la {\bf k},\alpha\ra\not\in 2\pi\Z.$$
On the other hand, the subspace $\Ker(d_0\circ d_0^*)_j=\IM(d_1^*\circ d_1)_j$ is
$${\rm Vect}\left\{\cos\la {\bf k},x \ra\partial_{x_m}, \sin\la {\bf k},x \ra \partial_{x_m}:\begin{array}{l}
     {\bf k}\in\Z^d, \ |{\bf k}|=j, \ \la {\bf k},\alpha\ra\in 2\pi\Z\\[1mm]
      1\leq m\leq d
    \end{array}
\right\}.$$
According to (\ref{d1nonres_cossin}), the only eigenvalue of $(d_1^*\circ d_1)_j$ is $n^2$.
Since there are at most $(\prod_{l=1}^d n_l-1)$ non-vanishing values in $\left\{4\sin^2\frac{\la {\bf k},\alpha \ra}{2}\right\}_{{\bf k}\in\Z^d}$, the Box operator $\square=d_0\circ d_0^*+d_1^*\circ d_1$ has at most $\prod_{l=1}^d n_l$ eigenvalues. Hence $\pi$ is $\square-$Diophantine.
\qed


\smallskip

According to Theorem \ref{thm-geom0}, the $n-$periodic $F\in {\rm Diff}^\infty(\T^d)$ (resp. ${\rm Diff}^\omega(\T^d)$) defined in (\ref{F_defi}) is smoothly (resp. analytically) conjugated to the $n-$periodic translation $\pi$.

\subsection{Some other groups}\label{sec_groups}

As a corollary of Theorem \ref{thmmain}, we have
\begin{cor}\label{corKazhdan}
Let $G$ be a discrete group with Kazhdan's property (T).
Let $\pi$ be $G-$action by smooth (resp. analytic) isometries on the smooth (resp. analytic) compact Riemannian manifold $M$.
Then any smooth (resp. analytic) perturbation $\pi_{P_0}$ of $\pi$ is smoothly (resp. analytically) conjugate to $\pi$.
\end{cor}
 \proof According \cite{BHV}[Remark 1.1.4], the property (T) of the finitely presented group $G$ implies that the set of generators $\CS$ is Diophantine in the sense of Definition \ref{dioph-dolg}, and condition (\ref{dioph-cond}) is satisfied with $\tau=0$.
Then, through Lemma \ref{lemma_d_0-dio} and Remark \ref{rmk_DioDol}, this $G-$action $\pi$ is $d_0-$Diophantine.

On the other hand, according to \cite{fisher-margulis-invent}[Proposition 6.1], any $G-$action by smooth diffeomorphisms that is $C^K$ close to $\pi$ for a given $K\geq 2$, is $C^\ell$ conjugate to $\pi$ for any $\ell\geq K$ through a sequence of $C^\infty$ transformations $(\psi_n)\subset {\rm Diff}^\infty(M)$ which are $C^{K-1}$ near-identity. In this sense, $\pi_{P_0}$ is conjugate to $\pi$ through a $C^0$ near-identity $C^{\widehat R}$ transformation.
Indeed, as in several smooth rigidity theorems, e.g., \cite{nikos-tore, Moser-cercle, petko-tore}, stating that $\pi_{P_0}$ is a smooth perturbation of $\pi$ means that $\pi_{P_0}$ is $C^R$ close to $\pi$ for some $R>2$, depending on $n$ and $\tau$ (see Section \ref{secKAMsmooth} for more details). On the other hand, declaring that $\pi_{P_0}$ is an analytic perturbation of $\pi$ means that $\pi_{P_0}$ is close to $\pi$ in a complexified neighborhood of $M$ (known as a Grauert tube, see Section \ref{sec_Grauert_Hardy}), which also implies that $\pi_{P_0}$ is $C^2$ close to $\pi$ on $M$.
Hence, $\pi_{P_0}$ is $C^\ell$ conjugate to $\pi$ for any $\ell\geq R$ through a $C^{0}$ near-identity transformation.
Applying Theorem \ref{thmimpressive}, the corollary is shown.\qed

\begin{remark}
Corollary \ref{corKazhdan} in the smooth context was previously shown by Fisher-Margulis \cite{fisher-margulis-invent}[Theorem 1.3] (see also \cite{Fisher}[Theorem 1.2]). In this paper, we 
also establish the result in the analytic setting.
\end{remark}

Studies of Kazhdan's property (T) contribute to strong rigidity results in geometry and group theory. We refer the reader to \cite{BHV} for the basic theory and notable examples.

\smallskip

As a corollary of Theorem \ref{thm-geom0}, we have
\begin{cor}
Let $\Gamma$ be an irreducible lattice in a semi-simple Lie group with rank at least $2$. Let $\pi$ be an action of $\Gamma$ by smooth (resp. analytic) isometries on a compact smooth (resp. analytic) Riemannian manifold $M$ having Diophantine relations in the sense of \rd{defi-dio-relation}. Then any smooth (resp. analytic) perturbation $\pi_{P_0}$ of $\pi$ is smoothly (resp. analytically) conjugate to $\pi$.
\end{cor}
\proof According to \cite{margulis-book}[Introduction, Theorem 3], for $j\in \N$, $H^1(\Gamma, E_{\lambda_j})=0$. Moreover, $\overline{ \langle \Gamma\rangle}$ is a semi-simple Lie group. As established by Dolgopyat \cite{Dolgo2002}[Theorem A.3], $\Gamma$ is Diophantine in the sense of Definition \ref{dioph-dolg}. Since $\pi$ also has Diophantine relations, it is $\square-$Diophantine.
We can therefore apply Theorem \ref{thm-geom0}.
\qed

\subsection{Rotations on the sphere $\SS^2$}\label{sec_sphere}

In spherical coordinates $(\theta, \phi)$, $0\leq \theta\leq \pi $, $0\leq \phi< 2\pi $, the eigenvectors associated with the eigenvalue $\sqrt{J(J+1)}$, $J\in \N$, of $|\Delta_{T\SS^2}|^\frac12$, are the {\it vector spherical harmonics} ${\bf Y}_{J m}^L(\theta,\phi)$ for $L=J,J\pm 1$ (with the only exception that $L = 1$ for $J = 0$), and $m=-J,\cdots, 0,\cdots, J$,
According to \cite{VMK-AngularMomentum}[Chapter 7, Section 7.3], the vector spherical harmonics constitute a complete orthonormal basis of $L^2(\SS^2, T\SS^2,d\Omega)$: for any vector field $u$ such that
$$\int_{\SS^2}|u(\theta,\phi)|^2 d\Omega<\infty,\qquad d\Omega:=\sin\theta \, d\theta \,  d\phi, $$
we have that
$u(\theta,\phi)=\sum_{J,L,m} u_{J m}^L {\bf Y}_{J m}^L(\theta,\phi)$.

Under a general rotation ${\rm R}:(\theta,\phi)\mapsto (\theta',\phi')$ on $\SS^2$, characterized by a ${\rm SO}(3)$ matrix, the vector spherical harmonics are transformed as
$${\bf Y}_{J m'}^L(\theta',\phi')= \sum_{m} D_{mm'}^J({\rm R}){\bf Y}_{J m}^L(\theta,\phi),$$
where $D_{mm'}^J({\rm R})$ is the complex conjugate of the element of Wigner D-matrix $D^J({\rm R})$ (see \cite{VMK-AngularMomentum}[Chapter 4]).
Hence, $u\in L^2(\SS^2, T\SS^2,d\Omega)$ is transformed into
$$u(\theta',\phi')=\sum_{J,L,m'} u_{J m'}^L {\bf Y}_{J m'}^L(\theta',\phi')= \sum_{J,L,m}\left(\sum_{m'}D_{mm'}^J({\rm R}) u_{J m'}^L\right){\bf Y}_{J m}^L(\theta,\phi).$$

Given rotations ${\rm R}_1$, $\cdots$, ${\rm R}_k$ on $\SS^2$, regarded as the generating isometries $(\pi(\gamma_l))_{1\leq l\leq k}$, the $d_0-$Diophantine condition in the sense of Definition \ref{defi-dio-gnrt} relies on the arithmetric properties of eigenvalues of
the unitary matrices $D^J({\rm R}_l)$. As a simple case,
for the rotations turning around the $z-$axis by angle $\alpha_l$ (hence commuting), their Wigner D-matrices are diagonal:
\begin{equation}\label{WignerD_diagonal}
D^J({\rm R}_l)={\rm Diag}\{e^{{\rm i} m \alpha_l}: \ -J\leq m\leq J \}, \quad l=1,\cdots,k.
\end{equation}
If $(\alpha_l)_{1\leq l\leq k}$ is simultaneous Diophantine \cite{Moser-cercle} in the sense that,
there exist $\sigma,\tau >0$ such that $\max_l |j \alpha_l|\geq \sigma |j|^{-\tau}$ for $j\in \Z^*$, then, for $u\in \IM(d_0^*\circ \PP_J)$, $J\in \N^*$, we have
$$\|d_0 u\|_{L^2} =\left(\sum_{l=1}^k \sum_{L=J,J\pm1}\sum_{-J\leq m\leq J \atop{m\neq 0}}|(e^{{\rm i} m \alpha_l}-1)u_{Jm}^L|^2\right)^\frac12\geq \frac{\sigma}{(J(J+1))^{\frac\tau2}}\|u\|_{L^2}.$$
The above inequality is indeed (\ref{d0-Dio-condition}) in the proof of Lemma \ref{lemma_d_0-dio}.
Hence these generating rotations satisfy the $(\sigma^2, 2\tau)-d_0-$Diophantine.
According to Theorem \ref{thmmain}, there exists $\widehat R>0$ such that any smooth or analytic perturbations which are simultaneously $C^{\widehat R}$ almost conjugate to $({\rm R}_l)_{1\leq l\leq k}$, are simultaneously smoothly or analytically conjugate to these rotations.

In view of (\ref{WignerD_diagonal}), we see that $\dim\Ker d_0=\infty$ in this simple situation, since $1$ is the eigenvalue of Wigner D-matrix $D^J({\rm R}_l)$ for every $J\in \N$. Hence, the Diophantine condition in Definition \ref{dioph-dolg} is not satisfied.

%
%



\section{Smooth KAM scheme}\label{secKAMsmooth}

Let $M$ be a smooth compact Riemannian manifold of dimension $n$ and $G$ a finitely presented group with $k=\# \CS$, $p=\# \CR$ as given in Section \ref{sec_GroupAction}.
Let $\pi$ be a $G-$action by smooth isometries on $M$. Let $\pi_0$ be a $G-$action with $\pi_{P_0}(\gamma)=\Exp\{P_0(\gamma)\}\circ\pi(\gamma)$ for $\gamma\in \CS$, where $P_0:\CS\to \Gamma^\infty$, with $\|P_0\|_{\CS,C^0}=\varepsilon_0$.

In Theorem \ref{thmmain}, we assume that, for fixed $(\sigma,\tau)$, $\pi$ is $(\sigma,\tau)-d_0-$Diophantine and $\pi_{P_0}$ is $\varepsilon_0^\frac34-C^{\widehat R}$ almost conjugate to $\pi$, where $\widehat R:=20(\tau+n+1)$ \footnote{We assume that the constant $\tau\geq 0$ in (\ref{diophantine-gnrt-intro}) and (\ref{dio-box}) satisfies $20\tau\in\N$. Otherwise, with $\frac1{20}\lfloor 20\tau\rfloor+1$ in the place of $\tau$, the inequalities (\ref{diophantine-gnrt-intro}) and (\ref{dio-box}) still hold.}
. In Theorem \ref{thm-geom0}, we assume that $\pi$ is $(\sigma,\tau)-\square-$Diophantine and the first cohomology group of the complex \re{complex} is vanishing, which implies, according to \rl{lem-H1}, that $ \Ker\square=0$.

Let us set $R_*:=60(\tau+n+1)$.
Let $\pi_{P_{m}}:G\to {\rm Diff}^\infty(M)$ be a $G-$action by diffeomorphisms of $M$ of the form
$$\pi_{P_{m}}(\gamma)=\Exp\{P_m(\gamma)\}\circ \pi(\gamma),\qquad \gamma\in \CS, $$
where $P_m:\CS\to \Ga^\infty(M,TM)$, and for $\CP_m:=(P_m(\gamma))_{\gamma\in \CS}\in {\bf G}_\pi$,
\begin{equation}\label{esti-Pm}
\|\CP_m\|_{C^0}=\|P_m\|_{\CS,C^0}<\varepsilon_m:=\varepsilon_0^{(\frac54)^m},\qquad  \|\CP_m\|_{C^{R_*}}<\varepsilon_m^{-1}.
\end{equation}
Corresponding to the hypothesis of Theorem \ref{thmmain} and \ref{thm-geom0} respectively, we consider $\pi_{P_m}$ in the one of two following situations: either
\begin{hypo}\label{hypo_d0}
  $\pi$ is $d_0-$Diophantine, and $\pi_{P_m}$ is $\varepsilon_m^\frac34-C^{\widehat R}$ almost conjugate to $\pi$,
\end{hypo}
or
\begin{hypo}\label{hypo_Box}
 $\pi$ is $\square-$Diophantine and the first cohomology group is vanishing.
\end{hypo}
\noindent In view of Remark \ref{rmk_Diobox-equi}, $\pi$ is $d_0-$Diophantine in both cases.

The goal in this $(m+1)-$th KAM step is to conjugate the $\pi_{P_m}(\gamma)$'ssimultaneously to
$\pi_{P_{m+1}}(\gamma)$ for $\gamma \in \CS$, with $\CP_{m+1}=(P_{m+1}(\gamma))_{\gamma\in \CS}\in {\bf G}_\pi$, so that
\beq \label{goal-ind}
\|\CP_{m+1}\|_{C^0}<\varepsilon_{m+1},\qquad  \|\CP_{m+1}\|_{C^{R_*}}< \varepsilon_{m+1}^{-1}.
\eeq
Moreover, we show that $\pi_{P_{m+1}}$ is $\varepsilon_{m+1}^\frac34-C^{\widehat R}$ almost conjugate to $\pi$
under {\bf H \ref{hypo_d0}}.

As in previous sections, the inequality with ``$\lesssim$" means boundedness from above by an implicit constant depending on the manifold $M$, the group $G$, and the Diophantine constants $\sigma$, $\tau$. The inequality ``$\lesssim_R$" means that the later constant depends also on the order $R\in\N$. 
Let $c_\lesssim >1$ be the maximum of these implicit constants for $R\leq R_*=60(\tau+n+1)$ (if depending on $R$). Assume that $\varepsilon_0$ is sufficiently small such that
\begin{equation}\label{smallvarepsilon0}
\left(60+\tau+n+ k+ p+ c_{\lesssim} \right)^{12+9n}<\varepsilon_0^{-\frac{1}{60(\tau+6n+1)}}.
\end{equation}

\subsection{Truncation operator} As a first procedure of one KAM step, we define the truncation operator on the $L^2$ space and based on the decomposition (\ref{generalFourier}).
\begin{defi}\label{def_truncation}
Given $\nu , N\in\N^*$, let us define the {\bf truncation operator} of degree N, ${\Pi}^{(\nu)}_{N}: L^2(M,TM)^\nu\to L^2(M,TM)^\nu$ as
$$
{\Pi}^{(\nu)}_{N}u:=\sum_{j\in \N\atop{\lambda_j\leq N}}\PP_j u=\left(\sum_{j\in \N\atop{\lambda_j\leq N}}\sum_{i\in I_j} \hat u_{l,i} {\bf e}_i\right)_{1\leq l\leq \nu}, \quad u=\left(\sum_{i\in\N} \hat u_{l,i}{\bf e}_i\right)_{1\leq l\leq \nu} \in L^2(M,TM)^\nu,
$$
and let ${\Pi}^{(\nu)\perp}_{N}:=\id-{\Pi}^{(\nu)}_{N}=\sum_{j\in \N\atop{\lambda_j> N}}\PP_j$.
\end{defi}

\noindent For convenience, the superscript ``$(\nu)$" in the notation of the truncation operator will be omitted in the sequel if there is no ambiguity.

Let $N_m:=\varepsilon_m^{-\frac{1}{8(\tau+n+1)}}$. In view of the fact that $\frac{\frac{3}{2}n+1}{8(\tau+n+1)}< \frac{3}{16}$, we have that
\beq\label{Nm3_16}
N_m^{\frac{3}{2}n+1}<\varepsilon_m^{-\frac{3}{16}}.
\eeq
For $\CP_m=(P_m(\gamma))_{\gamma\in \CS}\in \Ga^\infty(M,TM)^k$ satisfying (\ref{esti-Pm}), we have

\begin{lemma}\label{lemma_reste-smooth}
We have $\left\|{\Pi}^{\perp}_{N_m}\CP_m\right\|_{C^0}\lesssim \varepsilon_m^{\frac53} $, and for $R\in\N$, $$\left\|{\Pi}_{N_m}\CP_m\right\|_{C^R}, \ \left\|{\Pi}^{\perp}_{N_m}\CP_m\right\|_{C^R}\lesssim_R \varepsilon_m^{-\frac{3}{16}} \left\|\CP_m\right\|_{C^R}.$$
In particular, $\left\|{\Pi}_{N_m}\CP_m\right\|_{C^{R_*}}, \left\|{\Pi}^{\perp}_{N_m}\CP_m\right\|_{C^{R_*}}\lesssim \varepsilon_m^{-\frac{19}{16}}$.
\end{lemma}
\proof  According to Proposition \ref{propNorm-CrHr}, (\ref{esti-Pm}) implies that
\begin{equation}\label{esti-Pm-Hr*}
\|\CP_m\|_{L^2}=\|\CP_m\|_{\CH^0}\lesssim \varepsilon_m,\qquad  \|\CP_m\|_{\CH^{R_*}}\lesssim \varepsilon_m^{-1}.\end{equation}
Let $d=10(\tau+n+1)$. By the definition of Sobolev norm, together with (\ref{asymp-lamdak}), we have
\begin{eqnarray*}
\left\|{\Pi}^{\perp}_{N_m}\CP_m\right\|_{\CH^{\frac{3}{2}n+1}}&\lesssim&\left(\sum_{1\leq l\leq k} \sum_{j\in \N\atop{\lambda_j> N_m}}(1+\lambda_j)^{3n+2}\sum_{i\in I_j} |\hat P_{m,i}(\gamma_l)|^2\right)^\frac12\\
&\lesssim& \left\|\CP_m\right\|_{\CH^d} \left(\sum_{K= N_m}^\infty (K+1)^n \sum_{K< \lambda_j\leq K+1}(1+\lambda_j)^{-2d+3n+2}\right)^\frac12 \\
&\lesssim& \left\|\CP_m\right\|_{\CH^d} \left( \int_{N_m}^{+\infty} t^{4n+2-2d} \, dt\right)^\frac12 .
\end{eqnarray*}
By computations, we have
$$\int_{N_m}^{+\infty} t^{4n+2-2d} \, dt=\int_{N_m}^{+\infty}t^{-20\tau-16n-18} \, dt
=\frac{1}{(20\tau+16n+17) N_m^{20\tau+16n+17}}   \   \lesssim  \   \varepsilon_m^{2}. $$
On the other hand, applying the interpolation lemma \ref{interpol-hs} with (\ref{esti-Pm-Hr*}), we obtain
$$
\left\|\CP_m\right\|_{\CH^d}\leq \left\|\CP_m\right\|_{L^2}^{\frac{R_*-d}{R_*}} \left\|\CP_m\right\|_{\CH^{R_*}}^{\frac{d}{R_*}} \lesssim \varepsilon_m^{1-\frac{1}{6}}\varepsilon_m^{-\frac{1}{6}}=\varepsilon_m^{\frac23}.$$
Hence, according to \rp{propNorm-CrHr},
$\left\|{\Pi}^{\perp}_{N_m}\CP_m\right\|_{C^0}\lesssim \left\|{\Pi}^{\perp}_{N_m}\CP_m\right\|_{\CH^{\frac{3}{2}n+1}}\lesssim \varepsilon_m^{\frac53}$.

For any $R\in\N$, using Corollary \ref{cor-CR-truncation} and recalling (\ref{Nm3_16}), we have
\begin{equation}\label{smoothing-truncation}
\left\|{\Pi}^{\perp}_{N_m}\CP_m\right\|_{C^{R}}, \ \left\|{\Pi}_{N_m}\CP_m\right\|_{C^{R}} \lesssim_R N_m^{\frac32n +1}\left\|\CP_m\right\|_{C^{R}} < \varepsilon_m^{-\frac{3}{16}}\left\|\CP_m\right\|_{C^{R}}.
\end{equation}
In particular, according to (\ref{esti-Pm}), we have
$$\left\|{\Pi}_{N_m}\CP_m\right\|_{C^{R_*}}, \left\|{\Pi}^{\perp}_{N_m}\CP_m\right\|_{C^{R_*}} \lesssim\varepsilon_m^{-\frac{3}{16}} \left\|\CP_m\right\|_{C^{R_*}}\leq  \varepsilon_m^{-\frac{19}{16}}.    \qed$$

\subsection{Construction of conjugacy}\label{sec_cohomoeq-smooth}

Let us write the equation of the $(m+1)-{\rm th}$ conjugacy for every $\gamma\in \CS$, i.e., the conjugacy from $\pi_{P_m}(\gamma)=\Exp\{P_m(\gamma)\}\circ \pi(\gamma)$ to $\pi_{P_{m+1}}(\gamma)=\Exp\{P_{m+1}(\gamma)\}\circ \pi(\gamma)$, by $\Exp\{w_m\}\in {\rm Diff}^\infty(M)$, as the following:
\begin{equation}\label{applem-3-smooth}
\Exp\{w_m\}^{-1}\circ\Exp\{P_m(\gamma)\}\circ \pi(\gamma) \circ \Exp\{w_m\} = \Exp\{P_{m+1}(\gamma)\}\circ \pi(\gamma),\quad \forall \ \gamma\in \CS,
\end{equation}
with $w_m\in \Ga^\infty$
and $P_{m+1}:\CS\to \Ga^\infty$, with
$\CP_{m+1}=(P_{m+1}(\gamma))_{\gamma\in \CS}$  satisfying (\ref{goal-ind}).

Let us rewrite this conjugacy equation in a more tractable way. To do so, let us assume that $\|w_m\|_{C^1}$ is small enough so that, according to Lemma \ref{lem_exp-isometry-smooth}, Eq. \re{applem-3-smooth} reads
\begin{eqnarray*}
\Exp\{P_{m+1}(\gamma)\}&=& \Exp\{w_m\}^{-1} \circ \Exp\{P_m(\gamma)\}\circ \pi(\gamma) \circ \Exp\{w_m\}\circ \pi^{-1}(\gamma)\\
&=& \Exp\{w_m\}^{-1} \circ \Exp\{P_m(\gamma)\} \circ  \Exp\left\{\pi(\gamma)_*w_m\right\},\quad \forall \ \gamma\in \CS,
\end{eqnarray*}
and hence
\begin{equation}\label{applem-1-smooth}
\Exp\{P_{m}(\gamma)\}\circ \Exp\left\{\pi(\gamma)_*w_m\right\}= \Exp\{w_m\}\circ \Exp\{P_{m+1}(\gamma)\}.
\end{equation}
By (\ref{esti-Pm}) and (\ref{goal-ind}), we see that, through Lemma \ref{interpol-cs}, $\|P_m\|_{\CS,C^1}$ and $\|P_{m+1}\|_{\CS,C^0}$ would be sufficiently small so that we could apply Lemma \ref{lem_moser-s1} to both sides of (\ref{applem-1-smooth}), and obtain
\begin{eqnarray}
& & P_m(\gamma)+\pi(\gamma)_*w_m+s_1\left(P_{m}(\gamma), \pi(\gamma)_*w_m\right) \label{eq_cohomo-im-smooth}\\
   &= & w_m+P_{m+1}(\gamma)+ s_1\left(w_m, P_{m+1}(\gamma)\right), \quad \forall \ \gamma\in \CS,  \nonumber
\end{eqnarray}
where, for every $\gamma\in \CS$,
$s_1\left(P_{m}(\gamma), \pi(\gamma)_*w_m\right)$, $s_1\left(w_m, P_{m+1}(\gamma)\right)\in \Gamma^\infty(M,TM)$,
and for $1\leq l\leq k$, we would have, through Proposition \ref{prop_s1-smooth},
\begin{eqnarray}
\|s_1(w_m, P_{m+1}(\gamma_l))\|_{C^0}&\lesssim& \|w_m\|_{C^1}\|P_{m+1}(\gamma_l)\|_{C^0}, \label{error_s1-1}\\
\|s_1(P_{m}(\gamma_l),\pi(\gamma_l)_*w_m)\|_{C^0} &\lesssim& \|P_{m}(\gamma_l)\|_{C^1}\|w_m\|_{C^0}. \label{error_s1-2}
\end{eqnarray}


Recalling (\ref{defi_d0}),
Eq. \re{eq_cohomo-im-smooth} with $\gamma=\gamma_1,\cdots,\gamma_k\in \CS$ can be written as the {\bf cohomological equation} on $M$:
\begin{equation}\label{cohomo}
\CP_m-d_0 w_m=\CP_{m+1}+s_1(w_m, \CP_{m+1}) -s_1(\CP_{m}, \pi_*w_m),
\end{equation}
 with the last two terms defined as
 $$s_1(w_m, \CP_{m+1}):= \left(s_1(w_m, P_{m+1}(\gamma_l))\right)_{l},\quad  s_1(\CP_{m}, \pi_*w_m ):= \left(s_1(P_{m}(\gamma_l), \pi(\gamma_l)_*w_m )\right)_{l}.$$
Our goal is to find $w_m\in \Ga^\infty$ and $\CP_{m+1}\in (\Ga^\infty)^k$ satisfying this equation as well as the aforementioned estimates, then we also solve Eq. \re{applem-3-smooth}.

With the decomposition (\ref{L2k-decomp}) of $L^2(M,TM)^k$, we have the following lemma, which allows us to solve Eq. (\ref{cohomo}) approximately.

\begin{lemma}\label{lem_sol_cohomo-d0}
If $\pi$ is $(\sigma,\tau)-d_0-$Diophantine, then, for any $u\in (\Ga^\infty)^k$, there exists a unique $w\in \IM d_0^*\cap \bigoplus_{\lam_j\leq N_m} E_{\lam_j}$ with $d_0 w=(\D_0\circ\Pi_{N_m}) u$ and
\begin{eqnarray}
\|\PP_j w\|_{L^2}&\lesssim&(1+\lam_j)^\tau \|\PP_j u\|_{L^2},\qquad \lam_j\leq N_m,\label{cohomo_esti-L2proj}\\
\|w\|_{C^R}&\lesssim_R& \varepsilon_m^{-\frac{3}{16}} \|u\|_{C^{R}},\qquad R\in \N.\label{cohomo_esti-smooth}
\end{eqnarray}
\end{lemma}

\proof For any $j\in\N$ with $\lam_j\leq N_m$, let us find some $f_j\in\IM d_0\cap E_{\lam_j}^k$ such that $$(d_0\circ d_0^*)f_j=u_j:=(\D_0\circ\PP_j) u.$$
Noting that $\IM d_0=\IM (d_0\circ d_0^*)$, we see the existence of such $f_j$ according to the $d_0-$Diophantine condition (\ref{diophantine-gnrt-intro}), and $\|f_j\|_{L^2}\leq \sigma^{-1}(1+\lam_j)^\tau\|u_j\|_{L^2}$.
With $w_j:=d_0^* f_j$, we obtain that $d_0 w_j=u_j$. In view of (\ref{d0*}) and (\ref{invar-d0}), we see that $w_j\in \IM d_0^*\cap E_{\lam_j}$
and $$\|w_j\|_{L^2}\lesssim (1+\lam_j)^\tau \|u_j\|_{L^2}\lesssim (1+\lam_j)^\tau \|\PP_j u\|_{L^2}.$$
Let $w:=\sum_{\lam_j\leq N_m} w_j$, then $d_0 w=(\D_0\circ \Pi_{N_m}) u$ and (\ref{cohomo_esti-L2proj}) is shown.

For $R\in \N$, according to Proposition \ref{propNorm-CrHr} and recalling (\ref{Nm3_16}), we have,
\begin{eqnarray*}
\|w\|_{C^R} \ \lesssim_R \ \|w\|_{\CH^{R+\frac32n+1}}&=&\left(\sum_{\lam_j\leq N_m} (1+\lam_j)^{2(R+\frac32n+1)}\|w_j\|^2_{L^2} \right)^\frac12\\
   &\lesssim_R& N_m^{\tau+\frac32n+1} \left(\sum_{\lam_j\leq N_m}\|u_j\|^2_{\CH^{R}} \right)^\frac12\\
   &\lesssim_R& N_m^{\tau+\frac32n+1} \| u\|_{\CH^{R}} \ \lesssim_R \ \varepsilon_m^{-\frac{3}{16}} \|u\|_{C^{R}}. \qed
\end{eqnarray*}

\begin{remark}\label{rmk_w-analytic}
The above solution $w\in \IM d_0^*\cap \bigoplus_{\lam_j\leq N_m} E_{\lam_j}$ is $C^\infty$ smooth on $M$ since it has finitely many ``Fourier modes".
In the analytic context, the above construction yields an analytic $w$ on $M$.
\end{remark}

As in Lemma \ref{lem_sol_cohomo-d0}, let $w_m\in \bigoplus_{\lam_j\leq N_m}  E_{\lam_j}\cap\IM d_0^*$ be the unique solution with bound: 
\beq\label{esti_wCR}
d_0w_m=(\D_0\circ\Pi_{N_m})\CP_m,\quad \|w_m\|_{C^R}\lesssim_R  \varepsilon_m^{-\frac{3}{16}}\|\CP_m\|_{C^{R}},\qquad R\in\N,\eeq
and, in particular, recalling (\ref{esti-Pm}),
\begin{equation}\label{esti_wm-smooth}
\|w_m\|_{C^0}  \lesssim \varepsilon_m^\frac{13}{16},\qquad \|w_m\|_{C^{R_*}} \lesssim \varepsilon_m^{-\frac{19}{16}}.
\end{equation}
By interpolation in Lemma \ref{interpol-cs} between $0$ and $R_*=60(\tau+n+1)\geq120$, we have
\begin{eqnarray}
\|w_m\|_{C^1} \lesssim \|w_m\|^{1-\frac{1}{R_*}}_{C^0} \|w_m\|^{\frac{1}{R_*}}_{C^{R_*}}\lesssim  \varepsilon_m^{\frac{13\cdot 119 -19}{16\cdot 120}}\leq \varepsilon_m^{\frac34}, \label{esti_wm-C1}
\end{eqnarray}
According to (\ref{s1-C0-decomp}) and (\ref{s1-Cr}) in \rp{prop_s1-smooth}, we have
\beq\label{s1wP-C0}
\|s_1(w_m, \CP_{m+1})\|_{C^0} \lesssim \|w_m\|_{C^1}\|\CP_{m+1}\|_{C^0}   \lesssim \varepsilon_{m}^{\frac34}\|\CP_{m+1}\|_{C^0} , \eeq
Hence, by the fixed point theorem (see e.g., \cite{dieudonne1}[10.1.1]), we
obtain the existence $\CP_{m+1}\in (\Ga^\infty)^k$,
such that Eq. (\ref{cohomo}) holds. Furthermore, we have
\beq\label{s1wP-CR*}
 \|s_1(w_m,\CP_{m+1})\|_{C^{R^*}} \lesssim \|w_m\|_{C^{R^*}}+ \|w_m\|_{C^{1}}\|\CP_{m+1}\|_{C^{R^*}}
 \lesssim \varepsilon_{m}^{\frac34}\|\CP_{m+1}\|_{C^{R^*}}+ \varepsilon_{m}^{-\frac{19}{16}}.
\eeq
Therefore, with $\D_0^\perp:=\CH+ \D_1$ and \re{esti_wCR}, we have
\beq\label{Pm+1}
\CP_{m+1}= \Pi_{N_m}^\perp \CP_m  -s_1(w_m,\CP_{m+1})+s_1(\CP_m, \pi_* w_m)+(\Pi_{N_m}\circ\D_0^\perp)\CP_m.
\eeq

\begin{lemma}\label{lem_rough-estim-smooth}
$\|\CP_{m+1}\|_{C^{R_*}}\lesssim\varepsilon_{m}^{-\frac{19}{16}}$ and
$\|\CP_{m+1}\|_{C^0}\lesssim   \|(\Pi_{N_m}\circ\D_0^\perp)\CP_m\|_{C^0}+ \varepsilon_{m}^{\frac{21}{16}}\lesssim  \varepsilon_m^{\frac{13}{16}}$.
\end{lemma}
\proof
By the interpolation lemma \ref{interpol-cs}, we have
\beq\label{PmC1}
\|\CP_m\|_{C^{1}}\lesssim \|\CP_m\|^{1-\frac{1}{R_*}}_{C^0}\|\CP_m\|^{\frac{1}{R_*}}_{C^{R_*}} \lesssim \varepsilon_m^{1-\frac{2}{R_*}}\leq\varepsilon_m^{\frac34}.
\eeq
According to (\ref{s1-C0-decomp}), (\ref{s1-Cr}) in \rp{prop_s1-smooth}, together with (\ref{esti-Pm}), (\ref{esti_wm-smooth}) and (\ref{esti_wm-C1}), we have
$$
\|s_1(\CP_m, \pi_* w_m)\|_{C^0}\lesssim\|\CP_m\|_{C^1} \|w_m\|_{C^0}\lesssim\varepsilon_m^{\frac34}\cdot \varepsilon_m^{\frac{13}{16}}= \varepsilon_m^{\frac{25}{16}},
$$
$$\|s_1(\CP_m, \pi_* w_m)\|_{C^{R_*}} \lesssim \|\CP_m\|_{C^{R_*}} + \|\CP_m\|_{C^{1}}\|w_m\|_{C^{R_*}} \lesssim \varepsilon_{m}^{-1}.$$
In view of Corollary \ref{cor-CR-truncation} and (\ref{Nm3_16}), for $0\leq R\leq R_*$, $\|(\Pi_{N_m}\circ\D_0^\perp)\CP_m\|_{C^R}\lesssim \varepsilon_m^{-\frac{3}{16}} \|\CP_m\|_{C^{R}}$.
Hence, combining with (\ref{s1wP-C0}), (\ref{s1wP-CR*}) and
Lemma \ref{lemma_reste-smooth}, we obtain
\begin{eqnarray*}
 \|\CP_{m+1}\|_{C^0}&\lesssim&   \|(\Pi_{N_m}\circ\D_0^\perp)\CP_m\|_{C^0} + \varepsilon_{m}^{\frac{21}{16}} \ \lesssim \ \varepsilon_m^{\frac{13}{16}}, \\
\|\CP_{m+1}\|_{C^{R_*}}&\lesssim & \|(\Pi_{N_m}\circ\D_0^\perp)\CP_m\|_{C^{R_*}} + \varepsilon_{m}^{-\frac{19}{16}} \ \lesssim \ \varepsilon_{m}^{-\frac{19}{16}}.\qed
\end{eqnarray*}

\subsection{Refined $C^0$ estimate of $\CP_{m+1}$}\label{sec_refine-smooth}

Recalling that the aim of this $(m+1)^{\rm th}$ KAM step is to conjugate the $G-$action $\pi_m(\gamma)$ to $\pi_{P_{m+1}}(\gamma)=\Exp\{P_{m+1}(\gamma)\}\circ\pi(\gamma)$ with (\ref{goal-ind}) satisfied
for $\CP_{m+1}=(P_{m+1}(\gamma))_{\gamma\in \CS}$.

In view of (\ref{smallvarepsilon0}), all implicit coefficients in the inequalities with ``$\lesssim$" are smaller than $\varepsilon_0^{-\frac{1}{16}}\leq \varepsilon_m^{-\frac{1}{16}}$.
Then the $C^{R_*}$ estimate in Lemma \ref{lem_rough-estim-smooth} implies $\|\CP_{m+1}\|_{C^{R_*}}\leq \varepsilon_{m}^{-\frac54}= \varepsilon_{m+1}^{-1}$. Let us refine the $C^0$ estimate of $\CP_{m+1}$, or more precisely that of $(\Pi_{N_m}\circ\D_0^\perp)\CP_m$,
under either the assumption {\bf H \ref{hypo_d0}} or {\bf H \ref{hypo_Box}}.
 In both situations, we will show that
\beq\label{PNDperpP}\|(\Pi_{N_m}\circ\D_0^\perp)\CP_m\|_{C^0}\lesssim \varepsilon_m^{\frac{21}{16}}, \eeq
which refines the $C^0$ estimate in Lemma \ref{lem_rough-estim-smooth} as we finally show that $\|\CP_{m+1}\|_{C^0}\lesssim \varepsilon_m^{\frac{21}{16}}$.
Then (\ref{goal-ind}) is proved.

\subsubsection{Vanishing first cohomology and $\square-$Diophantineness}

Under {\bf H \ref{hypo_Box}}, we have $\Ker\square=0$.
Then $\D_0^\perp \CP_m =\D_1 \CP_m $ and $\square\circ \D_1=d_1^*\circ d_1$ is invertible on $\IM d_1^*$.
The $(\sigma,\tau)-\square-$Diophantine condition means that the eigenvalues of $(d_1^*\circ d_1)_j$ are bounded from below by $\sigma(1+\lam_j)^{-1}$. Noting that
$$(\Pi_{N_m}\circ\D_1)\CP_m= \sum_{\lam_j\leq N_m} (d_1^*\circ d_1)_j^{-1}(d_1^*\circ d_1)_j\CP_m,$$
we have, applying Proposition \ref{propNorm-CrHr} and Lemma \ref{lem-d1-quadra}, and recalling (\ref{Nm3_16}), (\ref{PmC1}),
\begin{eqnarray*}
\|(\Pi_{N_m}\circ\D_1)\CP_m\|_{C^0}
   &\leq&\left\|\sum_{\lam_j\leq N_m}(d_1^*\circ d_1)_j^{-1}(d_1^*\circ d_1)_j\CP_m\right\|_{\CH^{\frac32n+1}}\\
   &\lesssim& \left(\sum_{\lam_j\leq N_m} (1+\lam_j)^\tau \|(d_1^*\circ d_1)_j\CP_m\|_{\CH^{\frac32n+1}}^2\right)^\frac12\\
   &\lesssim& N_m^{\tau+\frac32n+1}\|(d_1^*\circ d_1)\CP_m\|_{L^2}
    \  \lesssim \ \varepsilon_m^{-\frac{3}{16}} \|\CP_m\|_{C^{1}} \|\CP_m\|_{C^{0}} \  \lesssim \   \varepsilon_{m}^{\frac{25}{16}}.
\end{eqnarray*}
Hence (\ref{PNDperpP}) is shown.

\subsubsection{Almost conjugacy and $d_0-$Diophantineness}
Under the hypothesis {\bf H \ref{hypo_d0}}, the $G-$action by diffeomorphisms $\pi_m=\Exp\{P_m\}\circ \pi$ is $\varepsilon_m^\frac12-C^{\widehat R}$ almost conjugate to $\pi$, with $\widehat R=20(\tau+n+1)$.
In view of Definition \ref{defi_almost_conj-intro}, we have, for any $0<\varepsilon < \varepsilon_m$, there exists $y_m^{\varepsilon}\in \Ga^{\widehat R}$ with $\|y_m^{\varepsilon}\|_{C^{0}}< \varepsilon_m^\frac34$ and
$\|y_m^{\varepsilon}\|_{C^{\widehat R}}< \varepsilon_m^{-\frac34}$ such that
$$
\Exp\{y_m^{\varepsilon}\}^{-1} \circ \Exp\{P_m(\gamma)\}\circ \pi(\gamma)\circ \Exp\{y_m^{\varepsilon}\}=\Exp\{z_m^{\varepsilon}(\gamma)\}\circ \pi(\gamma),\quad \gamma\in \CS,$$
where $z_m^{\varepsilon}:\CS\to  \Ga^{0}$ satisfies that
$\|z_m^{\varepsilon}\|_{\CS,C^0}<\varepsilon$. By interpolation, we have
\begin{equation}\label{esti_ym-C1}
\|y_m^{\varepsilon}\|_{C^{1}} \lesssim \|y_m^{\varepsilon}\|_{C^{0}}^{1-\frac{1}{\widehat R}}\|y_m^{\varepsilon}\|_{C^{\widehat R}}^{\frac{1}{\widehat R}}< \varepsilon_m^{\frac34\left(1- \frac{2}{20(\tau+n+1)}\right) }< \varepsilon_m^{\frac{1}{2}}.
\end{equation}
According to Lemma \ref{lemma_exp-inverse}, there exists $\tilde y_m^{\varepsilon}\in \Ga^{1}$ with $\|\tilde y_m^{\varepsilon}\|_{C^{1}}\lesssim \|y_m^{\varepsilon}\|_{C^{1}}< \varepsilon_m^\frac12$ such that
$$\Exp\{\tilde y_m^{\varepsilon}\}\circ \Exp\{P_m(\gamma)\}\circ \pi(\gamma)\circ \Exp\{\tilde y_m^{\varepsilon}\}^{-1}
= \Exp\{z_m^{\varepsilon}(\gamma)\}\circ \pi(\gamma).$$
With $\CZ_m^{\varepsilon}:=(z_m^{\varepsilon}(\gamma))_{\gamma\in \CS}$, $\pi_*\tilde y_m^{\varepsilon}:=(\pi(\gamma)_*\tilde y_m^{\varepsilon})_{\gamma\in \CS} $,
the above conjugation  implies that
 \beq\label{conj_almost}
\CP_m+d_0\tilde y_m^{\varepsilon}=  \CZ_m^{\varepsilon}-s_1(\tilde y_m^{\varepsilon},\CP_m)+s_1(\CZ_m^{\varepsilon},\pi_*\tilde y_m^{\varepsilon}).
\eeq
Projecting (\ref{conj_almost}) onto $(\IM d_0)^\perp=\IM d_1^* \bigoplus \Ker\square$, we have
$$
\D_0^\perp\CP_m=\D_0^\perp(\CZ_m^{\varepsilon}+s_1(\CZ_m^{\varepsilon},\pi_*\tilde y_m^{\varepsilon}))-\D_0^\perp s_1(\tilde y_m^{\varepsilon},\CP_m).$$
Assume that $\varepsilon< \varepsilon_m^2$.
In view of Corollary \ref{cor-CR-truncation}, we have, through (\ref{s1-C0-decomp}) in \rp{prop_s1-smooth},
\begin{eqnarray*}
\|(\Pi_{N_m}\circ\D_0^\perp)\CP_m\|_{C^0}&\leq&\|(\Pi_{N_m}\circ\D_0^\perp)( \CZ_m^{\varepsilon}+s_1(\CZ_m^{\varepsilon},\pi_*\tilde y_m^{\varepsilon}))\|_{C^0}\\
& & + \, \|(\Pi_{N_m}\circ\D_0^\perp) s_1(\tilde y_m^{\varepsilon},\CP_m)\|_{C^0}\\
&\lesssim& N_m^{\frac32n+1}(\|\CZ_m^{\varepsilon}\|_{C^0}+\|y_m^{\varepsilon}\|_{C^1}\|\CP_m\|_{C^0}) \ \lesssim \ \varepsilon_{m}^{-\frac{3}{16}}\left(\varepsilon_m^2 +\varepsilon_m^\frac32\right) \ \lesssim \ \varepsilon_m^{\frac{21}{16}},
\end{eqnarray*}
which shows (\ref{PNDperpP}).

\subsection{Almost conjugacy of $\pi_{P_{m+1}}$}\label{sec_verify_alconj}

In the previous subsection, under {\bf H \ref{hypo_d0}} or {\bf H \ref{hypo_Box}}, we performed one KAM step, which conjugated $\pi_{P_m}$, satisfying estimates (\ref{esti-Pm}), to $\pi_{P_{m+1}}$ with (\ref{goal-ind}) fulfilled.

Since {\bf H \ref{hypo_Box}} is imposed on the $G-$action by isometries $\pi$ and is independent of iteration step, we are able to continue the iteration with $\pi_{P_{m+1}}$.
Regarding {\bf H \ref{hypo_d0}}, because the almost conjugacy is assumed for $\pi_{P_m}$, we need to verify it for $\pi_{P_{m+1}}$ before proceeding to the next iteration.

\begin{prop}\label{prop-verify_alconj}
If $\pi_{P_{m}}$ is $\varepsilon_m^{\frac34}-C^{\widehat R}$ almost conjugate to $\pi$, then $\pi_{P_{m+1}}$ satisfying (\ref{applem-3-smooth}) is $\varepsilon_{m+1}^{\frac34}-C^{\widehat R}$ almost conjugate to $\pi$.
\end{prop}
\proof According to Definition \ref{defi_almost_conj-intro},
for any $\varepsilon>0$, there exists $y_{m}^{\varepsilon}\in \Ga^{\widehat R}$ such that $\|y_{m}^{\varepsilon}\|_{C^{0}}<\varepsilon_m^{\frac34}$ and $\|y_{m}^{\varepsilon}\|_{C^{\widehat R}}<\varepsilon_m^{-\frac34}$, and
$$ \Exp\{y_{m}^{\varepsilon}\}^{-1} \circ \Exp\{P_{m}(\gamma)\}\circ \pi(\gamma)\circ \Exp\{y_{m}^{\varepsilon}\}=\Exp\{z_{m}^{\varepsilon}(\gamma)\}\circ \pi(\gamma),\quad \gamma\in \CS,$$
with $z_{m}^{\varepsilon}:\CS\to  \Ga^1$ satisfies that $\|z_{m}^{\varepsilon}\|_{\CS,C^0}<\varepsilon$.
Combining with (\ref{applem-3-smooth}), we have
\begin{eqnarray*}
& &\Exp\{y_{m}^{\varepsilon}\}^{-1}\circ \Exp\{w_m\} \circ \Exp\{P_{m+1}(\gamma)\}\circ \pi(\gamma) \circ \Exp\{w_m\}^{-1}\circ \Exp\{y_{m}^{\varepsilon}\}\\
&=& \Exp\{z_{m}^{\varepsilon}(\gamma)\}\circ \pi(\gamma).
\end{eqnarray*}
By Lemma \ref{lemma_exp-inverse}, there exists $\tilde w_m\in \Ga^\infty$ with $\|\tilde w_m\|_{C^R}\lesssim\|w_m\|_{C^R} $ for $0\leq R\leq \widehat R$
 such that $\Exp\{\tilde w_m\}=\Exp\{w_m\}^{-1}$.
Then, for
$\breve{y}_{m+1}^{\varepsilon}:=\tilde w_m+y_{m}^{\varepsilon}+ s_1(\tilde w_m, y_{m}^{\varepsilon}) \in \Ga^{\widehat R}$ with $s_1(\cdot,\cdot)$ defined as in Lemma \ref{lem_moser-s1}, we have that
$$\Exp\{\tilde w_m\}\circ  \Exp\{y_{m}^{\varepsilon}\} =\Exp\{\breve{y}_{m+1}^{\varepsilon}\},$$
and hence, for every $\gamma\in \CS$,
\begin{equation}\label{alconj-breve}
\Exp\{\breve{y}_{m+1}^{\varepsilon}\}^{-1} \circ \Exp\{P_{m+1}(\gamma)\}\circ \pi(\gamma) \circ \Exp\{\breve{y}_{m+1}^{\varepsilon}\}=\Exp\{z_{m}^{\varepsilon}(\gamma)\}\circ \pi(\gamma).
\end{equation}

Applying the interpolation lemma \ref{interpol-cs} with (\ref{esti_wm-smooth}), we have
\beq\label{esti_wm-CR}
\|\tilde w_m\|_{C^{\widehat R}} \lesssim\|w_m\|_{C^{\widehat R}} \lesssim \|w_m\|^{1-\frac{\widehat R}{R_*}}_{C^0} \|w_m\|^{\frac{\widehat R}{R_*}}_{C^{R_*}}\lesssim  \varepsilon_m^{\frac{13\cdot 2 -19}{16\cdot 3}}\leq \varepsilon_m^{\frac{7}{48}}.\eeq
By (\ref{s1-C0-decomp}) and (\ref{s1-Cr}) in Proposition \ref{prop_s1-smooth} and (\ref{esti_wm-C1}), we have
\begin{eqnarray*}
& &\|s_1(\tilde w_m,y_{m}^{\varepsilon})\|_{C^0}\lesssim \|w_m\|_{C^1} \|y_{m}^{\varepsilon}\|_{C^0}\lesssim \varepsilon_m^{\frac34} \cdot\varepsilon_m^{\frac34}=\varepsilon_m^{\frac32},\\
& &\|s_1(\tilde w_m,y_{m}^{\varepsilon})\|_{C^{\widehat R}}\lesssim \|\tilde w_m\|_{C^{\widehat R}}+\|\tilde w_m\|_{C^1}\|y_{m}^{\varepsilon}\|_{C^{\widehat R}}\lesssim \varepsilon_m^{\frac{7}{48}}+\varepsilon_m^{\frac34}\cdot \varepsilon_m^{-\frac34} \lesssim 1.
\end{eqnarray*}
Hence, $\breve{y}_{m+1}^{\varepsilon}$ satisfies that
\begin{eqnarray}
\|\breve{y}_{m+1}^{\varepsilon}\|_{C^0} &\leq& \|\tilde w_m\|_{C^0}+\|y_{m}^{\varepsilon}\|_{C^0}+ \|s_1(\tilde w_m, y_{m}^{\varepsilon})\|_{C^0}\label{breveyC0} \\
  &\lesssim& \varepsilon_m^{\frac{13}{16}}+\varepsilon_m^{\frac34}+\varepsilon_m^{\frac32} \ \lesssim \ \varepsilon_m^{\frac34},\nonumber\\
\|\breve{y}_{m+1}^{\varepsilon}\|_{C^{\widehat R}} &\leq& \|\tilde w_m\|_{C^{\widehat R}}+\|y_{m}^{\varepsilon}\|_{C^{\widehat R}}+ \|s_1(\tilde w_m, y_{m}^{\varepsilon})\|_{C^{\widehat R}} \label{breveyCR}\\
   &\lesssim& \varepsilon_m^{\frac{7}{48}}+\varepsilon_m^{-\frac34}+1 \ \lesssim \ \varepsilon_m^{-\frac34}.\nonumber
\end{eqnarray}

Let $\CD_0$ and $\CD_0^\perp$ be the projection from $L^2(M,TM)$ onto $\IM d_0^*$ and $\Ker d_0$ respectively,  and recall $N_{m}=\varepsilon_{m}^{-\frac{1}{8(\tau+n+1)}}$.
Let us find $Y_m^{\varepsilon}$, $y_{m+1}^{\varepsilon}\in \Ga^{\widehat R}$ such that
\beq\label{conj_new-ym+1}
\Exp\{\breve{y}_{m+1}^{\varepsilon}\}\circ\Exp \{Y_m^{\varepsilon}\} = \Exp\{y_{m+1}^{\varepsilon}\},\quad
\begin{array}{l}
Y_m^{\varepsilon}\in \bigoplus_{\lam_j\leq N_{m}}\Ker (d_0\circ \PP_j),   \\[2mm]
y_{m+1}^{\varepsilon}\in \bigoplus_{\lam_j> N_{m}}\Ker (d_0\circ \PP_j)\bigoplus \IM d_0^*
\end{array}.\eeq
The above equation is equivalent to
\beq\label{eqyY_N}
\breve{y}_{m+1}^{\varepsilon}+Y_m^{\varepsilon}+s_1(\breve{y}_{m+1}^{\varepsilon},Y_m^{\varepsilon})=y_{m+1}^{\varepsilon}.\eeq Projecting onto $\bigoplus_{\lam_j\leq N_{m}} \Ker (d_0\circ \PP_j)$, we have
\beq\label{eqyY-ac}
( \Pi_{N_{m}}\circ \CD_0^{\perp}) \breve{y}_{m+1}^{\varepsilon}+Y_m^{\varepsilon}+( \Pi_{N_{m}}\circ \CD_0^{\perp}) s_1(\breve{y}_{m+1}^{\varepsilon},Y_m^{\varepsilon})= 0 .
\eeq
In view of Corollary \ref{cor-CR-truncation}, we have, for $0\leq R\leq \widehat R$,
\begin{eqnarray}
\left\|( \Pi_{N_{m}}\circ \CD_0^{\perp}) s_1(\breve{y}_{m+1}^{\varepsilon},Y_m^{\varepsilon})\right\|_{C^R}
&\lesssim& N_{m}^{\frac32n+1} \|s_1(\breve{y}_{m+1}^{\varepsilon},Y_m^{\varepsilon})\|_{C^R}\label{YN-CRnorm}\\
&\lesssim& N_{m}^{\frac32n+1} (\|\breve{y}_{m+1}^{\varepsilon}\|_{C^R}+\|\breve{y}_{m+1}^{\varepsilon}\|_{C^1}\|Y_m^{\varepsilon}\|_{C^R}),\nonumber
\end{eqnarray}
where we have used (\ref{s1-Cr}) of Proposition \ref{prop_s1-smooth} in the last inequality.
Interpolating with (\ref{breveyC0}) and (\ref{breveyCR}), we have
\beq\label{breveyC1}
\|\breve{y}_{m+1}^{\varepsilon}\|_{C^1}\lesssim \|\breve{y}_{m+1}^{\varepsilon}\|_{C^0}^{1-\frac{1}{\widehat R}} \|\breve{y}_{m+1}^{\varepsilon}\|_{C^{\widehat R}}^{\frac{1}{\widehat R}}\lesssim \varepsilon_m^{\frac34\frac{20(\tau+n)+18}{20(\tau+n+1)}} \leq \varepsilon_m^{\frac{11}{16}},\eeq
which implies, through (\ref{Nm3_16}), that
\beq\label{Nbry_small}
N_{m}^{\frac32n+1}\|\breve{y}_{m+1}^{\varepsilon}\|_{C^1} \lesssim \varepsilon_{m}^{-\frac{3}{16}} \varepsilon_m^{\frac{11}{16}}\leq \varepsilon_m^{\frac12} \ll 1.
\eeq
Hence, by the fixed point theorem, there is a unique $Y_{m}^{\varepsilon} \in \Gamma^{\widehat R}$ solving Eq. (\ref{eqyY-ac}) and which is $0$ when $\breve{y}_{m+1}^{\varepsilon}=0$.
Therefore, composing, on both sides of (\ref{alconj-breve}), by $\Exp\{Y_{m}^{\varepsilon}\}^{-1}$ from the left and by $\Exp\{Y_{m}^{\varepsilon}\}$ from the right, we obtain
\begin{eqnarray}
& &\Exp\{y_{m+1}^{\varepsilon}\}^{-1}  \circ \Exp\{P_{m+1}(\gamma)\}\circ \pi(\gamma)\circ \Exp\{y_{m+1}^{\varepsilon}\}\label{exp-tilde}\\
 &=&\Exp\{Y_m^{\varepsilon}\}^{-1} \circ\Exp\{z_m^{\varepsilon}(\gamma)\}\circ \pi(\gamma) \circ\Exp\{Y_m^{\varepsilon}\}\nonumber\\
  &=&\left(\Exp\{Y_m^{\varepsilon}\}^{-1} \circ\Exp\{z_m^{\varepsilon}(\gamma)\}\circ\Exp\{Y_m^{\varepsilon}\}\right)\circ \left(\Exp\{Y_m^{\varepsilon}\}^{-1}\circ \pi(\gamma)\circ \Exp\{Y_m^{\varepsilon}\} \right).\nonumber
\end{eqnarray}		
Let us write, for $\gamma\in \CS$,
\beq\label{tilde_z-pi}
\Exp\{Y_m^{\varepsilon}\}^{-1} \circ\Exp\{z_m^{\varepsilon}(\gamma)\}\circ\Exp\{Y_m^{\varepsilon}\}=:\Exp \{ z_{m+1}^{\varepsilon}(\gamma)\},
\eeq
and let us show that
\beq\label{z_0-pi}
\Exp\{Y_m^{\varepsilon}\}^{-1}\circ \pi(\gamma)\circ\Exp\{Y_m^{\varepsilon}\}=\pi(\gamma).
\eeq
According to Lemma \ref{lem_exp-isometry-smooth}, we have, for $\gamma\in \CS$,
$$  \pi(\gamma)\circ \Exp\{Y_m^{\varepsilon}\} \circ  \pi(\gamma)^{-1}=\Exp \{\pi(\gamma)_*Y_m^{\varepsilon}\}.$$
Noting that $Y_m^{\varepsilon}\in\Ker d_0$ implies
$\left(Y_m^{\varepsilon}-\pi(\gamma)_*Y_m^{\varepsilon}\right)_{\gamma\in \CS}=d_0Y_m^{\varepsilon}=0$,
we obtain (\ref{z_0-pi}) by
\begin{eqnarray*}
\Exp\{Y_m^{\varepsilon}\}^{-1}\circ \pi(\gamma)\circ\Exp\{Y_m^{\varepsilon}\}&=&\Exp\{Y_m^{\varepsilon}\}^{-1}\circ \left(\pi(\gamma)\circ\Exp\{Y_m^{\varepsilon}\}\circ \pi(\gamma)^{-1}\right)\circ\pi(\gamma)\\
&=&\Exp\{Y_m^{\varepsilon}\}^{-1}\circ \Exp \{Y_m^{\varepsilon}\} \circ\pi(\gamma) \  = \ \pi(\gamma).
\end{eqnarray*}	
Hence, combining (\ref{exp-tilde}) -- (\ref{z_0-pi}), we have
\begin{equation}\label{alconj-m+1}
\Exp\{y_{m+1}^{\varepsilon}\}^{-1}  \circ \Exp\{P_{m+1}(\gamma)\}\circ \pi(\gamma)\circ \Exp\{y_{m+1}^{\varepsilon}\}=\Exp \{z_{m+1}^{\varepsilon}(\gamma)\}\circ \pi(\gamma).
\end{equation}
Assuming $0<\varepsilon<\varepsilon_{m+1}$, let us show that
\beq\label{esti-yyz-alconj}
\|y_{m+1}^{\varepsilon}\|_{C^0}< \varepsilon_{m+1}^{\frac34},\quad \|y_{m+1}^{\varepsilon}\|_{C^{\widehat R}}< \varepsilon_{m+1}^{-\frac34}, \quad \|z_{m+1}^{\varepsilon}\|_{\CS,C^0}\lesssim  \|z_{m}^{\varepsilon}\|_{\CS,C^0},\eeq
which will complete the proof of Proposition \ref{prop-verify_alconj}.

In view of (\ref{eqyY-ac}), we see that
$$Y_m^{\varepsilon}=-  \left( \Pi_{N_{m}}\circ \CD_0^{\perp}\right)\left(\breve{y}_{m+1}^{\varepsilon} + s_1(\breve{y}_{m+1}^{\varepsilon},Y_m^{\varepsilon})\right).$$
According to Corollary \ref{cor-CR-truncation}, we have, for $0\leq R\leq \widehat R$,
$$
\|Y_m^{\varepsilon}\|_{C^{R}}
\lesssim N_{m}^{\frac32n+1} \left\|\breve{y}_{m+1}^{\varepsilon}+s_1(\breve{y}_{m+1}^{\varepsilon},Y_m^{\varepsilon})\right\|_{C^{R}}
\lesssim N_{m}^{\frac32n+1} \left(\|\breve{y}_{m+1}^{\varepsilon}\|_{C^{R}}+\|\breve{y}_{m+1}^{\varepsilon}\|_{C^1}\|Y_m^{\varepsilon}\|_{C^{R}}\right).
$$
Then, there exists a constant $c_{\widehat R}>0$ such that, for $0\leq R\leq \widehat R$,
$$\left(1-c_{\widehat R} N_{m}^{\frac32n+1}\|\breve{y}_{m+1}^{\varepsilon}\|_{C^1}\right)\|Y_m^{\varepsilon}\|_{C^{R}}\leq c_{\widehat R} N_{m}^{\frac32n+1} \|\breve{y}_{m+1}^{\varepsilon}\|_{C^{R}}, $$
which implies, through (\ref{Nbry_small}), $\|Y_m^{\varepsilon}\|_{C^{R}}\lesssim N_{m}^{\frac32n+1}  \|\breve{y}_{m+1}^{\varepsilon}\|_{C^{R}}$.
Therefore, according to (\ref{eqyY_N}) and (\ref{breveyC1}), we have,  for $0\leq R\leq \widehat R$,
\begin{eqnarray*}
\|y_{m+1}^{\varepsilon}\|_{C^{R}}&=&\|\breve{y}_{m+1}^{\varepsilon}+Y_m^{\varepsilon}+s_1(\breve{y}_{m+1}^{\varepsilon},Y_m^{\varepsilon})\|_{C^{R}} \\
&\lesssim & N_{m}^{\frac32n+1} \|\breve{y}_{m+1}^{\varepsilon}\|_{C^R}+  \|\breve{y}_{m+1}^{\varepsilon}\|_{C^1}\|Y_m^{\varepsilon}\|_{C^R}
\ \lesssim \ N_{m}^{\frac32n+1} \|\breve{y}_{m+1}^{\varepsilon}\|_{C^R} .
\end{eqnarray*}
In particular, since $N_{m}^{\frac32n+1}=\varepsilon_{m}^{\frac{\frac32n+1}{8(\tau+n+1)}}\leq\varepsilon_{m}^{-\frac{3}{16}+\frac{1}{16(\tau+n+1)}}$, we have, through (\ref{breveyCR}) and (\ref{breveyC1}),
\begin{eqnarray}
& &\|y_{m+1}^{\varepsilon}\|_{C^1} \lesssim \varepsilon_{m}^{-\frac{3}{16}}\varepsilon_m^{\frac{11}{16}}=\varepsilon_m^{\frac12},\label{ym+1C1}\\
& &\|y_{m+1}^{\varepsilon}\|_{C^{\widehat R}} \lesssim    \varepsilon_{m}^{-\frac{3}{16}+\frac{1}{16(\tau+n+1)}} \varepsilon_m^{-\frac34}  \ \lesssim \  \varepsilon_{m}^{\frac{1}{16(\tau+n+1)}}   \varepsilon_m^{-\frac{15}{16}} = \varepsilon_{m}^{\frac{1}{16(\tau+n+1)}}   \varepsilon_{m+1}^{-\frac34}.\label{ym+1CR}
\end{eqnarray}

According to (\ref{tilde_z-pi}), we have, for $\gamma\in \CS$,
$$ z_{m}^{\varepsilon}(\gamma)+ Y_{m}^{\varepsilon}+s_1(z_{m}^{\varepsilon}(\gamma),Y_{m}^{\varepsilon})=Y_{m}^{\varepsilon}+ z_{m+1}^{\varepsilon}(\gamma)+s_1(Y_m^{\varepsilon}, z_{m+1}^{\varepsilon}(\gamma)),$$
which implies $z_{m+1}^{\varepsilon}(\gamma)+s_1(Y_m^{\varepsilon}, z_{m+1}^{\varepsilon}(\gamma))= z_m^{\varepsilon}(\gamma)+s_1(z_m^{\varepsilon}(\gamma),Y_m^{\varepsilon})$.
Applying (\ref{s1-C0-decomp}) (with $w_1=Y_m^{\varepsilon}$ and $w_2=0$) and (\ref{s1-Cr+}) of Proposition \ref{prop_s1-smooth}, we have
$$\|s_1(Y_m^{\varepsilon}, z_{m+1}^{\varepsilon}(\gamma))\|_{C^0}\lesssim \|Y_m^{\varepsilon}\|_{C^1}\|z_{m+1}^{\varepsilon}(\gamma)\|_{C^0}
\lesssim N_{m}^{\frac32n+1}  \|\breve y_{m+1}^{\varepsilon}\|_{C^{1}}\|z_{m+1}^{\varepsilon}(\gamma)\|_{C^0} .$$
Applying (\ref{s1-C0-decomp}) of Proposition \ref{prop_s1-smooth} (with $w_1=0$ and $w_2=z_m^{\varepsilon}(\gamma)$), there is a constant $c>0$ such that
\begin{eqnarray*}
\left(1-cN_{m}^{\frac32n+1}\|\breve y_{m+1}^{\varepsilon}\|_{C^{1}} \right)\|z_{m+1}^{\varepsilon}(\gamma)\|_{C^0}&<&\| z_{m+1}^{\varepsilon}(\gamma)+s_1(Y_m^{\varepsilon}, z_{m+1}^{\varepsilon}(\gamma))\|_{C^0}\\
&=&\|z_m^{\varepsilon}(\gamma)+s_1(z_m^{\varepsilon}(\gamma),Y_{m}^{\varepsilon})\|_{C^0} \  \leq \ (1+c)\|z_m^{\varepsilon}(\gamma)\|_{C^0},
\end{eqnarray*}	
which implies, through (\ref{Nbry_small}), that $\|z_{m+1}^{\varepsilon}\|_{\CS,C^0}\lesssim\| z_m^{\varepsilon}\|_{\CS,C^0}$.

Recalling that $y_{m+1}^{\varepsilon}\in \bigoplus_{\lam_j> N_{m}}\Ker (d_0\circ \PP_j)\bigoplus \IM d_0^*$, we have
$$y_{m+1}^{\varepsilon}=\Pi_{N_{m}}^\perp y_{m+1}^{\varepsilon} + (\CD_0\circ \Pi_{N_{m}} )y_{m+1}^{\varepsilon},$$
where $\Pi_{N_{m}}^\perp y_{m+1}^{\varepsilon} $ satisfies that
\begin{eqnarray*}
\|\Pi_{N_{m}}^\perp y_{m+1}^{\varepsilon} \|_{C^0} \ \lesssim \  \|\Pi_{N_{m}}^\perp y_{m+1}^{\varepsilon} \|_{\CH^{\frac32n+1}}&=&\left( \sum_{j\in\N\atop{\lambda_j> N_{m}}} (1+\lambda_j)^{-2(\widehat R-\frac32n-1)}\|\PP_j y_{m+1}^{\varepsilon} \|^2_{\CH^{\widehat R}}\right)^\frac12\\
&\leq& N_{m}^{-(\widehat R-\frac32n-1)}\|y_{m+1}^{\varepsilon} \|_{\CH^{\widehat R}} \ \lesssim \ \varepsilon_m^{\frac{18(\tau+n+1)}{8(\tau+n+1)}} \varepsilon_m^{-\frac34} \  < \ \varepsilon_{m}^\frac32.
\end{eqnarray*}
Moreover, the conjugation (\ref{alconj-m+1}) implies that
$$\CP_{m+1}-d_0y_{m+1}^{\varepsilon}=\CZ_{m+1}^\varepsilon +s_1(y_{m+1}^{\varepsilon},\CZ_{m+1}^\varepsilon)-s_1(\CP_{m+1}, \pi_* y_{m+1}^{\varepsilon}).$$
According to Lemma \ref{lem_sol_cohomo-d0}, we have
\begin{eqnarray*}
\|(\CD_0\circ\Pi_{N_{m}})y_{m+1}^{\varepsilon}\|_{C^0}&\lesssim& N_{m}^{\tau+\frac32n+1}\|\CP_{m+1}- \CZ_{m+1}^\varepsilon -s_1(y_{m+1}^{\varepsilon},\CZ_{m+1}^\varepsilon)+s_1(\CP_{m+1}, \pi_* y_{m+1}^{\varepsilon}) \|_{C^0}\\
&\lesssim& \varepsilon_m^{-\frac{3}{16}}\left(\|\CP_{m+1}\|_{C^0}+(1+\|y_{m+1}^{\varepsilon}\|_{C^1})\|\CZ_{m+1}^\varepsilon\|_{C^0}\right)\\
&\lesssim& \varepsilon_m^{-\frac{3}{16}} \varepsilon_m^{\frac54} \  = \ \varepsilon_m^{\frac{17}{16}}.
\end{eqnarray*}
Hence, $\|y_{m+1}^{\varepsilon}\|_{C^0}\lesssim \varepsilon_m^{\frac{17}{16}}=\varepsilon_m^{\frac{1}{8}}\cdot\varepsilon_{m+1}^{\frac{3}{4}}$. Combining with (\ref{ym+1CR}), we obtain the estimates of $y_{m+1}^{\varepsilon}$ in (\ref{esti-yyz-alconj}) since, in view of (\ref{smallvarepsilon0}), all implicit coefficients in the above inequalities with ``$\lesssim$" are smaller than $\varepsilon_{0}^{-\frac{1}{16(\tau+n+1)}}\leq \varepsilon_{m}^{-\frac{1}{16(\tau+n+1)}}$.
Then (\ref{esti-yyz-alconj}) is proved.\qed


\subsection{Smoothness of $\CP_{m+1}$}\label{sec-smoothPm}

 With the above induction procedure, we obtain the sequence $\{\CP_{m}\}_{m\in\N}\subset {\bf G}_\pi$
satisfying
\beq\label{seqPm}
\|\CP_{m}\|_{C^{0}}\leq  \varepsilon_{m},\qquad \|\CP_{m}\|_{C^{R_*}}\leq  \varepsilon_{m}^{-1},\qquad \|\CP_{m}\|_{C^{R}}< \infty,\quad R>R_*,
\eeq
and $\{w_{m}\}_{m\in\N}\subset \bigoplus_{\lam_j\leq N_m}E_{\lam_j}\subset \Gamma^\infty(M,TM)$ satisfying (\ref{esti_wm-smooth}).
Following the methods of Zehnder \cite{Zehnder75} and of Moser \cite{Moser-cercle}, we show that
$\{\CP_{m}\}_{m\in\N}$ and $\{\Exp\{w_{m}\}\}_{m\in\N}$ converge in $\Gamma^\infty(M,TM)$.

For $R'\in\N^*$, let $\tilde c_{R'}>1$ be the maximal implicit constant depending on $R'\in \N$ with $R'\leq R$ in all propositions, lemmas and corollaries in Section \ref{section-geom} and Lemma \ref{lemma_reste-smooth}, \ref{lem_sol_cohomo-d0}, and let $\hat c_R:=\max_{R'\leq R}\{\tilde c_{R'}\}$.
For fixed $R\in\N^*$, let $m$ be large enough such that
 \begin{equation}\label{suffi-large-m}
 \hat c_R < \varepsilon_{m}^{-\frac1{400}}.\end{equation}
\begin{prop}\label{prop-large_r}
Given any $R\in\N^*$,  if $m$ is sufficiently large, then
$\|\CP_{m}\|_{C^R}< \varepsilon_{m}^{\frac{99}{100}}$.
\end{prop}

\proof Since for any $R\in\N^*$, $\|\CP_{m}\|_{C^R}< \infty$, we always have $\|\CP_{m}\|_{C^R}< D_{R,m}  \varepsilon_{m}^{-1}$ for some constant $D_{R,m}>1$. At first, let us show that, for fixed $R\in \N^*$, these coefficients are uniformly bounded (w.r.t. $m$), i.e., there is $m_R\in\N$ and $D_R>1$, such that
\begin{equation}\label{induction-large-m}
\|\CP_{m}\|_{C^R}< D_R\varepsilon_{m}^{-1}, \quad m\geq m_R.\end{equation}
In view of (\ref{seqPm}), it is true for $R\leq R_*$ with $D_R=1$ and $m_R=0$. Let us consider the situation $R>R_*$.
Suppose that, for some large enough $m\in\N^*$ such that (\ref{suffi-large-m}) is satisfied,
and some $D_R>1$, we have $\|\CP_{m}\|_{C^R}<  D_R \varepsilon_{m}^{-1}$.
It is sufficient to show that $\|\CP_{m+1}\|_{C^R}< D_R \varepsilon_{m}^{-\frac54}= D_R \varepsilon_{m+1}^{-1}$.

Recalling (\ref{Pm+1}), we consider the inequality
\begin{eqnarray}
\|\CP_{m+1}\|_{C^R}
&\leq&\|\Pi_{N_m}^\perp \CP_{m} + (\Pi_{N_m}\circ\D_0^\perp)\CP_m\|_{C^R}\label{pm+1Cr}\\
& &  + \, \|s_1(w_m,\CP_{m+1})\|_{C^R} +\|s_1(\CP_m, \pi_* w_m)\|_{C^R}.\nonumber
\end{eqnarray}
According to Corollary \ref{cor-CR-truncation} and recalling (\ref{Nm3_16}), we have that
$$\|\Pi_{N_m}^\perp \CP_{m}+ (\Pi_{N_m}\circ\D_0^\perp)\CP_m\|_{C^R} <2 \hat c_R \varepsilon_{m}^{-\frac{3}{16}}\|\CP_{m}\|_{C^R}
< \hat c_R D_R \varepsilon_m^{-\frac{19}{16}}
< D_R \varepsilon_{m}^{-\frac{39}{32}}, $$
noting that (\ref{suffi-large-m}) implies $2\hat c_R\varepsilon_m^{\frac{1}{32}}<1 $.

Recall that $\|\CP_m\|_{C^{1}}, \|w_m\|_{C^{1}}\lesssim \varepsilon_m^{\frac34}$,
 which implies that $\|\CP_m\|_{C^{1}}, \|w_m\|_{C^{1}}< \varepsilon_m^{\frac23}$.
Then, by \re{s1-Cr} in \rp{prop_s1-smooth} and (\ref{esti_wCR}), we have that
\begin{eqnarray*}
\|s_1(\CP_m, \pi_* w_m)\|_{C^R} &<&  \hat c_R \left(\|\CP_m\|_{C^R} + \|\CP_m\|_{C^{1}}\|w_m\|_{C^R}\right)\\
 &<& \left(\hat c_R + \hat c_R^2\varepsilon_m^{\frac23}\right)  \|\CP_m\|_{C^R} \ \leq \ D_R \varepsilon_{m}^{-\frac{101}{100}}, \\
\|s_1(w_m,\CP_{m+1})\|_{C^R} &\leq& \hat c_R \left(\|w_m\|_{C^R} + \|w_m\|_{C^{1}}\|\CP_{m+1}\|_{C^R}\right)\\
 &<& \hat c_R^2 \varepsilon_{m}^{-\frac{3}{16}} \|\CP_m\|_{C^R}+\hat c_R\varepsilon_m^{\frac23}\|\CP_{m+1}\|_{C^R}
  \ < \ D_R \varepsilon_{m}^{-\frac{39}{32}}+\varepsilon_m^{\frac12}\|\CP_{m+1}\|_{C^R}
 \end{eqnarray*}
Hence, collecting all the previous inequalities into \re{pm+1Cr}, we obtain
$$ \|\CP_{m+1}\|_{C^R} \leq D_R\varepsilon_m^{-\frac{101}{100}}+ 3D_R\varepsilon_m^{-\frac{39}{32}}+ \varepsilon_m^{\frac12}\|\CP_{m+1}\|_{C^R}, $$
which implies that
$$
\|\CP_{m+1}\|_{C^R}\leq  4\left(1-\varepsilon_{m}^{\frac12}\right)^{-1} D_R \varepsilon_m^{-\frac{39}{32}}
\leq D_R \varepsilon_{m}^{-\frac54}=D_R \varepsilon_{m+1}^{-1}, $$
since, $4\left(1-\varepsilon_{m}^{\frac12}\right)^{-1} \varepsilon_m^{\frac{1}{32}}<1$. Hence, (\ref{induction-large-m}) is shown.

Now let $R_*<R<R+2n+1<\tilde R$.
In view of (\ref{induction-large-m}) and \rp{propNorm-CrHr}, there exist $m_{\tilde R}\in\N$ and $D_{\tilde R}>1$ such that, for $m\geq m_{\tilde R}$, $\|\CP_{m}\|_{\CH^{\tilde R}}\leq \hat c_{\tilde R} D_{\tilde R} \varepsilon_{m}^{-1}$.
Hence, by interpolation inequality in Lemma \ref{interpol-hs}, we have
$$\|\CP_{m}\|_{C^R}\leq \hat c_R \|\CP_{m}\|_{\CH^{R+2n+1}} \leq \hat c_R  \|\CP_{m}\|_{\CH^0}^{1-\frac{R+2n+1}{\tilde R}} \|\CP_{m}\|_{\CH^{\tilde R}}^{\frac{R+2n+1}{\tilde R}} \leq  \varepsilon_{m}^{\frac{99}{100}},$$
as soon as
$\hat c_R (\hat c_R)^{1-\frac{R+2n+1}{\tilde R}}\left(\hat c_{\tilde R} D_{\tilde R}\right)^{\frac{R+2n+1}{\tilde R}}\varepsilon_{m}^{\frac{1}{100}} <\varepsilon_{m}^{\frac{2(R+2n+1)}{\tilde R}}$,
which certainly holds if $200(R+2n+1)<\tilde R$, say $\tilde R
=300(R+2n+1)$ and $m$ large enough, depending on only $R$, so that $\hat c_R^{2\tilde R-(R+2n+1)}( \hat c_{\tilde R} D_{\tilde R})^{R+2n+1}\varepsilon_{m}^{R}< 1$. The proposition is shown.\qed

\subsection{Convergence in $\Gamma^\infty(M,TM)$}\label{sec_cov-smooth}
With the sequence $\{w_m\}_{m\in\N}\subset \bigoplus_{\lam_j\leq N_m} E_{\lam_j}$ built above, satisfying $\|w_m\|_{C^1}\lesssim \varepsilon_m^\frac34$,
let $W_1= w_0+w_1+s_1(w_0,w_1)$. Then, by (\ref{s1-Cr}), we have,
\beq\label{W1C1}
\|W_1\|_{C^1}\lesssim \|w_0\|_{C^1}+\|w_1\|_{C^1} + \|s_1(w_0,w_1)\|_{C^1}\lesssim  \|w_0\|_{C^1}+\|w_1\|_{C^1} \lesssim \varepsilon_0^\frac34+ \varepsilon_1^\frac34\leq \varepsilon_1^\frac12. \eeq
According to Proposition \ref{prop_s1-smooth}, we have $\Exp\{w_0\}\circ\Exp\{w_1\}=\Exp\{W_1\}$, which implies that, for every $\gamma\in \CS$,
\begin{eqnarray}
& & \Exp\{W_1\}^{-1}\circ \Exp\{P_0(\gamma)\} \circ \pi(\gamma) \circ\Exp\{W_1\}\label{W1}\\
&=&\Exp\{w_1\}^{-1}\circ\left(\Exp\{w_0\}^{-1}\circ \Exp\{P_0(\gamma)\} \circ \pi(\gamma)\circ\Exp\{w_0\} \right)\circ\Exp\{w_1\} \nonumber\\
&=&\Exp\{w_1\}^{-1}\circ \Exp\{P_1(\gamma)\} \circ \pi(\gamma) \circ \Exp\{w_1\} \nonumber\\
&=& \Exp\{P_2(\gamma)\} \circ \pi(\gamma).\nonumber
\end{eqnarray}

For $R\in\N$, according to Proposition \ref{prop-large_r}, there exists $m_R\in\N^*$ such that
$$\|\CP_{m}\|_{C^{R}} \leq \varepsilon_{m}^{\frac{99}{100}},\quad m\geq m_R.$$
Suppose that there exists $W_m\in\Ga^\infty (M,TM)$ such that, for every $\gamma\in \CS$,
$$\Exp\{W_{m}\}^{-1}\circ \Exp\{P_0(\gamma)\} \circ \pi(\gamma) \circ\Exp\{W_{m}\}=\Exp\{P_{m+1}(\gamma)\} \circ \pi(\gamma).$$
With $W_{m+1}:=W_m+w_{m+1}+s_1(W_m, w_{m+1})$, we have, through Lemma \ref{lem_moser-s1}, that
$\Exp\{W_m\}\circ\Exp\{w_{m+1}\}=\Exp\{W_{m+1}\}$,
which implies that, for every $\gamma\in \CS$,
\begin{eqnarray}
& &\Exp\{W_{m+1}\}^{-1}\circ \Exp\{P_0(\gamma)\} \circ \pi(\gamma) \circ\Exp\{W_{m+1}\}\label{Wm+1}\\
&=& \Exp\{w_{m+1}\}^{-1}\circ\left(\Exp\{W_{m}\}^{-1}\circ \Exp\{P_0(\gamma)\} \circ \pi(\gamma)\circ\Exp\{W_{m}\}\right)\circ\Exp\{w_{m+1}\}\nonumber\\
&=& \Exp\{w_{m+1}\}^{-1}\circ \Exp\{P_{m+1}(\gamma)\} \circ \pi(\gamma) \circ\Exp\{w_{m+1}\}\nonumber\\
&=&\Exp\{P_{m+2}(\gamma)\} \circ \pi(\gamma).\nonumber
\end{eqnarray}
By (\ref{esti_wCR}), we have $\|w_{m+1}\|_{C^{R}}\leq \varepsilon_{m+1}^{-\frac14}  \|\CP_{m+1}\|_{C^{R}}\leq \varepsilon_{m+1}^{\frac{99}{100}-\frac14}$.
By (\ref{s1-Cr+}) in Proposition \ref{prop_s1-smooth},
$$
\|W_m+s_1(W_m, w_{m+1})\|_{C^R}\leq \|W_m\|_{C^{R}}+\hat c_R (\|W_m\|_{C^{R}}+\|W_m\|_{C^{1}}\|w_{m+1}\|_{C^{R}})
\leq3\hat c_R\|W_m\|_{C^{R}}.
$$
Hence, noting that $\varepsilon_{m+1}^{\frac{99}{100}-\frac14}\leq \varepsilon_{m+1}^{\frac23}$, we obtain an upper bound for $\|W_{m+1}\|_{C^{R}}$:
\begin{eqnarray*}
\|W_{m+1}\|_{C^{R}}&\leq&3\hat c_R\|W_m\|_{C^{R}}+ \varepsilon_{m+1}^{\frac23} \ \leq \ \cdots\\
&\leq&  (3\hat c_R)^{m}\left(\|W_{1}\|_{C^{R}}+\sum_{j=0}^{m-1}(3\hat c_R)^{j-m}\varepsilon_{m-j+1}^{\frac23}\right)\\
&\leq&  (3\hat c_R)^{m}\left(\|W_{1}\|_{C^{R}}+\sum_{j=0}^{m-1}\varepsilon_{m-j+1}^{\frac23}\right) \ =  \  E_R^m F_R,\end{eqnarray*}
with $E_R:=3\hat c_R$ and $F_R:=\|W_{1}\|_{C^{R}}+\sum_{j\geq 2}\varepsilon_{j}^{\frac23}$.
According to (\ref{s1-Cr+}), if $m\geq m_{R+1}$, then
\begin{eqnarray*}
\|W_{m+1}-W_m\|_{C^R}&\leq&\|w_{m+1}\|_{C^{R}}+\|s_1(W_m, w_{m+1})\|_{C^{R}}\\
&\leq& \hat c_R \left(\|w_{m+1}\|_{C^{R}} +  \|W_m\|_{C^{2}}\|w_{m+1}\|_{C^{R}}+\|W_m\|_{C^{R+1}}\|w_{m+1}\|_{C^{0}}\right)\\
&\leq& \hat c_R \left(1+E_{2}^mF_{2}+E_{R+1}^mF_{R+1}\right)\varepsilon_{m+1}^{\frac23}\\
	&=& \hat c_R \left(1+E_{2}^m F_{2}+E_{R+1}^mF_{R+1}\right) \varepsilon_{0}^{\frac23(\frac54)^{m+1}}  \ \leq \ \varepsilon_{0}^{\frac23(\frac98)^{m+1}} .
\end{eqnarray*}
Hence, $\sum_{m\in \N^*} \|W_{m+1}-W_m\|_{C^R} $ converges which implies the $C^R-$convergence of the sequence $\{W_m\}_{m\in\N^*}$.
In particular, $\|W_{m+1}-W_m\|_{C^1}\lesssim \varepsilon_{0}^{\frac23(\frac98)^{m+1}}$, combining with \re{W1C1}, $\|W_m\|_{C^1}$ is uniformly bounded by $ \varepsilon_{0}^{\frac14}$. Therefore, $s_1(W_m, w_{m+1})$ is well-defined.

For the limit $W:=\lim_{m\to\infty}W_m$ in $\Ga^\infty(M,TM)$, $\Exp\{W\}$ defines a smooth diffeomorphism of $M$ such that
\beq\label{rig-smooth}
\Exp\{W\}^{-1}\circ\pi_0(\gamma) \circ\Exp\{W\}= \pi(\gamma),\quad \gamma\in \CS.\eeq
Theorem \ref{thmmain} and \ref{thm-geom0} are proved in the smooth case.

\section{Grauert tube and Hardy space}\label{sec_Grauert_Hardy}

Now we assume further that the smooth Riemannian manifold $M$ is real analytic.
This section is dedicated to the analytic case and is irrelevant to the smooth case.

\subsection{Grauert tube for a real analytic Riemannian manifold}\label{sec_Grauert}

We recall (without proofs) some useful facts for real analytic Riemannian manifolds stated in \cite{GLS96}[Section 1].
First of all, according to Bruhat-Whitney theorem \cite{WB59} (see also \cite{GLS96}[Lemma 1.2]), a compact real analytic Riemannian manifold $M$ can be identified with a totally real submanifold of a complex analytic manifold $\tilde M$ of (real) dimension $2n$~: for all $m\in M$, there exists an open neighborhood $W$ of $m$ in $\tilde M$ and a holomorphic coordinate system $(z_1,\cdots,z_n)$ on $W$ such that
\beq\label{im}
W\cap M=\{q\in W : \IM z_1(q)=\cdots=\IM z_n(q)=0 \}.
\eeq
We also recall a well-known fact (see \cite{GLS96}[Corollary 1.3]).
\begin{prop}\label{holom-ext}\cite{Tom66} Let $M\hookrightarrow \tilde M$ be a totally real submanifold of a complex manifold $\tilde M$. Let $M'$ be a complex manifold and let  $f:M\rightarrow M'$ be a real analytic mapping. Then, there exists an open connected neighborhood $W$ of $M$ in $\tilde M$ and a unique holomorphic mapping $f^+:W\rightarrow M'$ such that $\left. f^+\right|_{M}=f$.
\end{prop}



Following \cite{GLS96}[Corollary 1.3], for $m\in M$, there exists an open connected neighborhood $W_m\subset T_mM\otimes\Bbb C$ and a unique holomorphic extension of $\Exp_m$ on $W_m$, still denoted by $\Exp_m$, to $\tilde M$.
Moreover, according to \cite{GLS96}[Theorem 1.5], there exists $0<r_*\leq\inf_{m\in M}r(m)$ such that for every $0<r<r_*$, the map
\beq\label{complex-exp}
\Phi : T^{r}M \rightarrow \tilde M,\quad \Phi(m,\xi)=\Exp_m\{{\rm i}\xi\}
\eeq
is an analytic diffeomorphism onto its image, where $$T^{r}M:=\left\{(m,\xi)\in TM : |\xi|_{g(m)}<r\right\}.$$
According to \cite{Szo91}[Theorem 2.2], \cite{GLS96}[Proposition 1.7] and \cite{LeSz91}, for any $0<r<r_*$, $T^{r}M$ admits a unique complex structure for which the complexified exponential
$$
T^{r}M \ni(m,\xi)  \mapsto \Exp_m\{{\rm i}\xi\} \in \Phi(T^{r}M)=: M_{r}
$$
is a biholomorphism. We shall write $TM^{\mathbb{C}}:=TM\otimes_{\mathbb{R}}\mathbb{C}$.

According to \cite{Gra58}(see also \cite{GS91,GS92}[Introduction]), there exists a non-negative smooth strictly plurisubharmonic function
\beq\label{rho}
\rho:  M_{r_*}\rightarrow [0,r_*[  \  \text{ with } \rho^{-1}(0)=M \text{ and } M_{r}=\rho^{-1}([0,r[), \ 0<r<r_*.
\eeq
Moreover, there exists an anti-holomorphic involution $\sigma :  M_{r_*}\rightarrow  M_{r_*}$ whose fixed point set is $M$ and $\rho(\sigma(q))=\rho(q)$ for all $q \in  M_{r_*}$.

Since the metric $g$ on $M$ is real analytic, it turns out that such a $M_{r}$ can be defined by a unique {\it real analytic} strictly plurisubharmonic function $\rho$ such that the K\"ahler form
$$
\omega:=\frac{\rm i}{2}\partial\bar\partial \rho= \frac{\rm i}{2}\sum_{1\leq i,j\leq n}\frac{\partial^2\rho}{\partial z_i\partial \bar z_j}dz_i\wedge d\bar z_j
$$ defines a K\"ahler metric on $M_{r}$, $0<r<r_*$,  \beq\label{kappa}\kappa:= \sum_{1\leq i,j\leq n}\frac{\partial^2\rho}{\partial z_i\partial \bar z_j}dz_i\otimes d\bar z_j,\eeq
which extends the Riemannian metric $g$ on $M$ according to the following theorem.

\begin{thm}\cite{GS91}[P.562]\label{GS-thm}
There exists a neighborhood $U$ of $M$ in $\tilde M$ and a unique real analytic solution $\rho$ on $U\setminus M$ of the complex Monge-Ampère equation
\beq\label{MA}
\det\left(\frac{\partial^2 \sqrt{\rho}}{\partial z_i\partial\bar z_j}\right)=0
\eeq
such that the inclusion map $(M,g)\hookrightarrow (\tilde M, \kappa)$ is an isometric embedding.
\end{thm}
\noindent Hence, the boundary $\partial M_{r}$ of $M_{r}$ is a compact real analytic manifold and $\overline{M}_r:=\rho^{-1}([0,r])$ is a compact K\"ahler manifold. The complex neighborhood $M_r$ of $M$ is called a {\bf Grauert tube} of width $r$. 


\smallskip

Let us first extend the exponential map $ \Exp_m $ introduced in Section \ref{secExp} w.r.t. the metric $g$ at $m\in M$,
  to the one w.r.t. the metric $\kappa$ at $q\in M_{r}$,
$$\Exp_{q} : B_{q}(0,r(q))\subset T_{q}^{(1,0)} M_{r}\rightarrow M_{r}.$$
Let $X$ be a real analytic vector field on $M$. According to \rp{holom-ext}, it extends to a holomorphic vector field $X$ on an open connected neighborhood $\CU$ of $M$ in $ M_{r}$, still denoted $ M_{r}$.
That is, $X$ is a holomorphic section of $T^{(1,0)} M_{r}$ over $ M_{r}$.
First of all, according to \cite{GLS96}[Proposition 1.9, 1.13], if $r$ is small enough, the analytic Riemannian metric $g$ uniquely extends to a non-degenerate holomorphic section $g^+\in \Ga^\om\left( M_{r}, BS(T^{(1,0)} M_{r})\right)$, where $BS(T^{(1,0)} M_{r})$ denotes the bundle of symmectric bilinear forms on the holomorphic vector fields $T^{(1,0)} M_{r}$ over $M_{r}$, that is $g^+$ defines a {\it holomorphic Riemannian metric} \cite{Lebrun}.
For each $q\in M_{r}$, we define $\Exp_{q}$ on the ball $B_q(0,r(q))\subset T^{(1,0)}_{q}M_{r}$ with respect to its K\"ahler metric $\kappa$ as follow~:
given a coordinate chart $(U,x)=(U,x_1,\ldots, x_n)$ of $M$ trivializing $TM$, let $(W,z_1,\ldots, z_n)$ be a holomorphic chart of $M_{r}$ extending a chart $(U,x)$ of $M$ as in $(\ref{im})$, with $W\cap M=U$ and $x_i=\RE z_i$ trivializing $T^{(1,0)} M_{r}$. Let us write
$$g(x(m))=\sum_{1\leq i,j\leq n} g_{i,j}(x(m))dx_i\otimes dx_j,\quad g^+(z(q))=\sum_{1\leq i,j\leq n} g_{i,j}^+(z(q))dz_i\otimes dz_j,
$$
where the matrices $(g_{i,j}(x(m)))_{1\leq i,j\leq n}$, $(g_{i,j}^+(z(q)))_{1\leq i,j\leq n}$ are invertible for each point  $m\in M$ and $q\in M_{r}$ respectively. We recall that the geodesics on $M$ are solutions of the (real time) differential equation, in a coordinate chart~:
$$
\ddot x_j = \sum_{1\leq k,l\leq n}\Ga^j_{k,l}(x)\dot x_k\dot x_l,\quad j=1,\cdots, n,
$$
where, $\Ga^j_{k,l}$ denotes the Christoffel symbol defined by
$$
\Ga^j_{k,l}(x):=\frac{1}{2}\sum_{1\leq m\leq n} g^{j,m}(x) \left(\frac{\partial g_{m,k}}{\partial x_l}-\frac{\partial g_{k,l}}{\partial x_m}+\frac{\partial g_{l,m}}{\partial x_k}\right),
$$
and $(g^{j,m}(x))$ denotes the inverse matrix of $(g_{j,m}(x))$. 
Following \cite{Lebrun}[1.17, 1.18], let us consider the {\it holomorphic differential equation} with complex time:
\beq\label{geod-complex}
\ddot z_j = \sum_{1\leq k,l\leq n}\Ga^{+j}_{k,l}(z)\dot z_k\dot z_l,\quad j=1,\cdots, n,
\eeq
with $\Ga^{+j}_{k,l}$ defined as
$$
\Ga^{+j}_{k,l}(z):=\frac{1}{2}\sum_{1\leq m\leq n} (g^+)^{j,m}(z) \left(\frac{\partial g_{m,k}^+}{\partial z_l}-\frac{\partial g_{k,l}^+}{\partial z_m}+\frac{\partial g_{l,m}^+}{\partial z_k}\right),
$$
and $\left((g^+)^{j,m}(z)\right)$ the inverse matrix of $\left(g_{j,m}^+(z)\right)$. For any $q_0\in W$ and $(q_0,\xi)\in T_{q_0}^{(1,0)}M_{r}$, such that $(z_0,\xi)\in \Delta_1^n\times \mathbb{C}^n$ with $z_0=z(q_0)$, there exists a unique complex curve, a {\it complex geodesic}, $t\in D_{z_0,\xi}\mapsto z(t)=\Phi(t,z_0,\xi)$ with $(z(0),\dot z(0))=(z_0,\xi)$,
  solution of $(\ref{geod-complex})$. Here $D_{z_0,\xi}$ denotes a complex neighborhood of $0$ in $\Bbb C$ that depends on the point $(z_0,\xi)$. As in the real case, the form of Eq. \re{geod-complex} allows us to write $z(t)=\Psi(z_0,t\xi)$; it is holomorphic for $z_0\in\Delta_1^n$ and $t$ small complex number. Hence, $\Psi(z_0,\xi)$ is holomorphic for  $z_0\in\Delta_1^n$, and  $\xi$ in the complex ball 
in $\mathbb{C}^n$, centered at $0$ and of sufficiently small radius $\del$ w.r.t the K\"ahler metric $\kappa$~: $|\xi|_{\ka(q_0)}<\del$, $z(q_0)=z_0$.
Furthermore, it satisfies
\beq\label{basic-flow}
\Psi(z_0,0)=z_0,\quad D_{\xi}\Psi(z_0,0)=\id.
\eeq
Hence, for some holomorphic map $\varphi(z_0,\xi)$ satisfying $D_{\xi}\varphi(z_0,0)=0$, we have
\beq\label{exp-def}
z(t)= \Psi(z_0,t\xi)=z_0+t\xi+\varphi(z_0,t \xi).
\eeq
Taking a finite covering of $M_r$ by open sets, there exists an $a>0$ such that the solution $z(t)=\Phi(t,z(q),\xi)=\Psi(z(q),t\xi)$ is holomorphic for $|t|< 2$, $q\in  M_{r}$ and $|\xi|_{\kappa(q)}<a$. 
Let $(q,\xi)\in  T^{(1,0)}_{q} M_{r}$, be such a point (i.e. $|\xi|_{\ka}< a$) with $q\in W$. We define the {\it complex exponential map} $\Exp_{q}\{\xi\}$ to be the time$-1$ of this complex flow. It is the point of $M_{r}$ whose expression in the coordinate chart $W$ is
\beq\label{complex-exp}
\Exp_{q}\{\xi\}:= \Psi(z(q),\xi).
\eeq

Let $q\in W$ of sufficiently small coordinate $z(q)$. Let $(\xi,\eta)\in \mathbb{C}^{n}\times \mathbb{C}^n$ be small enough so that, 
$\Psi(\Psi(z(q),\xi),\eta)$ is well defined. According to \re{basic-flow}, there is a unique holomorphic map $(z,\xi,\eta)\mapsto P(z,\xi,\eta)=:\zeta\in \mathbb{C}^n$ 
that solves the equation  $\Psi(\Psi(z,\xi),\eta)=:\Psi(z,\zeta)$ for $(\xi,\eta)$ in small neighborhood of $0$ in $\mathbb{C}^{2n}$ and $z$ in a neighborhood of $z(q)$. Furthermore, there is a holomorphic map $\varrho(z,\xi,\eta)\in \mathbb{C}^n$ such that
$$ \zeta =\xi+\eta+\varrho(z,\xi,\eta),\quad \varrho(z,0,\eta)=\varrho(z,\xi,0)=0,\quad |D_\xi\varrho(z,\xi,\eta)|_{\ka}\lesssim |\eta|_{\ka},$$
where, and afterwards in the analytic setting, the inequality with ``$\lesssim$" means boundedness from above by a positive constant depending only on the manifold $(M_{r_*},\ka)$ but independent of other factors.


Given $0<r<r_*$, let the set of holomorphic sections of $T^{(1,0)} M_{r}$ over $M_{r}$, that is holomorphic vector fields, be denoted by $\Gamma_{r} = \Gamma( M_{r},T^{(1,0)} M_{r})$,
equipped with the norm
$$
 |v|_{0,r}:=\sup_{q\in M_{r}}|v(q)|_{\ka},\quad v\in \Gamma_{r}.
$$
There is an analytic trivializing atlas with a finite covering patches $\left\{U_i,x^{(i)}\right\}_i$ of $M$ that extends to a holomorphic atlas $\left\{W_i,z^{(i)}\right\}_i$ of $M_{r_*}$ as in \re{im} and such that $z^{(i)}\in \Delta_1^n$ (In what follows, $z^{(i)}$ stands for $z^{(i)}(q)$ with $q\in W_i$). For $v\in\Gamma_{r} $, restricting to $W$, one of these coordinates patches on which $z\in \Delta_1^n$ and writing $\tilde v(z)=\sum_{1\leq j\leq n} \tilde v_j(z)\frac{\partial}{\partial z_j}$ the expression of $v$ in this coordinate patch, we set
\begin{eqnarray}
 \|v\|_{C^0,r}&:=&\max_i\sup_{q\in W_i\cap M_{r}}|\tilde v(z^{(i)}(q))|_{\kappa(q)},\label{norm0r} \\
  \|v\|_{C^1,r}&:=& \|v\|_{C^0,r}+\max_i\sup_{q\in W_i\cap M_{r}}\sup_{\substack{\zeta\in \mathbb{C}^n,\\|\zeta|\leq 1}}|D\tilde v(z^{(i)}(q))\zeta|_{\kappa(q)}.\label{norm1r}
\end{eqnarray}
In particular, if $v\in \Gamma^{\omega}\subset  \Gamma^{\infty}$, the above norms with $r=0$ are equivalent to the $\|\cdot\|_{C^0}$ and $\|\cdot\|_{C^1}$ norms defined in (\ref{normsCk}).
Since every $v\in \Gamma^{\omega}$ can be holomorphically extended to $M_r$ for some $0<r<r_*$, let $\Gamma_r^{\omega}\subset \Gamma_r$ be the set of holomorphic extensions to $M_r$ of elements in $\Gamma^{\omega}$.
It is obvious that $ \|v\|_{C^0}\leq \|v\|_{C^0,r}$ for $v\in \Gamma^\om_r$.

\begin{prop}\label{propMoser-anal}
Given $0<r<r_*$, for $v,w\in\Gamma_{r}$ with $\|w\|_{C^1,r}$ and $\|v\|_{C^0,r}$ sufficiently small (depending only on the manifold $(M_{r_*},\ka)$), there exists $s_1(w,v)\in \Gamma_{r'}$ for any $r'\in ]0,r[$ such that,
\begin{equation}\label{eq_lemMoser1}
\Exp \{w\}\circ \Exp \{v\} = \Exp \{w+v+s_1(w,v)\},
\end{equation}
with $s_1(w,0)=s_1(0,v)=0$, and for any $r'\in ]0,r[$,
$\|s_1(w,v)\|_{C^0,{r'}}\lesssim \|w\|_{C^1,r}\|v\|_{C^0,{r'}}$.
\end{prop}

\begin{remark}\label{rem_s1}
In the above proposition, if $\tilde v, \tilde w\in \Gamma_r^\om$ are holomorphic extensions to $M_r$ of $v,w\in\Gamma^\om$, then $s_1(\tilde w,\tilde v)\in\Gamma_r^\om$ is the holomorphic extension of $s_1(w,v)\in\Gamma^\om$.
\end{remark}

Proposition \ref{propMoser-anal} is the holomorphic version of Lemma \ref{lem_moser-s1} on $M_r$. It can be
deduced readily from the proof in \cite{Mos69}.
For completeness, we give a proof of this proposition in Appendix \ref{app_proof}.

\subsection{Hardy space and weighted $\bf L^2-$norm}\label{sec-Hardy}

Since for the real analytic manifold $M$, the Riemannian metric $g$ extends to a K\"ahler metric $\kappa$ on the Grauert tube $M_{r}$, $0<r<r_*$, for any holomorphic sections $v,w\in \Gamma^\om_{r}$, $$
\la v,w\ra:=\int_{M_{r}}\la v(z),w(z)\ra_{\ka} \, \frac{\om^n(z)}{n!}.
$$

Now we follow and recall the result of Boutet de Monvel \cite{BdM78}(see also \cite{GLS96}[Section 1 and 2] and \cite{Leb18}). Let us consider the elliptic analytic pseudo-differential operator of order $1$, $|\De_{TM}|^{\frac12}$. Due to the classical elliptic theory, the eigenvectors of $|\De_{TM}|^{\frac12}$ (which are the same as those of $\De_{TM}$) are, in fact, real analytic. They can be considered as restrictions to $M$ of holomorphic sections on a same neighborhood of the $M$ in a complexified manifold $\tilde M$ of $M$.
\begin{defi}\label{def-Hardy}
Let the {\bf Hardy space} $\tilde H_r^2=\tilde H^2(M_{r},T^{(1,0)}M_{r})$ be the space of holomorphic sections of $T^{(1,0)}M_{r}$ over $M_{r}$ whose restriction (in the sense of distribution) to $\partial M_{r}$ belongs to $L^2(\partial M_{r},T^{(1,0)}M_{r})$, associated to
the {\bf Hardy product}
\beq\label{h2product}
\la f,h \ra_{\tilde H_r^2}:=
\int_{\partial M_{r}}\la f(q),h(q)\ra_{\ka}d\mu_{r}(q),\qquad f, h \in  \tilde H_r^2,
\eeq
and the {\bf Hardy norm}
\beq\label{h2norm}
	\|f\|_{\tilde H_r^2}:=\|f|_{\partial M_{r}}\|_{L^2(\partial M_{r})}=\left(\int_{\partial M_{r}}\la f(q),f(q)\ra_{\ka}d\mu_{r}(q)\right)^{\frac12},\qquad f\in\tilde H_r^2.
\eeq
Here, $d\mu_{r}$ denotes the ``surface measure" obtained by restriction of $\frac{\om^n(z)}{n!}$ to the real analytic level set $\rho=r$.
More generally, for $\nu\in\N^*$, the Hardy product and the Hardy norm on $(H_r^2)^\nu$ are defined as
\begin{eqnarray}
\|f\|_{\tilde H_r^2}^2:= \sum_{1\leq l\leq \nu}\|f_l\|_{\tilde H_r^2}^2, \qquad f=(f_l)_{1\leq l\leq \nu}\in (\tilde H_r^2)^\nu.\label{h2norm-nu}
\end{eqnarray}
Moreover, let the subspace
$\left(H_r^2\right)^\nu\subset \left(\tilde H_r^2\right)^\nu$ be
$$\left(H_r^2\right)^\nu:=\left\{f\in \left(\tilde H_r^2\right)^\nu: \left. f\right|_M \in (\Gamma^{\omega})^\nu=(\Gamma^{\omega}(M,TM))^\nu\right\},$$
equipped with the induced Hardy product.
\end{defi}

In what follows, we shall use the following ``vector-valued" version of Boutet de Monvel's theorem. It is obtained verbatim from its proof given by Stenzel \cite{stenzel-poisson}(or by Lebeau \cite{Leb18}) using the Heat kernel on sections of the tangent bundle (see \cite{gilkey}[Section 1.6.4, P.54]) instead of on functions. Indeed, recalling the definition of $\alpha:TM\rightarrow T^*M$ (see the beginning of Section \ref{section-geom}), the kernel
$$
K(t,x,y):=\sum_{k\geq 0}e^{-t\tilde\lambda_k}{\bf e_k}(x)\otimes \alpha({\bf e_k}(y))
$$
satisfies the elliptic system $(-2\partial_t^2+\I\otimes \De_{T^*M}+\Delta_{TM}\otimes \I)K=0$. Hence, it is analytic on $\mathbb{R}_+^*\times M\times M$ \cite{treves-book-analytic}[4.1.4].


\begin{thm}\label{theo-boutet}\cite{BdM78}
Let $u\in L^2(M,TM)$ with the expansion (\ref{generalFourier}). For $0<r<r_*$, $u$ extends to a section $\ti u\in H_r^2$ if and only if
	\beq\label{exp-decay-thm}
	\sum_{i\geq 0} \left|\hat u_i\right|^2 e^{2 r\tilde\lam_i} (1+\tilde\lam_i)^{-\frac{n-1}{2}}<+\infty.
	\eeq
\end{thm}

\smallskip

In the sequel, the extension $\ti u\in H_r^2$ of $u\in L^2(M,TM)$ will be still denoted by $u$ since they are identified through the sequence of coefficients $(u_i)_{i\in\N}$ satisfying (\ref{exp-decay-thm}). 

 \begin{defi}\label{defi_norm_L2pon}
For $u\in L^2(M,TM)$, let the {\bf $\bf L^2-$norm} be
\begin{equation}\label{norme-L2}
  \|u\|_{L^2}:=  \left(\int_{M}\la u(x),u(x)\ra_{g} d{\rm vol}(x)\right)^{\frac12}.
 \end{equation}
 For $u\in H_r^2$, $0<r<r_*$, let the (exponentially) {\bf weighted $\bf L^2-$norm} be
\begin{equation}\label{norme-ponderee}
\|u\|_{r}:=\left(\sum_{i\geq 0}\left| \hat u_i \right|^2 e^{2 r\tilde\lam_i} (1+\tilde\lam_i)^{-\frac{n-1}{2}}\right)^{\frac12},\quad u=\sum_{i\geq 0} \hat u_i{\bf e}_i\in H^2_{r}.
\end{equation}
For the vector of sections, the norms are naturally defined as in \re{h2norm-nu}.
\end{defi}

\begin{remark}\label{rmk_L2} For $u=(u_l)_{1\leq l\leq\nu}\in (\Ga^\om)^\nu$, we have $u\in L^2(M,TM)^\nu$ with
$$\|u\|_{L^2}=  \left(\sum_{1\leq l\leq \nu} \int_{M}\la u_l(x),u_l(x)\ra_{g}d{\rm vol}(x)\right)^{\frac12}\lesssim  \left(\sum_{1\leq l\leq \nu}\|u_l\|^2_{C^0}\right)^{\frac12}  =: \|u\|_{C^0}.$$
\end{remark}

In the following, let $\CI$ be a closed sub-interval of $]0,r_*[$.
The inequality with ``$\lesssim$" means boundedness from above by a positive constant uniform on $\CI$ (independent of the choice of $r\in \CI$) depending only on the manifold $(M_{r_*},\ka)$, and the inequality with ``$\simeq$" means such boundedness from above and below.

\begin{prop}\label{prop_estim-norm}
 The following assertions hold true for any $r\in {\CI}$.
\begin{enumerate}
\item [(i)] For every $v\in H_r^2$, $\|v\|_{r}\simeq \|v\|_{H^2_{r}} $.

\smallskip

\item [(ii)] Given $v=\sum_{i\geq 0} \hat v_i{\bf e}_i\in H^2_{r}$, the coefficients $\{\hat v_i\}_{i\geq 0}$ satisfy that
\begin{equation}\label{decay_coeff}
|\hat v_i|\leq \|v\|_{r}e^{- r\lam_j} (1+\lam_j)^{\frac{n-1}{4}},\qquad \forall \ i\in I_j.
\end{equation}
On the other hand, if the sequence $\{\hat v_i\}_{i\geq 0}$ satisfies that
\begin{equation}\label{decay_coeff-inverse}
|\hat v_i|\leq \CD e^{- r\lam_i} (1+\lam_i)^{\frac{n-1}{4}},\qquad \forall \ i\in I_j
\end{equation}
for some constant $\CD>0$, then $v=\sum_{i\geq 0} \hat v_i{\bf e}_i\in H^2_{r'}$ for any $r'\in ]0,r[$ with
\begin{equation}\label{esti-norm-vr'}
\|v\|_{r'}\lesssim \frac{\CD}{(r-r')^{\frac{n}{2}}}.
\end{equation}

\smallskip

\item [(iii)] For every $v\in H_r^2$, we have
\begin{eqnarray}
\|v\|_{r'}&\lesssim&\|v\|_{C^0,r},\qquad \forall \  r'\in ]0,r[. \label{supnormequiv-1}\\
\|v\|_{C^0,r'}&\lesssim&\frac{\|v\|_{r}}{(r-r')^{3n}},\qquad \forall \  r'\in [0,r[,\label{supnormequiv-2}
\end{eqnarray}
	
\smallskip

\item [(iv)] There are natural continuous embeddings $H^2_{r}\hookrightarrow H^2_{r'}$, $r'\in ]0,r[$.
	\end{enumerate}
\end{prop}

\noindent
{\it Proof of (i).}
In view of the inequality (1.7) in Theorem 1.1 of \cite{Leb18}, we have that, for any $ r\in ]0,r_*[$, there exists a constant $c_{r}>0$ such that
$$c_{r}^{-1}\|v\|_{H_r^2}\leq \|v\|_{r}\leq c_{r}\|v\|_{H_r^2},\qquad \forall \ v\in H_r^2.$$
Then, according to Proposition 5.4 of \cite{Leb18}, the constant $c_r$ in the above inequality can be chosen uniformly for $r$ in the sub-interval $\CI\subset ]0,r_*[$, which means the uniform equivalence between Hardy norm and weighted $L^2-$norm.

\smallskip

\noindent
{\it Proof of (ii).} The estimate (\ref{decay_coeff}) follows immediately from Theorem \ref{theo-boutet}. On the other hand, from the definition (\ref{norme-ponderee}) of weighted $L^2-$norm $\|\cdot\|_{r}$, and, under the assumption (\ref{decay_coeff-inverse}), $v=\sum_{i\geq 0} \hat v_i{\bf e}_i$ satisfies that
$$\|v\|_{r'}^2  =   \sum_{i\geq 0}\left|\hat v_i\right|^2 e^{2r'\tilde\lam_i} (1+\tilde\lam_i)^{-\frac{n-1}{2}}
		\leq  \CD^2 \sum_{i\geq 0}  e^{-2(r-r')\tilde\lam_i}.$$	
According to the asymptotic estimate (\ref{asymp-lamdak}), there is a constant $a$, depending only on the Riemannian manifold $(M,g)$ such that $\tilde\lam_i\geq  a i^{\frac{1}{n}}$. Hence, by successive integration by parts, we obtain (\ref{esti-norm-vr'}) through
$$\sum_{i\geq 0} e^{-2(r-r')\tilde\lam_i}\leq \sum_{i\geq 0} e^{-2a(r-r') i^{\frac{1}{n}}}\lesssim \int_0^{+\infty}  e^{-2a(r-r') t^{\frac{1}{n}}} dt \lesssim (r-r')^{-n}.$$

\smallskip

\noindent
{\it Proof of (iii).} By the definition of Hardy norm in \re{h2norm}, as well as (i), we have
$$\|v\|_{r'}\simeq \|v\|_{H^2_{r'}}\leq   \CS_{\CI}  \sup_{z\in \partial M_{r'}}|v(z)|_{\kappa}\leq   \CS_{\CI}\|v\|_{C^0,r},\quad \forall \  r'\in ]0,r[$$
where $\CS_{\CI}:=\sup_{r\in\CI} \CS_{r}$ with $\CS_{r}$ the ``surface size" of $\partial M_{r}$.

According to \cite{GLS96}[Proposition 2.1],  for every $0\leq r'<r_*$, $\|{\bf e}_i\|_{C^0,r'}\lesssim (1+\tilde\lam_i)^{n+1}e^{r'\tilde\lam_i}$,\footnote{We also mention an improved estimate due to Zelditch \cite{Zel12}[Corollary 3] which allows to replace $(1+\tilde\lam_i)^{n+1}$ by $(1+\tilde\lam_i)^{\frac{n+1}{4}}$.}
where the $\|\cdot\|_{C^0,r'}-$norm with $r'=0$ means the $\|\cdot\|_{C^0}-$norm defined in \re{normsCk}.	
Therefore, by \re{decay_coeff}, for $v\in H_r^2$, for $r'\in [0,r[$,
$$\|v\|_{C^0,r'}  \leq \sum_{i\geq 0} |\hat v_i|\|{\bf e}_i\|_{C^0,r'}
\leq\|v\|_{r}  \sum_{i\geq0} e^{- (r-r')\tilde\lam_i} (1+\tilde\lam_i)^{\frac{5n+3}{4}}
\lesssim\frac{\|v\|_{r}}{(r-r')^{5n}}. $$
Indeed, in view of the asymptotic estimate (\ref{asymp-lamdak}), we have
\beq\label{sum}
\sum_{i\geq0}(1+\tilde\lam_i)^{\frac{5n+3}{4}} e^{-(r-r')\tilde\lam_i}\lesssim \sum_{i\geq1} j^{\frac{1}{n}\cdot \frac{5n+3}{4}} e^{-a(r-r') i^{\frac{1}{n}}},
\eeq
and, for the function $h_b:t\mapsto t^{\frac{c}{n}}e^{-bt^{\frac{1}{n}}}$, $b>0$ and $c=\frac{5n+3}{4}$, we have
$$\max_{t\in \R_+^*}h_b(t)=  h_b\left( \left(\frac{nc}{b}\right)^{n}\right) = \frac{d_n}{b^{\frac{5n}{4}+\frac{3}{4}}},$$
with a constant $d_n>0$ depending only on $n$. By successive integration by parts on $[0,+\infty[$, the sum \re{sum} is bounded as
$$
\sum_{i\geq1} i^{\frac{1}{n}\cdot \frac{5n+3}{4}} e^{-a(r-r') i^{\frac{1}{n}}}\lesssim \max_{t\in \R_+^*}h_{a(r-r')}(t) + \int_{0}^{+\infty} t^{\frac{5}{4}+\frac{3}{4n}} e^{-a(r-r')t^{\frac{1}{n}}} dt \lesssim \frac{1}{(r-r')^{3n}}.$$

\smallskip

\noindent
{\it Proof of (iv).}
By Definition \ref{def-Hardy} and Remark \ref{rmk_L2}, any $f=\sum_{i\geq 0} \hat f_i{\bf e}_i\in H^2_{r}$, $ r\in \CI$, is an element of $L^2(M,TM)$ when it is restricted to $M$.
Moreover, for $r'\in ]0,r[$,
$\|f\|_{r'}\leq \|f\|_{r}$.
 As a consequence of the assertion (i), we have $$\|f\|_{H^2_{r'}}\lesssim \|f\|_{r'}\leq \|f\|_{r}\lesssim \|f\|_{H^2_{r}},$$
which implies the continuous injection $H_{r}^2 \hookrightarrow H_{r'}^2$.\qed

\smallskip

\begin{lemma}\label{lem_norm-1}
For $r\in \CI$ and $r'\in [0,r[$, and $v\in H_{r}^2$,
we have, for $\tilde r:=\frac{r+r'}{2}$,
$$\|v\|_{C^1,r'}\lesssim  \frac{\|v\|_{C^0,\tilde r}}{r-r'} \lesssim  \frac{\|v\|_{r}}{(r-r')^{3n+1}}.$$
\end{lemma}

\proof Recall that there is a coordinate chart $\left\{\left(W_i,z^{(i)}\right)\right\}$ of $M_{r_*}$, such that $z^{(i)}\in \Delta_1^n$ and that in one of these charts, $\tilde v$ denotes the expression of $v$.
Recalling the definition \re{rho} of $\rho$ and its properties from \rt{GS-thm}, let us define
\begin{eqnarray*}
\|D(\rho\circ (z^{(i)})^{-1})\|_{0} &:=&	\sup_{\mathfrak z\in z^{(i)}(M_{r_*}\cap W_i)}\sup_{\substack{\zeta\in \mathbb{C}^n,\\|\zeta|\leq 1}}|D(\rho\circ (z^{(i)})^{-1})(\mathfrak z)\zeta|,\\
\delta&:=&  \frac{r-r'}{2\|D(\rho\circ (z^{(i)})^{-1})\|_{0}} \, , \\
z^{(i)}(M_{r'}\cap W_i)_\delta&:=& \left\{z\in\C^n:|z-z^{(i)}(q)|<\delta \ {\rm for \; some} \; q\in M_{r'}\cap W_i\right\}.
\end{eqnarray*}
Assume that $\delta$ is small enough so that the $\delta$-neighborhood of $W_i\cap M_{r'}$ is still in $W_i\cap M_{r_*}$~:
$$\left(z^{(i)}\right)^{-1}\left(z^{(i)}(M_{r'}\cap W_i)_\delta\right)\subset M_{r_*}\cap W_i.$$
Let us devise a Cauchy-like estimate relative to K\"ahler norm. First of all, we can assume that, on the trivialization,
\begin{equation}\label{norm-triv}
|\zeta|_{z,\kappa}^2\simeq \sum_{1\leq i\leq n}|\zeta_i|^2,\qquad \forall \ z\in\Delta_r^n, \quad  \zeta\in \mathbb{C}^n.
\end{equation}
Let $\tilde v$ be a holomorphic vector field on $\Delta_r^n$.
For $z\in \Delta_{r-\delta}^n$ and $\zeta$ in the unit ball of $\mathbb{C}^n$ (i.e., $\sum_{i=1}^n |\zeta_i|^2=1$), we have, through Cauchy-Schwarz inequality, that
\begin{equation}\label{D-tildev}
|D\tilde v(z)\zeta|_{z,\kappa}^2\lesssim \sum_{1\leq i\leq n} \left|\sum_{1\leq k\leq n}\frac{\partial \tilde v_i}{\partial z_k}(z)\zeta_k\right|^2
\lesssim\sum_{1\leq i\leq n}\sum_{1\leq k\leq n} \left|\frac{\partial \tilde v_i}{\partial z_k}(z)\right|^2.
\end{equation}
Since, through Cauchy estimate, we have that
$$\sum_{1\leq k\leq n}\left|\frac{\partial \tilde v_i}{\partial z_k}(z)\right|^2 \leq  \frac{n}{\delta^2} \sup_{w\in \Delta_r^n}| \tilde v_i(w)|^2\leq  \frac{n}{\delta^2} \sup_{w\in \Delta_r^n} \sum_{1\leq i\leq n} | \tilde v_i(w)|^2.$$
Hence, by \re{norm-triv} and (\ref{D-tildev}) we have~:
$$|D\tilde v(z)\zeta|_{z,\kappa}^2\lesssim \frac{n^2}{\delta^2}\sup_{w\in \Delta_r^n}\sum_{1\leq i\leq n} | \tilde v_i(w)|^2\leq \frac{n}{\delta^2}\sup_{w\in \Delta_r^n}|\tilde v|_{w,\kappa}^2.$$
 According to the definition of norm (\ref{norm1r}), we have, by Cauchy estimate, that
\beq\label{maj1}
\|v\|_{C^1,r'}=\max_{i} \sup_{q'\in W_i\cap M_{r'}}\|D\tilde v(z^{(i)}(q'))\|_{\kappa}
\lesssim \max_{i} \sup_{q'\in W_i\cap M_{r'}}\sup_{|z-z^{(i)}(q')|=\delta}\frac{|\tilde v(z)|_{\kappa}}{ \delta}.
\eeq
We recall that $q\in M_{r'}$ if and only if $0\leq \rho(q)<r'$.
With $\tilde r=\frac{r+r'}{2}$, let us show that
\beq\label{maj2}
\max_{i} \sup_{q'\in W_i\cap M_{r'}}\sup_{|z-z^{(i)}(q')|=\delta}|\tilde v(z)|_{\kappa} \leq \|v\|_{C^0,\tilde r}.
\eeq
Applying Taylor formula of order $1$, we obtain that
\begin{eqnarray*}
\left|\rho((z^{(i)})^{-1}(z))-\rho(q')\right|&\leq&\|D(\rho\circ (z^{(i)})^{-1})\|_{0} \cdot \left|z-z^{(i)}(q')\right|\\
&\leq & \delta  \|D(\rho\circ (z^{(i)})^{-1})\|_0 \ =  \ \frac{r-r'}{2} \ =  \ \tilde r -r'.
\end{eqnarray*}
Hence,
\begin{equation*}\label{maj3}
\rho((z^{(i)})^{-1}(z))\leq \rho(q')+ \left|\rho((z^{(i)})^{-1}(z))-\rho(q')\right|  \leq r'+ (\tilde r-r')=\tilde r.
\end{equation*}
The first estimate is shown. The second follows from the latter together with \re{supnormequiv-1}.\qed

\smallskip

%
%

As a corollary of Proposition \ref{propMoser-anal}, we have
\begin{cor}\label{cor_s1-hardi}
Given $w,v\in H^2_r$, there exists $s_1(w,v)\in H^2_{r'}$ for any $0\leq r'< r$, such that (\ref{eq_lemMoser1}) holds with $s_1(w,0)=s_1(0,v)=0$. Moreover,
we have,
 $$\|s_1(w,v)\|_{r'}\lesssim\frac{\|w\|_{r}\|v\|_{r}}{(r-r')^{6n+1}}.$$
 \end{cor}
 \proof According to Proposition \ref{propMoser-anal}, there exists such $s_1(w,v)\in \Ga_{\frac{r+2r'}{3}}^\omega$. Then, applying Proposition \ref{prop_estim-norm}-(iii) and Lemma \ref{lem_norm-1}, we have
 \begin{eqnarray*}
\|s_1(w,v)\|_{r'}&\lesssim& \|s_1(w,v)\|_{C^0, \frac{r+2r'}{3}} \\
    &\lesssim& \|w\|_{C^1, \frac{2r+r'}{3}} \|v\|_{C^0, \frac{r+2r'}{3}}\\
     &\lesssim& \frac{\|w\|_{r}}{(r-r')^{3n+1}}\frac{\|v\|_{r}}{(r-r')^{3n}} \ = \ \frac{\|w\|_{r}\|v\|_{r}}{(r-r')^{6n+1}}.\qed
 \end{eqnarray*}

\section{Analytic rigidity of $G-$action by analytic isometries}\label{secKAManal}

In this section, we shall prove Theorem \ref{thmmain} and \ref{thm-geom0} in the analytic setting, i.e., to show the conjugacy equation (\ref{rig-smooth}) with some $W\in \Ga^\omega(M,TM)$.

 Let $M$ be an analytic compact Riemannian manifold, and let $\pi$ be a
$G-$action by analytic isometries. Let $\pi_{P_0}$ be a $G-$action by analytic diffeomorphisms with $\pi_{P_0}(\gamma)=\Exp\{P_0(\gamma)\}\circ\pi(\gamma)$ for $\gamma\in \CS$, where $P_0:\CS\to H_{r_0}^2$, $r_0\in ]0,r_*[$, with
$$\|P_0\|_{\CS,r_0}=\left(\sum_{\gamma\in\CS}\|P_0(\gamma)\|_{r_0} \right)^\frac12=\varepsilon_0.$$
Assume that the hypotheses in Theorem \ref{thmmain} (for $\pi$ and $\pi_{P_0}$) or that in Theorem \ref{thm-geom0} (for $\pi$) on $M$ are satisfied. If $\varepsilon_0$ is sufficiently small such that (\ref{smallvarepsilon0}) is satisfied, then as shown in Section \ref{secKAMsmooth}, $\pi_{P_0}$ is smoothly conjugate to $\pi$.

Let ${\CI}:=[\frac{r_0}{2}, r_0]\subset ]0,r_*[$. Let $c_\lesssim$ be the constant as in (\ref{smallvarepsilon0}), and assume further that it is greater than all implicit constant in the inequalities with ``$\lesssim$" in Section \ref{sec_Grauert_Hardy},
depending on the manifold $(M_{r_*},\ka)$ and uniform in ${\CI}$.
Assume that $\varepsilon_0$ is sufficiently small such that (\ref{smallvarepsilon0}) is satisfied with this new $c_\lesssim$ and
$r_0^{-(\tau+6n+1)}\varepsilon_0^\frac{1}{60} <1$.

 Define the sequences of radii $\{r_m\}$ by
\begin{equation}\label{seqs-r}
r_{m+1}:= r_m- \frac{r_0}{2^{m+2}}.
\end{equation}
It is easy to see that, for every $m\in\N$, $r_m\in \CI=[\frac{r_0}{2}, r_0]$ and $r_m\to \frac{r_0}{2}$ as $m\to\infty$.
With the sequence $\varepsilon_m=\varepsilon_0^{(\frac54)^m}$, same as (\ref{esti-Pm}) in the smooth case, we have
\begin{equation}\label{vareps-m_small}
\frac{\varepsilon_m^{\frac{1}{60}}}{(r_{m}-r_{m+1})^{\tau+6n+1}}<\frac{\varepsilon_m^{\frac{1}{60}}}{(r_{m+1}-r_{m+2})^{\tau+6n+1}}<1.
\end{equation}
Indeed, according to the definition of $\{r_m\}$, we see that
$$\frac{1}{(r_{m+1}-r_{m+2})^{\tau+6n+1}}=\frac{2^{(m+3)(\tau+6n+1)}}{r_0^{\tau+6n+1}}=\left(\frac{8}{r_0}\right)^{\tau+6n+1}2^{m(\tau+6n+1)}.$$
Then, under the assumption (\ref{smallvarepsilon0}), we have
$$
\ln\left(\frac{1}{(r_{m+1}-r_{m+2})^{\tau+6n+1}} \right)= (\tau+6n+1)\left(\ln\left(\frac{8}{r_0}\right)  + m \ln 2 \right)< \frac{|\ln(\varepsilon_{0})|}{60}\left(\frac54\right)^m,
$$
which implies (\ref{vareps-m_small}).

\subsection{Sequence of $G-$action by analytic diffeomorphisms}\label{sec_seq-anal}
\begin{prop}\label{prop_scheme-anal}
With $\widehat\CP_0:=\CP_0$, there exist $\{\widehat\CP_m\}\subset{\bf G}_\pi $ with $\widehat\CP_m=(\widehat P_m(\gamma))_{\gamma\in \CS}\in (H_{r_m}^2)^k$ satisfying $\|\widehat\CP_m\|_{r_m}\leq  \varepsilon_m$, and
$\widehat w_m\in H_{r_{m+1}}^2$ with $\|\widehat w_m\|_{r_{m+1}}\leq  \varepsilon_m^\frac67$, such that
\beq\label{conj_hat}
\Exp\{\widehat w_m\}^{-1}\circ\Exp\{\widehat P_m(\gamma)\}\circ\pi\circ\Exp\{\widehat w_m\}= \Exp\{\widehat P_{m+1}(\gamma)\}\circ\pi,\quad \gamma\in \CS.  \eeq
\end{prop}

The rest of this subsection is devoted to the proof of Proposition \ref{prop_scheme-anal}.
Define the sequence $\{K_m\}\subset \N$ as $K_m=8m$.
Then we have $\varepsilon_{K_m}\leq\varepsilon_m^{4^m}$.
Indeed, for $m=0$, we have
$\varepsilon_{K_0}=\varepsilon_0$,
and, if $\varepsilon_{K_m}\leq\varepsilon_m^{4^m}$ for some $m\in \N$, then, noting that $(\frac54)^{7}>4$, we have
$$\varepsilon_{K_{m+1}}=\varepsilon_{K_{m}}^{(\frac54)^{K_{m+1}-K_m}}\leq \varepsilon_m^{(\frac54)^{8}\cdot 4^m}\leq \varepsilon_m^{\frac54\cdot 4^{m+1}}= \varepsilon_{m+1}^{4^{m+1}}.$$
Let $\{\CP_m\}$ be the sequence in ${\bf G}_\pi$ obtained in the smooth KAM scheme satisfying $\|\CP_m\|_{C^0}<\varepsilon_m$ (see Section \ref{sec-smoothPm}). Then the subsequence $\{\CP_{K_{m}}\}$ satisfies that
$$\|\CP_{K_{m}}\|_{L^2}\lesssim \|\CP_{K_{m}}\|_{C^0}< \varepsilon_{K_m}\leq\varepsilon_m^{4^m}.$$

\begin{lemma}For every $m\in \N$, $\CP_{K_{m}}\in (H_{r_m}^2)^k$ with $\|\CP_{K_{m}}\|_{r_m}\leq\varepsilon_m$.
\end{lemma}
\proof For $m=0$, we have $\|\CP_{K_{0}}\|_{r_0}=\|P_0\|_{S,r_0}=\varepsilon_0$. Assume that $\|\CP_{K_{m}}\|_{r_m}\leq\varepsilon_m$ for $m\in \N$,
let us show $\|\CP_{K_{m+1}}\|_{r_{m+1}}\leq\varepsilon_{m+1}$, which proves the lemma by induction.

Given $m\in \N$, define the intermediate radii between $r_m$ and $r_{m+1}$ by
\beq\label{radii_intermed}
\tilde r_m:= \frac{r_{m}+r_{m+1}}{2}, \qquad r_{m,l}=\frac{l\tilde r_{m}+(32-l)r_{m}}{32},\quad l=0,1,\cdots, 32.
\eeq
It is easy to verify that $\tilde r_m=r_{m,32}<\cdots  <r_{m,0}=r_{m}$ and for $l=0,1,\cdots, 31$
$$r_{m,l} - r_{m,l+1}= \frac{r_{m}-\tilde r_{m}}{32}=\frac{r_{m}- r_{m+1}}{64}.$$

 As shown in Section \ref{sec_cohomoeq-smooth}, there exist
 $$w_{K_m+q}\in \IM d_0^*\cap \bigoplus_{\lam_j\leq N_{K_{m}+q}} E_{\lam_j}, \quad q=0,1,\cdots, 7,$$
satisfying $d_0w_{K_m+q}=(\D_0\circ\Pi_{N_{K_{m}+q}})\CP_{K_{m}+q}$ for $ N_{K_{m}+q}=\varepsilon_{K_{m}+q}^{-\frac{1}{8(\tau+n+1)}}$, such that, for $\gamma\in \CS$,
 \beq\label{conj_Kmq}
 \Exp\{w_{K_m+q}\}^{-1}\circ\Exp\{P_{K_{m}+q}(\gamma)\}\circ\pi\circ\Exp\{w_{K_m+q}\}= \Exp\{P_{K_m+q+1}(\gamma)\}\circ\pi.\eeq
In view of Remark \ref{rmk_w-analytic}, every $w_{K_m+q}$ is real analytic on $M$ with finitely many ``Fourier modes". As such they extend holomorphically to any Grauert tube, and in particular, they extend to $M_{r_*}$.
 Assume that $\|\CP_{K_{m}+q}\|_{r_{m,4q}}\leq   \varepsilon_{m}^{1-\frac{q}{56}}$ for some $q=0,1,\cdots, 7$.
 Let us show that $\| \CP_{K_m+q+1}\|_{r_{m,4(q+1)}}\leq \varepsilon_{m}^{1-\frac{q+1}{56}}$.
In view of (\ref{cohomo_esti-L2proj}) in Lemma \ref{lem_sol_cohomo-d0}, we have
\begin{eqnarray}
\|w_{K_m+q}\|_{r_{m,4q+1}} &\lesssim& \left(\sum_{\lam_j\leq N_{K_m+q}}(1+\lam_j)^{2\tau-\frac{n-1}{2}}e^{2r_{m,4q+1}\lam_j} \|\PP_j \CP_{K_{m}+q}\|^2_{L^2}\right)^\frac12 \nonumber\\
&=& \left(\sum_{\lam_j\leq N_{K_m+q}}(1+\lam_j)^{2\tau}e^{-2(r_{m,4q}-r_{m,4q+1})\lam_j}\|\PP_j \CP_{K_{m}+q}\|^2_{r_m,4q}\right)^\frac12 \nonumber\\
&\lesssim& \frac{\|\CP_{K_{m}+q}\|_{r_m,4q}}{(r_m-r_{m+1})^\tau},\label{esti-wKm+q}
\end{eqnarray}
since the maximum of the function $x\to (1+x)^{2\tau}e^{-\frac1{32}(r_m-r_{m+1})x}$ on $\R_+$ is bounded by $e^{\frac1{32}(r_m-r_{m+1})}(\frac{64\tau}{r_m-r_{m+1}})^{2\tau}$.
 Recalling (\ref{Pm+1}), we have
\begin{eqnarray}
\CP_{K_m+q+1}&=& \Pi_{N_{K_m+q}}^\perp \CP_{K_m+q} +(\Pi_{N_{K_m+q}}\circ\D_0^\perp)\CP_{K_m+q}\label{analytic-PKM+q+1}\\
& &  - \,  s_1(w_{K_m+q},\CP_{K_m+q+1})+s_1(\CP_{K_m+q}, \pi_* w_{K_m+q}),\nonumber
 \end{eqnarray}
where, according to Proposition \ref{propMoser-anal}, we have
\begin{eqnarray*}
\|s_1(w_{K_m+q},\CP_{K_m+q+1})\|_{C^0,r_{m,4q+3}}  &\lesssim& \|w_{K_m+q}\|_{C^1,r_{m,4q+2}}  \|\CP_{K_m+q+1}\|_{C^0,r_{m,4q+3}} \\
&\lesssim& \frac{\|w_{K_m+q}\|_{r_{m,4q+1}}}{(r_m-r_{m+1})^{3n+1}}  \|\CP_{K_m+q+1}\|_{C^0,r_{m,4q+3}}\\ &\lesssim&  \frac{\|\CP_{K_m+q}\|_{r_{m,4q}}}{(r_m-r_{m+1})^{\tau+3n+1}}  \|\CP_{K_m+1}\|_{C^0,r_{m,4q+3}},\\
\|s_1(\CP_{K_m+q}, \pi_* w_{K_m+q})\|_{C^0,r_{m,4q+3}}&\lesssim&\|\CP_{K_m+q}\|_{C^1, r_{m,4q+2}}  \|w_{K_m+q}\|_{C^0, r_{m,4q+3}}\\
&\lesssim&\frac{\|\CP_{K_m+q}\|_{r_{m,4q}}}{(r_m-r_{m+1})^{3n+1}} \frac{\|w_{K_m+q}\|_{r_{m,4q+1}}}{(r_m-r_{m+1})^{3n}}\\
&\lesssim&\frac{\|\CP_{K_m+q}\|^2_{r_{m,4q}}}{(r_m-r_{m+1})^{\tau+6n+1}} .
\end{eqnarray*}
Moreover, we have
$$\| \Pi_{N_{K_m+q}}^\perp \CP_{K_m+q} +(\Pi_{N_{K_m+q}}\circ\D_0^\perp)\CP_{K_m+q}\|_{C^0, r_{m,4q+3}}  \lesssim \frac{\|\CP_{K_m+q}\|_{r_{m,4q}}}{(r_m-r_{m+1})^{3n}}. $$
Recalling that $\|\CP_{K_{m}+q}\|_{r_{m,4q}}<   \varepsilon_{m}^{1-\frac{q}{56}}$, in view of (\ref{vareps-m_small}), we have
$$\frac{\|\CP_{K_m+q}\|_{r_{m,4q}}}{(r_m-r_{m+1})^{\tau+3n+1}} \leq \varepsilon_{m}^{\frac{59}{60}-\frac{q}{56}} $$
Hence, taking the $C^0, r_{m,4q+3}$-norm of \re{analytic-PKM+q+1}, there exists a constant $c>0$ such that
\begin{eqnarray*}
\left(1-c  \varepsilon_{m}^{\frac{59}{60}-\frac{q}{40}} \right)\|\CP_{K_m+q+1}\|_{C^0,r_{m,4q+3}}\leq c\left(\frac{\|\CP_{K_m}\|_{r_{m}}}{(r_m-r_{m+1})^{3n}} +\frac{\|\CP_{K_m}\|^2_{r_{m}}}{(r_m-r_{m+1})^{\tau+6n+1}}\right)\leq 2 c \varepsilon_{m}^{\frac{59}{60}-\frac{q}{56}}.\end{eqnarray*}
According to (\ref{supnormequiv-1}), we have $\| \CP_{K_m+q+1}\|_{r_{m,4(q+1)}}\lesssim \| \CP_{K_m+q+1}\|_{C^0,r_{m,4q+3}}$. Then we obtain
$$\| \CP_{K_m+q+1}\|_{r_{m,4(q+1)}}\leq \varepsilon_{m}^{\frac{55}{56}-\frac{q}{56}}=\varepsilon_{m}^{1-\frac{q+1}{56}}.  $$
As $q=7$, we have $\| \CP_{K_m+8}\|_{r_{m,32}}=\| \CP_{K_{m+1}}\|_{\tilde r_{m}}\leq \varepsilon_{m}^{\frac67}$.

Now, let us apply the interpolation lemma \ref{interpol-hardy}, with
$$\| \CP_{K_{m+1}}\|_{L^2}\lesssim \varepsilon_{m+1}^{4^{m+1}}=\varepsilon_{m}^{\frac54\cdot 4^{m+1}},\quad  \| \CP_{K_{m+1}}\|_{\tilde r_{m}}< \varepsilon_{m}^{\frac67}.$$
Since $\frac{r_0}{2}<\tilde r_m<r_0$, we have
$$\frac{\tilde r_m-r_{m+1}}{\tilde r_m}=\frac{r_{m}-r_{m+1}}{2\tilde r_m}=\frac{r_0}{2^{m+3} \tilde r_m}\in \left[\frac{1}{2^{m+3}},\frac{1}{2^{m+2}}\right],
 $$
which implies that
$$\frac{r_{m+1}}{\tilde r_m}=1-\frac{\tilde r_m-r_{m+1}}{\tilde r_m} \geq 1- \frac{1}{2^{m+2}}, $$
we obtain that
$$\|\CP_{K_{m+1}}\|_{r_{m+1}}\leq  \|\CP_{K_{m+1}}\|_{L^2}^{\frac{\tilde r_m-r_{m+1}}{\tilde r_m}} \|\CP_{K_{m+1}}\|_{\tilde r_m}^{\frac{r_{m+1}}{\tilde r_m}}
\lesssim \varepsilon_{m}^{\frac54 \cdot \frac{4^{m+1}}{2^{m+3}} } \cdot \varepsilon_{m}^{\frac67(1- \frac{1}{2^{m+2}})}, $$
which implies that $\|\CP_{K_{m+1}}\|_{r_{m+1}}\leq \varepsilon_{m+1}$, since
$$\frac54\cdot\frac{4^{m+1}}{2^{m+3}} + \frac67\left(1- \frac{1}{2^{m+2}}\right)\geq \frac54\cdot\frac12+\frac67\cdot\frac34=\frac{71}{56}>\frac54. \qed$$

\medskip

\noindent{\bf Proof of Proposition \ref{prop_scheme-anal}.} Let $\hat \CP_m:=\CP_{K_{m}}\in {\bf G}_\pi$, and for $q=1,\cdots,7$,
define $\widehat w_{m,q}\in \Ga^\infty(M,TM)$ by
$$\Exp\{\widehat w_{m,q}\}= \Exp\{w_{K_m}\}\circ \Exp\{w_{K_m+1}\}\circ \cdots \circ \Exp\{w_{K_m+q}\}.$$
In particular, for $\widehat w_{m}:=\widehat w_{m,7}$, we obtain (\ref{conj_hat}) according to (\ref{conj_Kmq}).

It remains to estimate the $\|\cdot\|_{r_{m+1}}-$norm of $\widehat w_{m}$. For $q=1$, we have
$$\widehat w_{m,1}=w_{K_m}+w_{K_m+1}+s_1(w_{K_m}, w_{K_m+1}).$$
In view of (\ref{esti-wKm+q}), and recalling that  $\|\CP_{K_{m}+q}\|_{r_{m,4q}}<   \varepsilon_{m}^{1-\frac{q}{56}}$, we have $w_{K_m}\in H_{r_{m,4}}^2$ with
$$ \|w_{K_m}\|_{r_{m,4}}\leq \|w_{K_m}\|_{r_{m,1}} \lesssim  \frac{\varepsilon_{m}}{(r_m-r_{m+1})^\tau},\quad  \|w_{K_m+1}\|_{r_{m,5}}   \lesssim  \frac{ \varepsilon_{m}^{\frac{55}{56}}}{(r_m-r_{m+1})^\tau}.$$
Then, according to Corollary \ref{cor_s1-hardi}, $\widehat w_{m,1}\in H_{r_{m,8}}^2$ with
\begin{eqnarray*}
\|\widehat w_{m,1}\|_{r_{m,8}}&\lesssim&\|w_{K_m}\|_{r_{m,4}}+\|w_{K_m+1}\|_{r_{m,5}}+\frac{\|w_{K_m}\|_{r_{m,4}}\|w_{K_m+1}\|_{r_{m,5}}}{(r_m-r_{m+1})^{6n+1}}\\
&\lesssim& \frac{ \varepsilon_{m}^{\frac{55}{56}}}{(r_m-r_{m+1})^\tau}+ \frac{ \varepsilon_{m}^{\frac{111}{56}}}{(r_m-r_{m+1})^{2\tau+6n+1}} \ \lesssim \  \frac{ \varepsilon_{m}^{\frac{55}{56}}}{(r_m-r_{m+1})^\tau}.
\end{eqnarray*}
For some $q\leq 6$, let us assume that $\widehat w_{m,q}\in H^2_{r_{m,4(q+1)}}$ with
\beq\label{esti-hatwmq}
\|\widehat w_{m,q}\|_{r_{m,4(q+1)}}\lesssim  \frac{ \varepsilon_{m}^{1-\frac{q}{56}}}{(r_m-r_{m+1})^{\tau}}.
\eeq
By Corollary \ref{cor_s1-hardi},
$\widehat w_{m,q+1} =\widehat w_{m,q}+w_{K_m+q+1}+s_1(\widehat w_{m,q}, w_{K_m+q+1})\in H^2_{r_{m,4(q+2)}}$ with
\begin{eqnarray*}
\|\widehat w_{m,q+1}\|_{r_{m,4(q+2)}}&\lesssim&\|\widehat w_{m,q}\|_{r_{m,4(q+1)}}+\|w_{K_m+q+1}\|_{r_{m,4q+5}}\\
& &+ \, \frac{\|\widehat w_{m,q}\|_{r_{m,4(q+1)}}\|w_{K_m+q+1}\|_{r_{m,4q+5}}}{(r_m-r_{m+1})^{6n+1}}\\
&\lesssim&\frac{ \varepsilon_{m}^{1-\frac{q+1}{56}}}{(r_m-r_{m+1})^{\tau}}+ \frac{ \varepsilon_{m}^{2-\frac{2q+1}{56}}}{(r_m-r_{m+1})^{2\tau+6n+1}} \ \lesssim \  \frac{ \varepsilon_{m}^{1-\frac{q+1}{56}}}{(r_m-r_{m+1})^\tau}.
\end{eqnarray*}
Hence, we have (\ref{esti-hatwmq}) for $q=1,\cdots,7$, which implies that $\|\widehat w_{m,q}\|_{r_{m,4(q+1)}}\leq \varepsilon_{m}^{1-\frac{q+1}{56}}$.
In particular, for $q=7$, $r_{m,32}=\tilde r_m$ (recalling (\ref{radii_intermed})), and $\widehat w_{m,7}= \widehat w_m\in H^2_{\tilde r_{m}}$,
$$\|\widehat w_{m}\|_{r_{m+1}}\leq \|\widehat w_{m,7}\|_{\tilde r_{m}}\leq \varepsilon_{m}^{\frac67}.\qed$$

\subsection{Convergence in $\Gamma^\om(M,TM)$}\label{sec_cov-anal}

With the sequences $\{\widehat w_m\}$, $\{\widehat \CP_m\}$ constructed in Proposition \ref{prop_scheme-anal}, satisfying $\widehat w_m\in H_{r_{m+1}}^2$, $\widehat\CP_m\in (H_{r_{m}}^2)^k$ and
$$\|\widehat w_m\|_{r_{m+1}}< \varepsilon_{m}^{\frac67}, \quad \|\widehat\CP_m\|_{r_{m}}< \varepsilon_{m},$$
the proof of the analytic part of Theorem \ref{thmmain} and \ref{thm-geom0} will be completed by a convergence argument on a certain Grauert tube.

At first, let us define $\widehat W_1:=\widehat w_0+\widehat w_1+s_1(\widehat w_0,\widehat w_1)$ with $s_1(\widehat w_0,\widehat w_1)\in H_{r_3}^2$, and according to Corollary \ref{cor_s1-hardi} and (\ref{vareps-m_small}),
\begin{eqnarray*}
\|s_1(\widehat w_0,\widehat w_1)\|_{r_3}\lesssim \frac{\|\widehat w_0\|_{r_1}\|\widehat w_1\|_{r_2}}{(r_2-r_3)^{6n+1}}\lesssim \frac{\varepsilon_0^{\frac67}\varepsilon_1^{\frac67}}{(r_2-r_3)^{6n+1}}.
\end{eqnarray*}
which implies that $\|s_1(\widehat w_0,\widehat w_1)\|_{r_3} < \varepsilon_1^{\frac34}$.
Then we have $\Exp\{\widehat w_0\}\circ\Exp\{\widehat w_1\}=\Exp\{\widehat W_1\}$, and, similar to (\ref{W1}) in Section \ref{sec_cov-smooth},
$$ \Exp\{\widehat W_1\}^{-1}\circ \Exp\{\widehat P_0(\gamma)\} \circ \pi(\gamma) \circ\Exp\{\widehat W_1\}= \Exp\{\widehat P_2(\gamma)\} \circ \pi(\gamma),\quad \gamma\in \CS. $$
We also have
$$\|\widehat W_1-\widehat w_0\|_{r_{3}}\leq \|\widehat w_1\|_{r_{2}}+\|s_1(\widehat w_0,\widehat w_1)\|_{r_{3}}<\varepsilon_1^{\frac67}+\varepsilon_1^{\frac34}<2\varepsilon_1^{\frac34}.$$

Assume that there exists $\widehat W_m\in H^2_{r_{m+2}}$ with $\|\widehat W_m\|_{r_{m+2}}<\varepsilon_0^{\frac34}$ such that
$$\Exp\{\widehat W_{m}\}^{-1}\circ \Exp\{\widehat P_0(\gamma)\} \circ \pi(\gamma) \circ\Exp\{\widehat W_{m}\}=\Exp\{\widehat P_{m+1}(\gamma)\} \circ \pi(\gamma),\quad \gamma\in \CS.$$
Let $\widehat W_{m+1}:=\widehat W_m+\widehat w_{m+1}+s_1(\widehat W_m, \widehat w_{m+1})$, with $s_1(\widehat W_m, \widehat w_{m+1})\in H^2_{r_{m+3}}$, and according to Corollary \ref{cor_s1-hardi},
$$
\|s_1(\widehat W_m, \widehat w_{m+1})\|_{r_{m+3}}\lesssim\frac{\|\widehat W_m\|_{r_{m+2}} \|\widehat w_{m+1}\|_{r_{m+2}}}{(r_{m+2}-r_{m+3})^{6n+1}}
 <  \frac{\varepsilon_0^{\frac34} \varepsilon_{m+1}^{\frac67} }{(r_{m+2}-r_{m+3})^{6n+1}},
$$
which implies, through (\ref{vareps-m_small}), that $\|s_1(\widehat W_m, \widehat w_{m+1})\|_{r_{m+3}} < \varepsilon_{m+1}^{\frac34}$.
Similar to (\ref{Wm+1}),
$$\Exp\{\widehat W_{m+1}\}^{-1}\circ \Exp\{\widehat P_0(\gamma)\} \circ \pi(\gamma) \circ\Exp\{\widehat W_{m+1}\}= \Exp\{\widehat P_{m+2}(\gamma)\} \circ \pi(\gamma).$$
We also have
\beq\label{error-W-m}
\|\widehat W_{m+1}-\widehat W_m\|_{r_{m+3}}\leq\|\widehat w_{m+1}\|_{r_{m+2}} + \|s_1(\widehat W_m,\widehat w_{m+1})\|_{r_{m+3}}
\leq\varepsilon_{m+1}^{\frac67}+\varepsilon_{m+1}^{\frac34} \leq 2 \varepsilon_{m+1}^{\frac34}.
\eeq

As $m\to\infty$, $r_{m}\to \frac{r_0}{2}$, then we have the convergence of $\{\widehat W_m\}$ in $H_{\frac{r_0}{2}}^2$ from (\ref{error-W-m}), and for every $m\in \N^*$,
$$ \|\widehat W_{m+1}\|_{r_{m+2}}\leq \|\widehat w_0\|_{r_1} + \sum_{j=0}^m \|\widehat W_{j+1}-\widehat W_j\|_{r_{j+3}}<\varepsilon_0^{\frac67}+2\sum_{j=0}^m  \varepsilon_{j+1}^{\frac34}<\varepsilon_0^{\frac34}.$$
Hence, for the limit $\widehat W:=\lim_{m\to\infty}\widehat W_m\in H_{\frac{r_0}{2}}^2$, we have $\|\widehat W\|_{\frac{r_0}2}<\varepsilon_0^{\frac34}$. Since
$$\|\widehat P_m\|_{S,\frac{r_0}2}\leq \|\widehat P_m\|_{S,r_m}<\varepsilon_m\to 0,$$
we have $\Exp\{\widehat W\}^{-1} \circ \pi_0(\gamma)\circ \Exp\{\widehat W\}= \pi(\gamma)$, for every $\gamma\in \CS$.

\appendix

\section{Proof of Proposition \ref{prop_s1-smooth}.}\label{app_proof-s1smooth}

According to Appendix of \cite{Mos69},
the vector field $s_1(w,v)$ has the following expression, in local chart~:
\begin{equation}\label{expression-s1}
s_1(w,v)(x)= \underbrace{\left(w\left(x+v(x)+\phi(x,v(x))\right)-w(x)\right)}_{=:\Psi_{w,v}(x)} +\underbrace{\varrho(x,v(x),w(x))}_{=:\Upsilon_{w,v}(x)},
\end{equation}
where $\phi=\phi(x,\xi)$ and $\varrho=\varrho(x,\xi,\eta)$ are $C^\infty-$vector functions with
\begin{equation}\label{property-phi-rho}
\phi(x,0)=\phi_\xi (x,0)=0,\qquad \varrho(x,0,\eta)=\varrho(x,\xi,0)=0,\quad |\varrho_\xi|\lesssim |\eta|.
\end{equation}
\begin{remark}\label{rem-0-deriv}
We emphasize the above property of $\phi$ implies that, for all non negative appropriate multi-indices,  $(\partial_x^{P'}\partial_{\xi}\phi)(x,0)=0$ and $(\partial_x^{P'}\phi)(x,0)=0$ for $P'\in \N^n$.
As to the property of $\varrho$, it implies that $(\partial_x^{P'}\partial^{P''}_{\xi}\varrho)(x,\xi,0)=0$ and $(\partial_x^{P'}\partial^{P'''}_{\eta}\varrho)(x,0,\eta)=0$ for $P', P'', P'''\in \N^n$.
\end{remark}

\begin{lemma}
For $w\in\Gamma^1$, $v\in \Gamma^0$ with $\|w\|_{C^1}$ and $\|v\|_{C^0}$ sufficiently small, we have
\begin{equation}\label{Psi-C0}
\|\Psi_{w,v}\|_{C^0}\lesssim \|w\|_{C^1}\|v\|_{C^0}.
\end{equation}
For $w_1,w_2,v\in \Gamma^\infty(M,TM)$ with $\|w_1\|_{C^1}$, $\|w_2\|_{C^0}$ and $\|v\|_{C^0}$ sufficiently small,
\begin{equation}\label{Psi-C0-decomp}
\|\Psi_{w_1+w_2,v}\|_{C^0}\lesssim \|w_1\|_{C^1}\|v\|_{C^0}+\|w_2\|_{C^0}.\end{equation}
\end{lemma}
\proof In view of the expression of $\Psi_{w,v}$ in \re{expression-s1}, we have $\|\Psi_{w,v}\|_{C^0}\lesssim \|w\|_{C^0}$. By writing $\Psi_{w,v}$ as
\begin{equation}\label{Psi-integral}
\Psi_{w,v}(x)=\int_{0}^1 w^{(1)}\left(x+t(v(x)+\phi(x,v(x)))\right) \left(v(x)+\phi(x,v(x))\right) dt,
\end{equation}
we have (\ref{Psi-C0}). Since $\Psi_{w_1+w_2,v}=\Psi_{w_1,v}+\Psi_{w_2,v}$, we obtain (\ref{Psi-C0-decomp}).\qed

\smallskip

 According to (\ref{property-phi-rho}), for $w_1,w_2,v\in \Gamma^\infty(M,TM)$ with $\|w_1\|_{C^0}$, $\|w_2\|_{C^0}$ and $\|v\|_{C^0}$ sufficiently small, we have
\begin{eqnarray*}
|\Upsilon_{w_1+w_2,v}(x)|&=&|\varrho(x,v(x),w_1(x)+w_2(x))-\varrho(x,0,w_1(x)+w_2(x))|\\
&\leq& \sup_{x,\xi}\left|\left.\varrho_\xi(x,\xi,\eta)\right|_{\eta=w_1(x)+w_2(x)}\right|\|v\|_{C^0} \ \lesssim \ (\|w_1\|_{C^0} +\|w_2\|_{C^0}) \|v\|_{C^0}.
\end{eqnarray*}
Together with (\ref{Psi-C0-decomp}), we obtain (\ref{s1-C0-decomp}).

\medskip

To show (\ref{s1-Cr}) and (\ref{s1-Cr+}), let us recall {\it Fa\`a di Bruno's formula} regarding the derivatives of compositions, under the formulation due to Constantine-Savits \cite{cs-faa}. Let $f$ and $g$ be smooth functions in some domains in $\R^m$ and $\R^n$ respectively. Let $x$ be a point at which $h(x):=(f\circ g)(x)$ is well defined. Here,
$$g(x)=(g_1(x),\cdots, g_m(x)),\quad x=(x_1,\cdots,x_n),\quad y=(y_1,\cdots, y_m).$$
Let ${\bf P}=(p_1,\cdots,p_n)\in\N^n\setminus\{0\}$ and set $R=|{\bf P}|:=p_1+\cdots+p_n$, ${\bf P}!:=p_1!\cdots p_n!$.
For ${\bf l}=(L_1,\cdots,L_n)\in\N^n$ and ${\bf k}=(K_1,\cdots,K_m)\in \N^m$, set $\partial_x^{\bf l}:=\partial_{x_1}^{L_1}\cdots \partial_{x_n}^{L_n}$, and
$$  {\bf g}_{\bf l}:=(\partial_x^{\bf l}g_1(x),\ldots,\partial_x^{\bf l} g_m(x)),\quad  {\bf g}_{\bf l}^{\bf k}:=(\partial_x^{\bf l}g_1(x))^{K_1}\cdots(\partial_x^{\bf l} g_m(x))^{K_m}.$$
 For ${\bf l}=(L_1,\cdots,L_n), {\bf l}'=(L'_1,\cdots,L'_n)\in\N^n$, we say ${\bf l}\prec {\bf l}'$, if either $|\bf l|<|{\bf l}'|$, or $|{\bf l}|=|{\bf l}'|$ with $L_1< L'_1$ or $|{\bf l}|=|{\bf l}'|$ with $L_1=L'_1$, $\cdots$, $L_{k}=L'_{k}$ and $L_{k+1}< L'_{k+1}$ for some $1\leq k<n$.
\begin{thm}\label{faadibruno}\cite{cs-faa} With the above notations, for $h(x)=(f\circ g)(x)$, we have
\beq\label{partialcompo}
\partial^{\bf P}_xh(x)=\mathbf{P!}\sum_{{\bf Q}\in \N^m\atop{1\leq |{\bf Q}|\leq|{\mathbf P}|}}\partial^{\bf Q}_yf(g(x))\sum_{s=1}^{|{\bf P}|} \sum_{({\bf k},{\bf l})\in p_s({\bf P},{\bf Q})}\prod_{j=1}^s\frac{\mathbf{g}_{{\bf l}_j}^{{\bf k}_j}}{{\bf k}_j!({\bf l}_j!)^{|{\bf k}_j|}}
\eeq
where $p_s({\bf P},{\bf Q}):=\left\{({\bf k},{\bf l})\in(\N^m)^s\times(\N^n)^s:
\begin{array}{ll}
 |{\bf k}_i|>0, & 0\prec {\bf l}_1\prec\cdots\prec {\bf l}_s\\[2mm]
 \sum_{i=1}^s{\bf k}_i={\bf Q}, & \sum_{i=1}^s |{\bf k}_i|{\bf l}_i={\bf P}
\end{array}\right\}$.
\end{thm}
As a corollary of Faa di Bruno formula above, we have
\begin{lemma}\label{lem-compoCk}($C^R$-norm of composition)\cite{deW23}[Lemma 46], \cite{hoemander-geodesy}[Theorem A.8]
	For $i=1,2,3$, let $B_i$ be a compact convex domain in $\R^{n_i}$ with interioir points. Let $R\geq 1$. There exists $C_R>0$ such that if $g:B_1\rightarrow B_2$ and $f:B_2\rightarrow B_3$ are both $C^R$, then $f\circ g$ is $C^R$, and
$$
	\|f\circ g\|_{C^R}\leq C_R\left(\|f\|_{C^R}\|g\|_{C^1}^R+\|f\|_{C^1}\|g\|_{C^R}+\|f\circ g\|_{C^0}\right).
$$
\end{lemma}

\begin{lemma}\label{estim-phiCr} For the $C^\infty-$vector function $\phi=\phi(x,\xi)$ in (\ref{expression-s1}), we have, for $R\in\N$, $\|\phi(\cdot,v(\cdot))\|_{C^R}\lesssim_R \|v\|_{C^R}$.
\end{lemma}
\proof We apply Fa\`a di Bruno's formula in \rt{faadibruno} to $\partial^{\bf P}_x\phi(x,v(x))=\partial^{\bf P}_x(f\circ g)(x)$ for $f=\phi$ and $g(x)=(x,v(x))$ with $m=2n$ and $|{\bf P}|=R$. All the inequalities with ``$\lesssim$" in the proof means boundedness from above by an implicit constant depending on $R$.
In the sum
$$
\partial^{\bf P}_x\phi(x,v(x))=\mathbf{P!}\sum_{{\bf Q}\in \N^m\atop{1\leq |{\bf Q}|\leq R}}\left.\left(\partial^{\bf Q}_y\phi\right)(y)\right|_{y=(x,v(x))}\sum_{s=1}^{R} \sum_{({\bf k},{\bf l})\in p_s({\bf P},{\bf Q})}\prod_{j=1}^s\frac{{\bf g}_{{\bf l}_j}^{{\bf k}_j}}{{\bf k}_j!({\bf l}_j!)^{|{\bf k}_j|}},
$$
with $p_s({\bf P},{\bf Q})$ defined as in \rt{faadibruno}, we have
\begin{equation}\label{phi_Cr}
\left|\left.\left(\partial^{\bf Q}_y\phi\right)(y)\right|_{y=(x,v(x))}\right|\leq \|\phi\|_{C^R} .\end{equation}
Let us write ${\bf Q}\in \N^{2n}$ as
${\bf Q}=({\bf Q}',{\bf Q}'')$ with $\partial^{\bf Q}_y\phi=\partial^{{\bf Q}'}_x \partial^{{\bf Q}''}_\xi\phi$.
 If ${\bf Q}''=0$, then, according to \rrem{rem-0-deriv},
\begin{equation}\label{case_r2nul-phi}
\left|\left.\left(\partial^{\bf Q}_y\phi\right)(y)\right|_{y=(x,v(x))}\right|=\left|\left.\left( \partial^{{\bf Q}'}_x \phi\right)(y)\right|_{y=(x,v(x))}\right|\lesssim \|v\|_{C^0} .
\end{equation}

For any $({\bf k},{\bf l})\in p_s({\bf P},{\bf Q})$ with ${\bf k}=({\bf k}_j)_{1\leq j\leq s}=:({\bf k}_j',{\bf k}_j'')_{1\leq j\leq s}\in (\N^{2n})^s$, we have
\begin{equation}\label{decomp_Q-phi}
 \sum_{j=1}^s {\bf k}'_j={\bf Q}',\quad \sum_{j=1}^s {\bf k}''_j={\bf Q}'',\quad {\bf g}_{{\bf l}_j}^{{\bf k}_j}=\left(\partial_x^{{\bf l}_j}v\right)^{{\bf k}''_j}.\end{equation}
Define $l_j:=|{\bf l}_j|$, $k_j:=|{\bf k}_j|$, $k'_j:=|{\bf k}'_j|$, $k''_j:=|{\bf k}''_j|$, and decompose $R$ as
\begin{equation}\label{decomp_r-phi}R=R'+R'',\quad R':=\sum_{j=1}^s k'_j l_j  ,\quad  R'':=\sum_{j=1}^s k''_j l_j.   \end{equation}
It is easy to see that $\# p_s({\bf P},{\bf Q})\lesssim 1$, and
$$
\left| \prod_{j=1}^s\frac{\mathbf{g}_{{\bf l}_j}^{{\bf k}_j}}{{\bf k}_j!({\bf l}_j!)^{|{\bf k}_j|}}\right|\lesssim  \prod_{j=1}^s\left| \mathbf{g}_{{\bf l}_j}^{{\bf k}_j}\right|
\lesssim \prod_{j=1}^s \|v\|_{C^{l_j}}^{k''_j},\qquad ({\bf k},{\bf l})\in p_s({\bf P},{\bf Q}).$$
If $R''>0$, then, applying \rl{interpol-cs}  to $\|v\|_{C^{l_j}}$ with $a=0$, $b=R''$, we have
\beq\label{v-prod-cr}
\prod_{j=1}^s \|v\|_{C^{l_j}}^{k''_j}\lesssim \prod_{j=1}^s \|v\|_{C^0}^{k''_j(1-\frac{l_j}{R''})} \prod_{j=1}^s\|v\|_{C^{R''}}^{k''_j\frac{l_j}{R''}}\leq\prod_{j=1}^s \|v\|_{C^{R''}}^{k''_j\frac{l_j}{R''}}= \|v\|_{C^{R''}},
\eeq
since $\|v\|_{C^0}$ is sufficiently small.
Hence, combining with (\ref{phi_Cr}),  we have
$$
\left| \left.\left(\partial^{\bf Q}_y\phi\right)(y)\right|_{y=(x,v(x))}  \prod_{j=1}^s\frac{\mathbf{g}_{{\bf l}_j}^{{\bf k}_j}}{{\bf k}_j!({\bf l}_j!)^{|{\bf k}_j|}}\right| \lesssim \|v\|_{C^{R''}}.
$$
If $R''=0$, which means that ${\bf Q}''=0$, then, according to (\ref{case_r2nul-phi}) -- (\ref{decomp_r-phi}), we have
$$\left| \left.\left(\partial^{\bf Q}_y\phi\right)(y)\right|_{y=(x,v(x))} \prod_{j=1}^s\frac{\mathbf{g}_{{\bf l}_j}^{{\bf k}_j}}{{\bf k}_j!({\bf l}_j!)^{|{\bf k}_j|}}\right|\lesssim  \|v\|_{C^{0}}.$$
Therefore, for any ${\bf Q}\in \N^{2n}$ with $1\leq |{\bf Q}|\leq R$, and any $({\bf k},{\bf l})\in p_s({\bf P},{\bf Q})$, $1\leq s \leq R$,
$$\left| \left.\left(\partial^{\bf Q}_y\phi\right)(y)\right|_{y=(x,v(x))} \prod_{j=1}^s\frac{\mathbf{g}_{{\bf l}_j}^{{\bf k}_j}}{{\bf k}_j!({\bf l}_j!)^{|{\bf k}_j|}}\right|\lesssim \|v\|_{C^{R''}} \leq  \|v\|_{C^{R}} .$$
The lemma is shown.\qed

%
%
%

\smallskip

\begin{lemma}\label{lemma_Psi-r}
For $w, v\in \Gamma^\infty(M,TM)$ with $\|w\|_{C^1}$ and $\|v\|_{C^1}$ sufficiently small, we have, for $R\in \N^*$,
$$\|\Psi_{w,v}\|_{C^R}\lesssim_R \|w\|_{C^R}+\|w\|_{C^1}\|v\|_{C^R},\qquad \|\Psi_{w,v}\|_{C^R}\lesssim_R \|w\|_{C^2}\|v\|_{C^R} + \|w\|_{C^{R+1}}  \|v\|_{C^0} .$$
\end{lemma}
\proof For $R\in \N^*$, by \rl{lem-compoCk}, we have
$$\|\Psi_{w,v}\|_{C^R}\lesssim_R  \|w\|_{C^R} \left( 1+\|v\|_{C^1} \right)^R +  \|w\|_{C^1} \|v\|_{C^R} + \|w\|_{C^0}
\lesssim_R \|w\|_{C^R} + \|w\|_{C^1} \|v\|_{C^R}.$$
On the other hand, in view of the expression (\ref{Psi-integral}) of $\Psi_{w,v}$ and combining with \rl{lem-compoCk}, we have
\begin{eqnarray*}
\|\Psi_{w,v}\|_{C^R}&\lesssim_R& \|w^{(1)}\|_{C^0}\|v\|_{C^R}  \\
& &   + \,  \left(\|w^{(1)}\|_{C^R} \left( 1+\|v\|_{C^1} \right)^R+  \|w^{(1)}\|_{C^1} \|v\|_{C^R} + \|w^{(1)}\|_{C^0}\right) \|v\|_{C^0}    \\
   &\lesssim_R&  \|w^{(1)}\|_{C^0}\|v\|_{C^R} + \|w^{(1)}\|_{C^R}  \|v\|_{C^0}  + \|w^{(1)}\|_{C^1} \|v\|_{C^R}  \|v\|_{C^0} \\
   &\lesssim_R&  \|w\|_{C^1}\|v\|_{C^R} + \|w\|_{C^{R+1}}  \|v\|_{C^0}  + \|w\|_{C^2} \|v\|_{C^R}  \|v\|_{C^0}\\
   &\lesssim_R&  \|w\|_{C^2}\|v\|_{C^R} + \|w\|_{C^{R+1}}  \|v\|_{C^0} .\qed
\end{eqnarray*}

\begin{lemma}\label{lemma_Upsilon-r}
 For $w, v\in \Gamma^\infty(M,TM)$ with $\|w\|_{C^0}$, $\|v\|_{C^0}$ sufficiently small, we have
\begin{equation}\label{esti-Upsilon}
\|\Upsilon_{w,v}\|_{C^R}\lesssim_R \|w\|_{C^R} \|v\|_{C^{0}}+\|w\|_{C^0} \|v\|_{C^{R}},\quad \forall \ R\in \N^*.\end{equation}
\end{lemma}

\proof It is sufficient to show that, for any ${\bf P}\in\N^n$ with $|{\bf P}|=R$,
\begin{equation}\label{esti-Upsilon-partial}
 \left\|\partial^{\bf P}_x\Upsilon_{w,v}\right\|_{C^{0}} \lesssim_R \|w\|_{C^R} \|v\|_{C^{0}}+\|w\|_{C^0} \|v\|_{C^{R}}.
\end{equation}
As in the proof of Lemma \ref{estim-phiCr}, the inequalities with ``$\lesssim$" in the proof means boundedness from above by an implicit constant depending on $R$.

With the expression of $\Upsilon_{w,v}$ in (\ref{expression-s1}), we apply Fa\`a di Bruno's formula to $\partial^{\bf P}_x\Upsilon_{w,v}(x)=\partial^{\bf P}_x(f\circ g)(x)$ for $f=\varrho$ and $g(x)=(x,v(x),w(x))$ with $m=3n$. In the sum
$$
\partial^{\bf P}_x\Upsilon_{w,v}(x)=\mathbf{P!}\sum_{{\bf Q}\in \N^m\atop{1\leq |{\bf Q}|\leq R}}\left.\left(\partial^{\bf Q}_y\varrho\right)(y)\right|_{y=(x,v(x),w(x))}\sum_{s=1}^{r} \sum_{({\bf k},{\bf l})\in p_s({\bf P},{\bf Q})}\prod_{j=1}^s\frac{\mathbf{g}_{{\bf l}_j}^{{\bf k}_j}}{{\bf k}_j!({\bf l}_j!)^{|{\bf k}_j|}},
$$
with $p_s({\bf P},{\bf Q})$ defined as in \rt{faadibruno}, we have
\begin{equation}\label{rho_Cr}
\left|\left.\left(\partial^{\bf Q}_y\varrho\right)(y)\right|_{y=(x,v(x),w(x))}\right|\leq \|\varrho\|_{C^R} .\end{equation}
Let us write ${\bf Q}\in \N^{3n}$ as
${\bf Q}=({\bf Q}',{\bf Q}'',{\bf Q}''')$ with $\partial^{\bf Q}_y\varrho=\partial^{{\bf Q}'}_x \partial^{{\bf Q}''}_\xi\partial^{{\bf Q}'''}_\eta\varrho$.
Recalling \rrem{rem-0-deriv}, if ${\bf Q}''=0$, then
\begin{equation}\label{case_r2nul}
\left|\left.\left(\partial^{\bf Q}_y\varrho\right)(y)\right|_{y=(x,v(x),w(x))}\right|=\left|\left.\left( \partial^{{\bf Q}'}_x \partial^{{\bf Q}'''}_\eta\varrho\right)(y)\right|_{y=(x,v(x),w(x))}\right|\lesssim \|v\|_{C^0} .
\end{equation}
If ${\bf Q}'''=0$, then
\begin{equation}\label{case_r3nul}
\left|\left.\left(\partial^{\bf Q}_y\varrho\right)(y)\right|_{y=(x,v(x),w(x))}\right|=\left|\left.\left( \partial^{{\bf Q}'}_x \partial^{{\bf Q}''}_\xi\varrho\right)(y)\right|_{y=(x,v(x),w(x))}\right|\lesssim \|w\|_{C^0} .\end{equation}
If ${\bf Q}''={\bf Q}'''=0$, then
\begin{equation}\label{case_r2r3nul}
\left|\left.\left(\partial^{\bf Q}_y\varrho\right)(y)\right|_{y=(x,v(x),w(x))}\right|=\left|\left.\left( \partial^{{\bf Q}'}_x\varrho\right)(y)\right|_{y=(x,v(x),w(x))}\right|\lesssim  \|v\|_{C^0} \|w\|_{C^0} .\end{equation}

For $({\bf k},{\bf l})\in p_s({\bf P},{\bf Q})$ with ${\bf k}=({\bf k}_j)_{1\leq j\leq s}=:({\bf k}_j',{\bf k}_j'',{\bf k}_j''')_{1\leq j\leq s}\in (\N^{3n})^s$, we have
\begin{equation*}\label{decomp_Q}
\sum_{j=1}^s {\bf k}'_j={\bf Q}',\quad \sum_{j=1}^s {\bf k}''_j={\bf Q}'',\quad \sum_{j=1}^s {\bf k}'''_j={\bf Q}''',\qquad {\bf g}_{{\bf l}_j}^{{\bf k}_j}=\left(\partial_x^{{\bf l}_j}v\right)^{{\bf k}''_j}\left(\partial_x^{{\bf l}_j}w\right)^{{\bf k}'''_j}.\end{equation*}
Define $l_j:=|{\bf l}_j|$, $k_j:=|{\bf k}_j|$, $k'_j:=|{\bf k}'_j|$, $k''_j:=|{\bf k}''_j|$, $k'''_j:=|{\bf k}'''_j|$, and decompose $R$ as
\begin{equation*}\label{decomp_r}
R=R'+R''+R''',\quad R':=\sum_{j=1}^s k'_j l_j  ,\quad  R'':=\sum_{j=1}^s k''_j l_j,\quad R''':=\sum_{j=1}^s k'''_j l_j.   \end{equation*}
It is easy to see that $\# p_s({\bf P},{\bf Q})\lesssim 1$, and
$$
\left| \prod_{j=1}^s\frac{\mathbf{g}_{{\bf l}_j}^{{\bf k}_j}}{{\bf k}_j!({\bf l}_j!)^{|{\bf k}_j|}}\right|\lesssim  \prod_{j=1}^s\left| \mathbf{g}_{{\bf l}_j}^{{\bf k}_j}\right|
\lesssim \prod_{j=1}^s \|v\|_{C^{l_j}}^{k''_j}\cdot \prod_{j=1}^s \|w\|_{C^{l_j}}^{k'''_j},\qquad ({\bf k},{\bf l})\in p_s({\bf P},{\bf Q}).$$
Since $\|v\|_{C^0}$ and $\|w\|_{C^0}$ are sufficiently small, similar to (\ref{v-prod-cr}), we have
$$\prod_{j=1}^s \|v\|_{C^{l_j}}^{k''_j}\lesssim  \|v\|_{C^{R''}}  \;\  {\rm if}  \;\  R''>0,\qquad  \prod_{j=1}^s \|w\|_{C^{l_j}}^{k'''_j}\lesssim  \|w\|_{C^{R'''}}   \;\  {\rm if}  \;\  R'''>0. $$
Noting that $R''=0$ (resp. $R'''=0$) means that ${\bf Q}''=0$ (resp. ${\bf Q}'''=0$), we have, in view of (\ref{rho_Cr}) -- (\ref{case_r2r3nul}),  for any ${\bf Q}\in \N^{3n}$ with $1\leq |{\bf Q}|\leq R$, and any $({\bf k},{\bf l})\in p_s({\bf P},{\bf Q})$, $1\leq s \leq R$,
\begin{eqnarray*}
\left| \left.\left(\partial^{\bf Q}_y\varrho\right)(y)\right|_{y=(x,v(x),w(x))} \prod_{j=1}^s\frac{\mathbf{g}_{{\bf l}_j}^{{\bf k}_j}}{{\bf k}_j!({\bf l}_j!)^{|{\bf k}_j|}}\right|&\lesssim& \|v\|_{C^{R''}}  \|w\|_{C^{R'''}}\\
&\leq& \|v\|_{C^{R''}}  \|w\|_{C^{R-R''}}\\ &\lesssim& \|w\|_{C^{0}} \|v\|_{C^{R}}+\|w\|_{C^{R}} \|v\|_{C^{0}},
\end{eqnarray*}
where we obtain the above last inequality through interpolation and the concavity of logarithm: for any $0\leq k\leq R$,
\begin{eqnarray*}
\|v\|_{C^{k}}  \|w\|_{C^{R-k}}&\lesssim& \|v\|^{1-\frac{k}{R}}_{C^{0}}\|v\|^{\frac{k}{R}}_{C^{R}} \|w\|_{C^{0}}^{\frac{k}{R}}\|w\|^{1-\frac{k}{R}}_{C^{R}}\\
&=& \exp\left\{\frac{k}{R}\ln(\|w\|_{C^{0}} \|v\|_{C^{R}})+\frac{R-k}{R}\ln(\|w\|_{C^{R}} \|v\|_{C^{0}}) \right\}\\
&\leq&\exp\left\{\ln\left(\frac{k}{R}\|w\|_{C^{0}} \|v\|_{C^{R}}+\frac{R-k}{R}\|w\|_{C^{R}} \|v\|_{C^{0}}\right) \right\}\\
&\leq&\|w\|_{C^{0}} \|v\|_{C^{R}}+\|w\|_{C^{R}} \|v\|_{C^{0}}.
\end{eqnarray*}
Then the inequality (\ref{esti-Upsilon-partial}) is shown.\qed

\smallskip

Combining Lemma \ref{lemma_Psi-r} and \ref{lemma_Upsilon-r}, we obtain (\ref{s1-Cr}) and (\ref{s1-Cr+}).

\section{Proof of Proposition \ref{propMoser-anal}.}\label{app_proof}
Let $W=W_i$ be a trivializing coordinate patch as in Section \ref{sec_Grauert}. 
Given $q\in M_{r'}\cap W$, with $z=z(q)$, let $\eta:=\tilde w(z)$, $\xi:=\tilde v(z)$.
If $|\xi|_{\kappa}$ is small enough, then, recalling (\ref{complex-exp}), $\Psi(z,\eta)$ defines the coordinates of a point in $M_{r}$ and we apply the composition of flows~:
$$\zeta=P(z,\xi, \tilde w(\Psi(z,\xi)))
	= \xi+\tilde w(\Psi(z,\xi))+ \varrho(z,\xi, \tilde w(\Psi(z,\xi))). $$
Let us set
\begin{equation}\label{s1}
s_1(w,v)(z) := \zeta-\tilde w(z)-\xi
= \left(\tilde w(\Psi(z,\xi))-\tilde w(z)\right)+ \varrho(z,\xi, \tilde w(\Psi(z,\xi))).
\end{equation}
All these quantities are well defined if $\|w\|_{0,r}$ and $\|v\|_{0,r}$ are sufficiently small. 
Let $v_1$, $v_2$ be two holomorphic small enough vector fields on $W\cap{M}_{r'}$ and let us set
$$\om_1(z):= \tilde w(\Psi(z,\tilde v_1(z))),\quad \om_2(z):= \tilde w(\Psi(z,\tilde v_2(z))).$$
We have that
\begin{eqnarray}
& &	\sup_{z=z(q)\in\Delta_1^n\atop{q\in M_{r'}\cap W}} |\om_1(z)-\om_2(z)|_{\ka}\label{omom'}\\
&\leq& \sup_{\tilde q\in M_{r_*}\cap W}\sup_{\zeta\in \C^n\atop{|\zeta|\leq 1}} |D_{z} \tilde w(z(\tilde q))\zeta|_{\ka}
\sup_{z=z(q)\in  \Delta_1^n\atop{q\in M_{r'}\cap W}}| \Psi(z,\tilde v_1(z))- \Psi(z,\tilde v_2(z))|_{\ka}\nonumber\\
&\lesssim&\sup_{\tilde q\in M_{r_*}\cap W}\sup_{\zeta\in \mathbb{C}^n\atop{|\zeta|\leq 1}}|D_{z} \tilde w(z(\tilde q))\zeta|_{\ka}
\sup_{z=z(q)\in  \Delta_1^n\atop{q\in M_{r'}\cap W}}|\tilde v_1(z)-\tilde v_2(z)|_{\ka}.\nonumber
\end{eqnarray}
On the other hand, we have
\begin{eqnarray*}
(s_1(w,v_1)-s_1(w,v_2))(z) &=& (\tilde w(\Psi(z,\tilde v_1(z)))-\tilde w(z)) - (\tilde w(\Psi(z,\tilde v_2(z)))-\tilde w(z))  \\
   & & + \, \varrho(z,\tilde v_1(z), \tilde w(\Psi(z,\tilde v_1(z))))- \varrho(z,\tilde v_2(z), \tilde w(\Psi(z,\tilde v_2(z)))).
\end{eqnarray*}
	According to \re{s1}, we have
\begin{eqnarray*}
& & \|s_1(w,v_1)-s_1(w,v_2)\|_{C^0,r'}\\
&\leq& \|\om_1-\om_2\|_{C^0,r'}+\|\varrho(z,\tilde v_1(z), \om_1(z))-\varrho(z,\tilde v_2(z), \om_2(z))\|_{C^0,r'}\\
&\leq& \|\om_1-\om_2\|_{C^0,r'}\\
& &  + \, \sup_{t\in[0,1]}\|\partial_\xi \varrho(z,t\tilde v_1(z)+(1-t)\tilde v_2(z), t\om_1(z)+(1-t)\om_2(z))\|_{C^0,r}\|v_1-v_2\|_{C^0,r'} \\
& &+ \, \sup_{t\in[0,1]}\|\partial_\eta\varrho(z,t\tilde v_1(z)+(1-t)\tilde v_2(z), t\om_1(z)+(1-t)\om_2(z))\|_{C^0,r}\|\om_1-\om_2\|_{C^0,r'}.
\end{eqnarray*}
As $\|v_1\|_{C^0,r}$, $\|v_2\|_{C^0,r}$, $\|\om_1\|_{C^0,r}$, $\|\om_2\|_{C^0,r}$ are uniformly bounded, it is deduced from \re{omom'} that
$\|s_1(w,v_1)-s_1(w,v_2)\|_{C^0,r'}\lesssim \|w\|_{C^1,r}\|v_1-v_2\|_{C^0,r'}$, which completes the proof.
\qed

\section{Interpolation inequalities}\label{app_interpolation}

\begin{lemma}\label{interpol-cs} (Interpolation of $C^r-$norms, \cite{hoemander-geodesy}) For $0\leq a\leq b<\infty$, $0<\lambda <1$,
$$\|u\|_{C^{\lambda a +(1-\lambda)b}}\lesssim_{\lambda, a, b} \|u\|^\lambda_{C^{a}} \|u\|^{1-\lambda}_{C^{b}}. $$
\end{lemma}

\begin{lemma} (Interpolation of Sobolev norms)\label{interpol-hs}  For $0\leq a\leq b<\infty$, $0<\lambda <1$,
$$\|u\|_{\CH^{\lambda a +(1-\lambda)b}}\leq \|u\|^\lambda_{\CH^{a}}  \|u\|^{1-\lambda}_{\CH^{b}},\qquad u\in\CH^{b}.$$
\end{lemma}
\proof For $u=\sum_{j\in\N} u_j {\bf e}_j$, it is sufficient to show that
$$\sum_{j\in\N} (1+\tilde\lambda_j)^{2(\lambda a +(1-\lambda)b)} |u_j|^2 \leq \left(\sum_{j\in\N} (1+\tilde\lambda_j)^{2 a} |u_j|^2 \right)^{\lambda} \left(\sum_{j\in\N} (1+\tilde\lambda_j)^{2 b} |u_j|^2 \right)^{1-\lambda}.$$
Applying H\"older's inequality $\|fg\|_{\ell^1}\leq \|f\|_{\ell^p} \|g\|_{\ell^q}$, $\frac{1}{p}+\frac{1}{q}=1$, with
$$f_j= (1+\tilde\lambda_j)^{2\lambda a} |u_j|^{2\lambda},\quad g_j= (1+\tilde\lambda_j)^{2(1-\lambda) b} |u_j|^{2(1-\lambda)},\quad p=\frac{1}{\lambda},\quad q=\frac{1}{1-\lambda},$$
the above inequality is shown.\qed

With a similar proof as Lemma \ref{interpol-hs}, we have

\begin{lemma} (Interpolation of Hardy norms)\label{interpol-hardy}  For $0\leq r'\leq r<r_*$, $0<\lambda <1$,
$$\|u\|_{\lambda r' +(1-\lambda)r}\leq \|u\|^\lambda_{r'} \|u\|^{1-\lambda}_{r},\quad u\in H_r^2.$$
\end{lemma}

\end{document}